\documentclass[11pt,reqno]{amsart}
\usepackage[foot]{amsaddr}

\usepackage{times}
\usepackage[T1]{fontenc}
\usepackage{mathrsfs}
\usepackage{latexsym}
\usepackage[dvips]{graphics}
\usepackage{epsfig}
\usepackage{amsmath,amsfonts,amsthm,amssymb,amscd}
\input amssym.def
\input amssym.tex
\usepackage{color}

\usepackage{enumitem}
\usepackage{xcolor,color}
\usepackage{mathtools,stackrel,bm}
\usepackage{float,caption,subcaption,grffile,wrapfig}
\usepackage{array,multirow,tabularx,tabulary,booktabs,diagbox}

\usepackage[colorlinks,citecolor=blue,urlcolor=blue,linkcolor=black]{hyperref}%% uncomment this for coloring bibliography citations and linked URLs
\usepackage{xurl} % add URL line breaks if available
\DeclareRobustCommand{\SkipTocEntry}[5]{}

\usepackage{adjustbox}
\allowdisplaybreaks[4]

\usepackage{cleveref}
\numberwithin{equation}{section}
\newtheorem{theorem}{Theorem}[section]

\newtheorem{corollary}[theorem]{Corollary}
\newtheorem{definition}[theorem]{Definition}

\newtheorem{lemma}[theorem]{Lemma}
\newtheorem{proposition}[theorem]{Proposition}

\newtheorem*{claim}{Claim}
\newenvironment{remark}{{\bfseries Remark.\quad}}{}
\renewenvironment{proof}[1][Proof]{\leavevmode\newline\par\noindent
   {\bfseries\itshape\underline{#1.}\quad}
   }{
   \hfill\rule{2mm}{2mm}
}

\crefname{theorem}{Theorem}{Theorems}
\crefname{section}{Section}{Sections}
\crefname{figure}{Figure}{Figures}
\crefname{table}{Table}{Tables}
\crefname{lemma}{Lemma}{Lemmas}
\crefname{assumption}{Assumption}{Assumptions}
\crefname{corollary}{Corollary}{Corollaries}
\crefname{definition}{Definition}{Definitions}
\crefname{example}{Example}{Examples}
\crefname{proposition}{Proposition}{Propositions}
\crefname{conjecture}{Conjecture}{Conjectures}
\crefname{appendix}{Appendix}{Appendices}
\crefname{subappendix}{Appendix}{Appendices}

\newlist{theoremenum}{enumerate}{1} 
\setlist[theoremenum]{label=(\alph*), ref=\thetheorem\alph*,labelindent=1em}
\crefalias{theoremenumi}{theorem} 
\newlist{propositionenum}{enumerate}{1} 
\setlist[propositionenum]{label=(\alph*), ref=\theproposition\alph*,labelindent=1em}
\crefalias{propositionenumi}{proposition} 
\newlist{corollaryenum}{enumerate}{1} 
\setlist[corollaryenum]{label=(\alph*), ref=\thecorollary\alph*,labelindent=1em}
\crefalias{corollaryenumi}{corollary} 
\newlist{lemmaenum}{enumerate}{1} 
\setlist[lemmaenum]{label=(\alph*), ref=\thelemma\alph*,labelindent=1em}
\crefalias{lemmaenumi}{lemma} 
\newlist{casesenum}{enumerate}{1} 
\setlist[casesenum]{label=\arabic*.,ref=\arabic*,labelindent=1em}%
\crefname{casesenumi}{case}{cases} 
\Crefname{casesenumi}{Case}{Cases}

\newcommand*{\dif}{\mathop{}\!\mathrm{d}}
\newcommand*{\comp}{\mathbin{\vcenter{\hbox{\scalebox{.6}{$\circ$}}}}}

\begin{document}

\title[Wasserstein-$P$ Bounds in CLT Under Local Dependence]{Wasserstein-$\boldsymbol P$ Bounds in the Central Limit Theorem Under Local Dependence}

\author{Tianle Liu} \email{tianleliu@fas.harvard.edu}
\author{Morgane Austern} \email{maustern@fas.harvard.edu} 
\address{Department of Statistics, Harvard University, Cambridge, MA 02138}

\subjclass[2020]{60F05.}

\keywords{Normal approximation; Wasserstein distance; dependency neighborhood; U-statistics; Stein's method}

% \date{\today}

\begin{abstract}
  The central limit theorem (CLT) is one of the most fundamental results in probability; and establishing its rate of convergence has been a key question since the 1940s. For independent random variables, a series of recent works established optimal error bounds under the Wasserstein-$p$ distance (with $p\ge 1$). In this paper, we extend those results to locally dependent random variables, which include $m$-dependent random fields and U-statistics. Under conditions on the moments and the dependency neighborhoods, we derive optimal rates in the CLT for the Wasserstein-$p$ distance. Our proofs rely on approximating the empirical average of dependent observations by the empirical average of \emph{i.i.d.} random variables. To do so, we expand the Stein equation to arbitrary orders by adapting the Stein's dependency neighborhood method. Finally we illustrate the applicability of our results by obtaining efficient tail bounds. 
\end{abstract}

\maketitle

\tableofcontents

\section{Introduction}

The central limit theorem (CLT) is one of the most fundamental theorems in probability theory. Initially formulated for independent and identically distributed random variables, it has since then been generalized to triangular arrays \cite{feller1945fundamental}, martingales \cite{levy1935proprietes}, U-statistics \cite{hoeffding1948class}, locally dependent random variables \cite{hoeffding1948central,heinrich1982method,petrovskaya1983central}, and mixing random fields \cite{rosenblatt1956central,bolthausen1982central}. Let $(I_n)$ be an increasing sequence of subsets $I_1\subseteq I_2\subseteq\cdots \subseteq I$, whose sizes increase to infinity $|I_n|\rightarrow \infty$. Set $(X_i)_{i\in I}$ to be (dependent) centered random variables. Under certain conditions on the moments of $(X_i)$ and on its dependence structure, the CLT states that the scaled sum is asymptotically normal, i.e.,
$$W_n:=\sigma_n^{-1}\sum_{i\in I_n} X_i\xrightarrow{d}\mathcal{N}(0,1),$$ 
where we write $\sigma_n^2=\operatorname{Var}\bigl(\sum_{i\in I_n}X_i\bigr)$. Starting with the work of Berry and Esseen in 1940s, there is a long history of quantifying how far $W_n$ is from being normally distributed. One of the most important metrics to do so is the Wasserstein-$p$ distance originated in optimal transport theory \cite{villani2009optimal}. For two probability measures $\nu$ and $\mu$ over the real line $\mathbb{R}$, we denote by $\Gamma(\nu,\mu)$ the set of all couplings of $\nu$ and $\mu$, and the Wasserstein-$p$ distance between $\nu$ and $\mu$ is defined as 
\begin{equation*}
\mathcal{W}_p(\nu,\mu):=\inf_{\gamma\in \Gamma(\nu,\mu)}\Bigl(\mathbb{E}_{(X,Y)\sim \gamma}\left[\lvert X-Y \rvert^p\right]\Bigr)^{1/p}.
\end{equation*}
When the observations $(X_i)$ are independent, \cite{agnew1957estimates,esseen1958mean} established that for $p=1$ the convergence rate for the CLT is $\mathcal{O}\bigl(\lvert I_{n} \rvert^{-\frac{1}{2}}\bigr)$. Extending such results to $p>1$ remained for a while an open question. The first bounds for $p> 1$ obtained by \cite{bartfai1970entfernung} dating back to the 1970s were sub-optimal in terms of the sample size $|I_n|$ as they decrease at a slower rate of $\mathcal O\bigl(|I_n|^{-\frac{1}{2}+\frac{1}{p}}\bigr)$. \cite{rio2009upper} obtained that, for $1\leq p\leq 2$, the Wasserstein distance converges at the optimal rate $\mathcal{O}\bigl(\lvert I_{n} \rvert^{-\frac{1}{2}}\bigr)$ under some additional necessary moment conditions, and they conjectured that such a rate would be extendable to arbitrary $p\ge 1$. This was recently proven to be true by \cite{bobkov2018berry,bonis2020stein} using a series of methods including the Edgeworth expansion and the exchangeable pair method. They showed that if  $\max_i \|X_i\|_{p+2}<\infty$ and if $\operatorname{Var}(X_1)=\operatorname{Var}(X_i)=1$, then there is a constant $K_p<\infty$ such that
\begin{equation*}
\mathcal{W}_p\bigl(\mathcal{L}(W_n),\mathcal{N}(0,1)\bigr)\le \frac{K_p\lVert X_1\rVert_{p+2}^{1+2/p}}{\sqrt{\lvert I_{n} \rvert}},
\end{equation*}
where $\mathcal{L}(\,\cdot\,)$ designates the distribution of the given random variable. It is however crucial to note that these rates were obtained under the key assumption of independence of the $(X_i)$. In this paper, we aim to generalize this beyond the assumption of independence which is restrictive for many applications. 

An important class of dependent observations $(X_i)$ are locally dependent random variables. Intuitively, we say that $(X_i)$ are locally dependent if for every finite group of random variables $(X_i)_{i\in J}$, where $J\subset I$, there exists a subset $N(J)\subset I$ such that $(X_i)_{i\in J}$ is independent from $(X_i)_{i\in I\setminus N(J)}$. The subset $N(J)$ is often called the dependency neighborhood of $J$. Examples of such random variables include $m$-dependent random fields, U-statistics, and subgraph count statistics in the \text{Erd\H{o}s--R\'{e}nyi} random graphs. Under general conditions on the sizes of the dependency neighborhoods the central limit theorem is known to hold and its rate of convergence in Wasserstein-$1$ distance was established by \cite{baldi1989normal,barbour1989central}. This was extended to Wasserstein-$2$ bounds by \cite{fang2019wasserstein} by relating it to Zolotarev's metrics and cleverly exploiting Stein's method. Drawing inspiration from \cite{bonis2020stein}, sub-optimal rates were also achieved in \cite{fang2022p} for arbitrary $p\ge 1$ under more technical conditions. Nevertheless, an optimal rate bound for general Wasserstein-$p$ distances ($p\geq 1$) remains unknown. This is the gap that we fill in this paper. We consider locally dependent (not necessarily identically distributed) random variables $(X_i)$, and consider the empirical average $W_n:=\sigma_n^{-1}\sum_{i\in I_n}X_i$ where $\sigma_n^2:=\sum_{i\in I_n}X_i$. For all $p\ge 1$ we obtain bounds for the $\mathcal{W}_p$ distance $\mathcal{W}_p(\mathcal{L}(W_n),\mathcal{N}(0,1))$. We do so under the assumption that the variances $(\sigma_n)$ are nondegenerate, and under moment conditions and on the sizes of dependency neighborhoods. Notably if the size of the dependency neighborhoods is uniformly bounded we obtain bounds that decrease at the optimal rate (see \cref{THM:LOCALWP2})
\begin{equation*}
\mathcal{W}_p(\mathcal{L}(W_n),\mathcal{N}(0,1))=\mathcal O\Bigl(\frac{1}{\sqrt{|I_n|}}\Bigr).
\end{equation*}
We further generalize our results to triangular arrays where the random variables $\bigl(X^{\scalebox{0.6}{$(n)$}}_i\bigr)$ are allowed to change with $n$. Finally, we demonstrate how those bounds can be exploited to obtain non-uniform Berry--Esseen type bounds that have polynomial decay. 

The key idea of our proofs is to approximate the empirical average $W_n$ by an empirical average $V_n$ of \textit{i.i.d.} random variables for which Wasserstein's bounds are already known. To do this we establish an Edgeworth-type expansion of the Stein equation in terms of the cumulants of the $W_n$. Indeed, in \cref{THM:BARBOURLIKE} we prove that if $h$ is a function smooth enough (made precise later) and $Z\sim\mathcal{N}$ is a standard random variable then
\begin{equation}\label{eq:demo2}
\begin{aligned}
  &\mathbb{E}[ h(W_{n})]-\mathbb{E}[h(Z)]=\mathbb{E} [f_h'(W_{n})-W_{n}f_h(W_{n})]\\
  = &\sum_{(r,s_{1:r})\in \Gamma (\lceil p\rceil-1)}(-1)^{r}\prod_{j=1}^{r}\frac{\kappa _{s _{j}+2}(W_{n})}{(s _{j}+1)!}\mathcal{N}\ \Bigl[\prod_{j=1}^{r}(\partial ^{s _{j}+1}\Theta)\ h\Bigr] +\text{ Remainders},
\end{aligned}
\end{equation} 
where $f_h$ is the solution of the Stein equation \cref{eq:stein} and where $(\kappa _{j}(W_{n}))$ designates the cumulants of $W_n$ (the other notations will be made explicit in the next few sections). This generalizes a similar well-known result for \emph{i.i.d.} observations established in \cite{barbour1986asymptotic}. To guarantee that our choice of $V_n$ is a good approximation of $W_n$ we utilize this expansion and exploit the Hamburger moment problem to choose $V_n$ to be such that its first $\lceil p \rceil+1$ cumulants match the ones of $W_n$.

\subsection{Related Literature}
\cite{agnew1957estimates,esseen1958mean} established that the convergence rate in the central limit theorem is $\mathcal{O}\bigl(\lvert I_{n} \rvert^{-\frac{1}{2}}\bigr)$ in terms of the Wasserstein-$1$ distance. Since then it has been tightened and generalized to dependent observations. 
Notably, the Stein's method offers a series of powerful techniques for obtaining Wasserstein-$1$ bounds in the dependence setting. See \cite{ross2011fundamentals} for a survey of those methods. \cite{baldi1989normal,barbour1989central} obtained Wasserstein-$1$ bounds under local dependence conditions.
%  and \cite{sunklodas2007normal} showed a similar bound for strongly mixing sequences, which is $\mathcal{O}\left(1/\sqrt{\lvert I_{n} \rvert}\right)$ when the mixing coefficients decay fast enough.

\cite{bartfai1970entfernung} proposed a rate of $\mathcal O\bigl(|I_n|^{-\frac{1}{2}+\frac{1}{p}}\bigr)$ for the Wasserstein-$p$ distance under the hypothesis that the random variables have finite exponential moments. \cite{sakhanenko1985estimates} obtained a similar rate but only required the existence of $p$-th moments. \cite{rio1998distances,rio2009upper} showed that in order to obtain a convergence rate of $\mathcal{O}\bigl(\lvert I_{n} \rvert^{-\frac{1}{2}}\bigr)$, it is necessary to require finite ($p$+$2$)-th moments of the random variables. They also obtained the optimal rate for $1\leq p\leq 2$ and conjectured that a similar rate should be valid for any arbitrary $p> 2$. This conjecture was demonstrated to be true by \cite{bobkov2018berry,bonis2020stein}. Those two papers took different approaches. \cite{bobkov2018berry} used an Edgeworth expansion argument. \cite{bonis2020stein}, on the other hand, used the Ornstein-Uhlenbeck interpolation combined with a Stein exchangeable pair argument and their methods further applied to multivariate settings. Previous to that, \cite{ledoux2015stein} had already obtained the optimal rate for the Wasserstein-$p$ distance using the Ornstein-Uhlenbeck interpolation but needed significantly stronger assumptions on the distribution of the random variables by requiring the existence of a Stein kernel. Moreover, for the special case $p=2$, the celebrated HWI inequality \cite{otto2000generalization} and Talagrand quadratic transport inequality \cite{talagrand1996transportation} can help obtain Wasserstein-$2$ bounds by relating it to the Kullback-Leibler divergence. 

Contrary to the independent case, much less is known for the general Wasserstein-$p$ distance for dependent data. \cite{fang2019wasserstein} adapted the Stein's method to obtain Wasserstein-$2$ bounds for locally dependent variables. \cite{fang2022p} modified the approach of \cite{bonis2020stein} and obtained a sub-optimal rate $\mathcal{O}\bigl(\lvert I_{n} \rvert^{-\frac{1}{2}}\log\, \lvert I_{n} \rvert\bigr)$ for the Wasserstein-$p$ distance under local dependence. Our results propose significant extensions to both of those results by generalizing the optimal rate to arbitrary $p\geq 1$.
%  and mixing random fields. When the mixing coefficients decrease fast enough we obtain that the Wasserstein distance decreases at the optimal rate $\mathcal{O}\left(1/\sqrt{\lvert I_{n} \rvert}\right)$. Otherwise, the rate of convergence is slower and depends on the mixing coefficients.% However, contrary to \cite{fang2019wasserstein} we obtain a bound that is valid for any $p\ge 1$. Moreover, the method in \cite{fang2022p} lead to sub-optimal rates of convergence and was also not generalized to mixing processes.

Our proofs also rely on the Stein's method and a result of \cite{rio2009upper} that allows to upper the Wasserstein-$p$ distance by an integral probability metric \cite{zolotarev1983probability}. As those metrics are defined as the supremum of expected differences over a certain class of functions, the Stein's method lends itself nicely to this problem. %Notably it allows, for a large class of function $f$, to obtain local developments of $\mathbb{E}(W_nf(W_n))$ in terms of the cumulants of $W_n$ and the derivatives of $f$. This 
The Stein's method was first introduced in \cite{stein1972bound} as a new method to obtain a Berry--Esseen bound and prove the central limit theorem for weakly dependent data. It has since then become one of the most popular and powerful tools to prove asymptotic normality for dependent data, and different adaptations of it have been proposed, notably the dependency neighborhoods, the exchangeable pairs, the zero-bias coupling, and the size-bias coupling \cite{ross2011fundamentals}. In addition to being used to prove the central limit theorem, it has also been adapted to obtain limit theorems with the Poisson distribution \cite{chen1975poisson} or the exponential distribution \cite{chatterjee2011exponential,pekoz2011new}. Moreover, it has been used for comparing different univariate distributions \cite{ley2017stein}. Our use of the Stein's method is closely related to the dependency neighborhood method described in \cite{ross2011fundamentals}.

% We further remark that the theory we have developed has an interesting by-product. We prove upper bounds on the absolute values of the cumulants of $W_{n}$ (see \cref{thm:corocumubd}). Previously, \cite{janson1988normal} showed a similar bound under the dependency graph conditions and \cite{heinrich1990some,gotze1995m} tightened the results for $m$-dependent random fields. \cite{gotze1983asymptotic,lahiri1993refinements,lahiri1996asymptotic} obtained similar cumulant bounds for the strongly mixing sequences in an effort to obtain Edgeworth expansions. Note that we also provide the bounds on the cumulants of $W_{n}$ for strongly mixing random fields in \cref{sec:mixingmainpart}, which is more general than the results mentioned above. Furthermore, \cite{doring2013moderate,doring2022method} showed that cumulant bounds can be useful in problems including the analysis of moderate deviations.

\subsection{Paper Outline}

In \cref{Notations} we clarify some notations that we use throughout the paper. Then we present our results under two different local dependence conditions in \cref{SEC:LOCAL}. In \cref{m_dependent} and \cref{sec:ustatistic} we respectively apply our results to $m$-dependent random fields and to U-statistics. In \cref{pol_decay} we apply our results to obtain non-uniform Berry--Esseen bounds with polynomial decay. 
In \cref{sec:proofoutline}, we make an overview of our proof techniques. In \cref{SEC:LOCALthmpf} we present the main lemmas (notably \cref{THM:BARBOURLIKE,THM:EXISTENCEXI}) and use them to prove the main result \cref{THM:LOCALWP}. Those lemmas and additional results are proved in \cref{sec:lemma1,sec:lemma2,sec:lemma3}.

\section{General Notations}\label{Notations}
\addtocontents{toc}{\SkipTocEntry}
\subsection*{Notations concerning integers and sets}
In this paper, we will write $\lceil x\rceil$ to denote the smallest integer that is bigger or equal to $x$ and $\lfloor x\rfloor$ denotes the largest integer smaller or equal to $x$. We use $\mathbb{N}$ to denote the set of non-negative integers and let $\mathbb{N}_{+}$ be the set of positive integers. For any $n\in \mathbb{N}_{+}$, denote $[n]:=\{ \ell\in \mathbb{N}_{+}:1\leq \ell \leq n\}$. 
Moreover, for a finite set $B$ we denote by $|B|$ its cardinality.
\addtocontents{toc}{\SkipTocEntry}
\subsection*{Notations for sequences}
Given a sequence $(x_i)$ we will shorthand $x_{1:\ell}=(x_1,\cdots,x_\ell)$ and similarly for any subset $B\subseteq\mathbb{N}_{+}$ we denote $x_B:=(x_i)_{i\in B}$.

\addtocontents{toc}{\SkipTocEntry}
\subsection*{Notations for functions}
For any real valued functions $f(\,\cdot\,),g(\,\cdot\,):\mathbb{N}_{+}\rightarrow \mathbb{R}$, we write $f(n)\lesssim g(n)$ or $f(n)=\mathcal{O}(g(n))$ if there exists some constant $C$ (with dependencies that are fixed in the contexts) and an integer $N>0$ such that the inequality $f(n)\leq C g(n)$ holds for all $n\geq N$. We further write $f(n)\asymp g(n)$ as shorthand for $f(n)\lesssim g(n)$ and $g(n)\lesssim f(n)$.
\addtocontents{toc}{\SkipTocEntry}
\subsection*{Notations for probability distributions}
For a random variable $X$ we write by $\mathcal{L}(X)$ the distribution of $X$.

\section{Main Theorems}\label{SEC:LOCAL}
Let $p\geq 1$ be a positive real number, we write $\omega :=p+1-\lceil p\rceil\in [0,1]$. We choose $I$ to be an infinite index set and $(I_n)_{n=1}^{\infty}$ to be an increasing sequence of finite subsets of $I_1\subseteq I_2\subseteq \cdots \subsetneq I$ that satisfy $|I_n|\xrightarrow{n\rightarrow\infty}\infty$.  \\\noindent Let $\bigl(X^{\scalebox{0.6}{$(n)$}}_{i}\bigr)_{i\in I_{n}}$ be a triangular array of random variables, each row indexed by $i\in I_{n}$ ($n=1,2,\cdots$), we define $W_n$ to be the following empirical average $$W_n:=\sigma_n^{-1}\sum_{i\in I_n}X^{\scalebox{0.6}{$(n)$}}_{i},\quad \text{ with }~\sigma_n^{2}:=\operatorname{Var} \Bigl(\sum_{i\in I_n} X^{\scalebox{0.6}{$(n)$}}_{i}\Bigr).$$ Under the hypothesis that the random variables $\bigl(X^{\scalebox{0.6}{$(n)$}}_i\bigr)$ are locally dependent we will, in this section, bound the Wasserstein-$p$ distance between $W_n$ and its normal limit. The bound we obtain depends on the size of the index set $I_{n}$, the moments of the random variables and the structure of local dependence in question.

% We would like to study the asymptotical properties of $W$ in terms of the Wasserstein metrics under weak dependence conditions on $X_{i}$'s. More specifically, we obtain upper bounds for the Wasserstein-$p$ distance between the distribution of $W$ and a standard normal as the cardinality of the index set goes to infinity.

  %\textbf{[LD-$1$]} to \textbf{[LD-$(\lceil p\rceil+1)$]}, which can be inductively written as follows:
{
  To formally state our conditions on the dependency structure of $(X_i^{(n)})$, we first define the notion of dependency neighborhoods similarly as in \cite{ross2011fundamentals}. 
  \\\noindent Given random variables $(Y_{i})_{i\in I}$ and given $J\subseteq I$, we say that $N(J)\subset I$ is a \textbf{dependency neighborhood} of $J$ if $\{ Y_{j}:j\notin N(J) \}$ is independent of $\{ Y_{j}: j\in J\}$. To state our theorem, we impose that such dependency neighborhoods can be defined for $\bigl(X_{i}^{\scalebox{0.6}{$(n)$}}\bigr)$. More formally, we assume that there is a sequence $(N_{n}(i_{1:q}))_{q}$ of subsets of $I_{n}$ that satisfy the following conditions:
}
\begin{enumerate}[align=left,leftmargin=20pt,itemindent=25pt]
  \item[\textbf{[LD-$\boldsymbol{1}$]}:] For each $i_{1}\in I_{n}$, the subset $N_n(i_{1})\subseteq I_{n}$ is such that $\bigl\{ X^{\scalebox{0.6}{$(n)$}}_{j}:j\notin N_{n}(i_{1}) \bigr\}$ is independent of $X^{\scalebox{0.6}{$(n)$}}_{i_{1}}$.
  \item[\textbf{[LD-$\boldsymbol{q}$]} ($q\geq 2$):] For each $i_{1}\in I_{n}$,~ $i_{2}\in N_n(i_{1})$, $\cdots$,~ $i_{q}\in N_n(i_{1:(q-1)})$, the subset $N_n(i_{1:q})\subset I_n$ is such that $\bigl\{ X^{\scalebox{0.6}{$(n)$}}_{j}:j\notin N_n(i_{1:q}) \bigr\}$ is independent of $\bigl(X^{\scalebox{0.6}{$(n)$}}_{i_{1}},\cdots,X^{\scalebox{0.6}{$(n)$}}_{i_{q}}\bigr)$.
\end{enumerate}
We remark that the sequence of subsets $(N_n({i_{1:q}}))_q$ is increasing, i.e., $N_n(i_{1:(q-1)})\subseteq N_n(i_{1:q})$ in $q$; and that the neighborhoods $N_n(i_{1:q})$ are allowed to be different for different values of $n$--which reflects the triangular array structure of our problem. The condition of dependency neighborhoods here generalizes the one in \cite{ross2011fundamentals} and was also adopted in \cite{fang2019wasserstein}, inspired by \cite{barbour1989central,chen2004normal}. \cite{barbour1989central} obtained a Wasserstein-$1$ bound under ``decomposable'' conditions similar to {[LD-$1$]} and {[LD-$2$]}, and \cite{chen2004normal} showed a Berry--Esseen type result under slightly stronger assumptions for local dependence, while finally \cite{fang2019wasserstein} obtained a Wasserstein-$2$ bound.

  % \textcolor{red}{i think you can delete this. it does not bring a lot of information to the reader}Before presenting the main theorem of this section, we introduce a few notions which are originated in combinatorics \textcolor{blue}{what exactly do you mean by that? what were originated there? do you have refs?} and will be useful later\textcolor{blue}{ in the proof?}.
  In order to define the remainder terms that will appear in our bounds, we introduce the following notions. Given $t\in\mathbb{N}_{+}$, and $\ell\in \mathbb{N}_{+}$  such that $t\geq 2$, we say that the tuple $(\eta_{1},\eta_{2},\cdots,\eta_{\ell})$ is an \textbf{integer composition} of $t$ if and only if $\eta_{1:\ell}$ are positive integers such that $\eta_{1}+\eta_{2}+\cdots+\eta_{\ell}=t$. We denote by $C(t)$ the set of all possible integer compositions
  $$\textstyle C (t):=\bigl\{ \ell,\eta_{1:\ell}\in \mathbb{N}_{+}:\sum_{j=1}^{\ell}\eta_{j}=t \bigr\}.$$
  %$Here for ease of notation we have denoted by $\eta_{1:\ell}$ the variables $\eta_{1},\cdots,\eta_{\ell}$.

  Moreover, for any random variables $(Y_{i})_{i=1}^{t}$, we define the order-$t$ \textbf{compositional expectation} with respect to $\eta_{1:\ell}$ as
  \begin{equation}\label{eq:compexp}
    \begin{aligned}
    &[\eta_{1},\cdots,\eta_{\ell}]\triangleright (Y_{1},\cdots,Y_{t}):=\\
    &\quad\mathbb{E} \bigl[Y_{1}\cdots Y_{\eta_{1}}\bigr]\ \mathbb{E} \bigl[Y_{\eta_{1}+1}\cdots Y_{\eta_{1}+\eta_{2}}\bigr]\ \cdots \ \mathbb{E} \bigl[Y_{\eta_{1}+\cdots+\eta_{\ell-1}+1}\cdots Y_{t}\bigr].
    \end{aligned}
  \end{equation}
Note that if $\eta_{\ell}=1$, the last expectation reduces to $\mathbb{E} [Y_{t}]$. For any positive integer $k$ and real value $\omega\in (0,1]$, we define
\begin{equation}\label{eq:defrmomega}
  \begin{aligned}
  R_{k,\omega,n}:=&\sum_{(\ell,\eta_{1:\ell})\in C^{*}(k+2)}\sum_{i_{1}\in I_{n}}\sum_{i_{2}\in N_n(i_{1})}\cdots\sum_{i_{k+1}\in N_n(i_{1:k})}\\
  &\  [\eta_{1},\cdots,\eta_{\ell}]\triangleright\left(\bigl\lvert X^{\scalebox{0.6}{$(n)$}}_{i_{1}}\bigr\rvert,\cdots,\bigl\lvert X^{\scalebox{0.6}{$(n)$}}_{i_{k+1}}\bigr\rvert,\left(\sum_{i_{k+2}\in N_n(i_{1:(k+1)})}\bigl\lvert X^{\scalebox{0.6}{$(n)$}}_{i_{k+2}}\bigr\rvert\right)^{\omega }\right),
  \end{aligned}
\end{equation}
where $C^{*}(k+2)$ is given by
$$\textstyle C^{*}(t):=\bigl\{(\ell,\eta_{1:\ell})\in \textstyle C(t):~\eta_{j}\geq 2~\text{ for }1\leq j\leq \ell-1,\bigr\}\subseteq \textstyle C(t).$$
The terms $(R_{k,\omega,n})$ are remainder terms that appear in our bound of the Wasserstein-$p$ distance between $W_n$ and its normal limit.

\begin{theorem}\label{THM:LOCALWP}
Let $\bigl(X^{\scalebox{0.6}{$(n)$}}_{i}\bigr)_{i\in I_{n}}$ be a triangular array of mean zero random variables and suppose that they satisfy \emph{[LD-$1$]} to \emph{[LD-($\lceil p\rceil$+$1$)]}. Let $\sigma_n^{2}:=\operatorname{Var} \left(\sum_{i\in I_n} X^{\scalebox{0.6}{$(n)$}}_{i}\right)$ and define $W_n:=\sigma_n^{-1}\sum_{i\in I_n}X^{\scalebox{0.6}{$(n)$}}_{i}$. Further suppose for any $j\in \mathbb{N}_{+}$ such that $j\leq \lceil p\rceil -1$, it holds that $R_{j,1,n}\overset{n\rightarrow \infty}{\longrightarrow} 0$ as $n\to \infty$. Then there exists an integer $N\in \mathbb{N}_{+}$ such that for all $n\ge N$, we have the following Wasserstein bounds:
  \begin{equation}\label{eq:localwp}
    \mathcal{W}_{p}(\mathcal{L}(W_n), \mathcal{N}(0,1)) \leq C_{p} \left (\sum_{j=1}^{\lceil p\rceil-1}R _{j,1,n}^{1/j}+\sum_{j=1}^{\lceil p\rceil}R _{j,\omega,n }^{1/(j+\omega -1)}\right ),
  \end{equation}
  where $\omega=p+1-\lceil p\rceil$ and $C_{p}$ is a constant that only depends on $p$.%\textcolor{red}{can we still talke $C_p$ only depending on $p$ if we removed the $N$? (this is asking if %}
\end{theorem}
\begin{remark}
We note that the condition that the remainder terms $R_{j,1,n}$ shrink to $0$ for all $j\leq \lceil p\rceil -1$ impose an implicit constraint on the size of the sets $N_n(i_{1:q})$.
\end{remark}
In particular, for $p=1,2$ we have
\begin{align}
  &\mathcal{W}_{1}(\mathcal{L}(W_n), \mathcal{N}(0,1))\leq C_{1}R_{1,1,n},\label{eq:localw1}\\
  &\mathcal{W}_{2}(\mathcal{L}(W_n), \mathcal{N}(0,1))\leq C_{2}\bigl(R_{1,1,n}+R_{2,1,n}^{1/2}\bigr).\label{eq:localw2}
\end{align} where the remainders are given by
\begin{align*}
  R_{1,1,n}=&\sigma_n^{-3}\sum_{i\in I_n}\sum_{j\in N_{n}(i)}\sum_{k\in N_{n}(i,j)}\Bigl(\mathbb{E} \bigl[\bigl\lvert X^{\scalebox{0.6}{$(n)$}}_{i}X^{\scalebox{0.6}{$(n)$}}_{j}X^{\scalebox{0.6}{$(n)$}}_{k} \bigr\rvert\bigr]+\mathbb{E} \bigl[\bigl\lvert X^{\scalebox{0.6}{$(n)$}}_{i}X^{\scalebox{0.6}{$(n)$}}_{j} \bigr\rvert\bigr] \ \mathbb{E} \bigl[\bigl\lvert X^{\scalebox{0.6}{$(n)$}}_{k} \bigr\rvert\bigr]\Bigr),\\
  R_{2,1,n}=&\sigma_n^{-4}\sum_{i \in I_n}\sum_{j\in N_{n}(i)}\sum_{k\in N_{n}(i,j)}\sum_{\ell\in N_{n}(i,j,k)}\Bigl(\mathbb{E} \bigl[\bigl\lvert X^{\scalebox{0.6}{$(n)$}}_{i}X^{\scalebox{0.6}{$(n)$}}_{j}X^{\scalebox{0.6}{$(n)$}}_{k}X^{\scalebox{0.6}{$(n)$}}_{\ell} \bigr\rvert\bigr]\\
  &\ +\mathbb{E} \bigl[\bigl\lvert X^{\scalebox{0.6}{$(n)$}}_{i}X^{\scalebox{0.6}{$(n)$}}_{j}X^{\scalebox{0.6}{$(n)$}}_{k} \bigr\rvert\bigr]\ \mathbb{E} \bigl[\bigl\lvert X^{\scalebox{0.6}{$(n)$}}_{\ell} \bigr\rvert\bigr]+\mathbb{E} \bigl[\bigl\lvert X^{\scalebox{0.6}{$(n)$}}_{i}X^{\scalebox{0.6}{$(n)$}}_{j} \bigr\rvert\bigr]\ \mathbb{E} \bigl[\bigl\lvert X^{\scalebox{0.6}{$(n)$}}_{k}X^{\scalebox{0.6}{$(n)$}}_{\ell} \bigr\rvert\bigr]\Bigr).
\end{align*}

Note that \eqref{eq:localw1} was proven by \cite{barbour1989central} and \eqref{eq:localw2} is a corollary of Theorem 2.1, \cite{fang2019wasserstein}. The bound \eqref{eq:localwp} with an integer $p$ was also proposed as a conjecture in \cite{fang2019wasserstein}. As $p$ grows, the right-hand side of \eqref{eq:localwp} becomes more and more complicated, which suggests the necessity of new assumptions in order to obtain a simplified result. We further remark that the choice of $N_n (i_{1:q})$ might not be unique (even if we require that it has the smallest cardinality among all possible index sets that fulfill the assumption [LD-$q$]).% For instance, consider the case when $I=\{ 1,2,3 \}$ and $X_{1},X_{2},X_{3}$ are pairwise independent but not jointly independent. Then $N(1)$ is either $\{ 1,2 \}$ or $\{ 1,3 \}$ by \textbf{[LD-$1$]}. This ambiguity does not hinder the calculation of $R_{k,\omega }$ given specific $(X_{i})_{i\in I}$ since we can arbitrarily fix one possible version of $N(i_{1:q})$'s. However, it might lead to trouble when we try to analyze the summations in general situations because there might be cases where
 ~Therefore, to be able to obtain more interpretable upper-bounds for the remainder terms $(R_{j,\omega,n})$
% \begin{itemize}
%   \item $N(i_{1},i_{2})$ does not coincide with $N(i_{2},i_{1})$;
%   \item $i_{2}\in N(i_{1})$ but $i_{1}\notin N(i_{2})$.
% \end{itemize}
%In order to avoid these troublesome situations and obtain a tractable upper bound for $R_{k,\omega }$
, we impose a slightly stronger assumption on the dependence structure:

\begin{enumerate}[align=left,leftmargin=20pt,itemindent=25pt]
\item[\textbf{[LD*]}:] We suppose that there exists a graph $G_n=(V_n,E_n)$, with $V_{n}:=I_n$ being the vertex set and $E_n$ being the edge set, such that for any two disjoint subsets $J_{1},J_{2}\subseteq I_n$ if there is no edge between $J_{1}$ and $J_{2}$, then $\bigl\{ X^{\scalebox{0.6}{$(n)$}}_{j}:j\in J_{1}\bigr\}$ is independent of $\bigl\{ X^{\scalebox{0.6}{$(n)$}}_{j}:j\in J_{2}\bigr\}$.
\end{enumerate}

{Introduced by \cite{petrovskaya1983central} the graph $G_{n}$ defined above is known as the \textbf{dependency graph} and was later adopted in \cite{janson1988normal,baldi1989normal,ross2011fundamentals}. Please refer to \cite{feray2016mod} for a detailed discussion.}

If \textup{[LD*]} is satisfied, for any subset $J\subseteq I_n$, we define $N_n(J)$ to be the set of vertices in the neighborhood of $J\subseteq I_{n}$ in the graph $G$.
To be precise, this is
\begin{equation*}
  N_n(J):=J\cup \{ i\in I_n: e(i,j)\in E_n\text{ for some }j\in J \},
\end{equation*}
where $e(i,j)$ denotes an edge between the vertices $i$ and $j$.
To simplify the notations, we further denote $N_n(J)$ by $N_n(i_{1:q})$ if $J=\{ i_{1},\cdots,i_{q} \}$ for any $1\leq q\leq \lceil p\rceil+1$.
Then $(N_{n}(i_{1:q}))$ not only satisfies {[LD-$1$]} to {[LD-($\lceil p\rceil$+$1$)]}, but has the following properties as well:
\begin{enumerate}[label=(\alph*)]
  \item $N_n(i_{1:q})=N_n\bigl(i_{\pi(1)},\cdots,i_{\pi(q)}\bigr)$ for any permutation $\pi$ on $\{ 1,\cdots,q \}$;
        \item\label{itm:permute} $i_{q}\in N_n(i_{1:(q-1)})\Leftrightarrow i_{1}\in N_n(i_{2:q})$.
\end{enumerate}
We point out that by definition of the dependency graph even if $\bigl\{ X^{\scalebox{0.6}{$(n)$}}_{j}:j\in J_{1}\bigr\}$ is independent of $\bigl\{ X^{\scalebox{0.6}{$(n)$}}_{j}:j\in J_{2}\bigr\}$, there can still be edges between the vertex sets $J_{1}$ and $J_{2}$. In fact, there might not exist $G_{n}$ such that there is no edge between $J_{1}$ and $J_{2}$ as long as $\bigl\{ X^{\scalebox{0.6}{$(n)$}}_{j}:j\in J_{1}\bigr\}$ is independent of $\bigl\{ X^{\scalebox{0.6}{$(n)$}}_{j}:j\in J_{2}\bigr\}$ since pairwise independence does not imply joint dependence.% In particular, the independence between $X_{i_{1}}$ and $X_{i_{2}}$ does not imply that $i_{2}$ lies outside the neighborhood of $i_{1}$. In the previous case with $X_{1},X_{2},X_{3}$ being pairwise independent but jointly dependent, the graph needs to contain at least two edges between $i_{1},i_{2}$, and $i_{3}$.

%\textcolor{red}{i feel like we can cut this iwthout dramatically changing things}We use the same example to illustrate the difference between the two assumptions. Under \textbf{[LD-$1$]}, $N(i_{1})$, $N(i_{2})$, and $N(i_{3})$ can all be chosen to have cardinality $2$. However, under \textbf{[LD*]}, one of them needs to be $\{ i_{1},i_{2},i_{3} \}$ since there are at least two edges, giving rise to larger $R_{k,\omega}$'s on the right-hand-side of \eqref{eq:localwp}. Nevertheless,
The condition \textup{[LD*]} provides us with a tractable bound on $R_{k,\omega,n}$, which is applicable in most of the commonly encountered settings, including $m$-dependent random fields and U-statistics.

\begin{proposition}\label{THM:LEMMACONTROLBRACKET}
  Given $M\in \mathbb{N}_{+}$ and real number $\omega\in (0,1]$, suppose that $\bigl(X^{\scalebox{0.6}{$(n)$}}_{i}\bigr)_{i\in I_{n}}$ satisfies \textup{[LD*]} and that the cardinality of $N_n(i_{1:(k+1)})$ is upper-bounded by $M<\infty$ for any $i_{1},\cdots,i_{k+1}\in I_n$. Then there exists a constant $C_{k+\omega }$ only depending on $k+\omega$ such that
  \begin{equation*}
    R_{k,\omega,n }\leq C_{k+\omega }  M^{k+\omega }\sum_{i\in I_{n}}\sigma_n^{-(k+1+\omega)}\mathbb{E} \bigl[\bigl\lvert X^{\scalebox{0.6}{$(n)$}}_{i} \bigr\rvert^{k+1+\omega }\bigr].
  \end{equation*}
\end{proposition}
We remark that the upper bound on $(R_{k,\omega,n})$ depends on the moments of the random variables $\bigl(X^{\scalebox{0.6}{$(n)$}}_{i}\bigr)$ and the maximum size of the dependency neighborhoods. The results of \cref{THM:LEMMACONTROLBRACKET} can be used to propose a more interpretable upper bound for the Wasserstein-$p$ distance.
\begin{theorem}\label{THM:LOCALWP2}
Suppose that $\bigl(X^{\scalebox{0.6}{$(n)$}}_{i}\bigr)$ is a triangular array of mean zero random variables satisfying \textup{[LD*]}, and that the cardinality of index set $N_n\bigl(i_{1:(\lceil p\rceil+1)}\bigr)$ is upper-bounded by $M_n<\infty$ for any $i_{1},\cdots,i_{\lceil p\rceil +1}\in I_n$. Furthermore, assume that 
  \begin{equation*}
    M_n^{1+\omega }\sigma_n^{-(\omega+2)}\sum_{i\in I_n}\mathbb{E} \bigl[\bigl\lvert X^{\scalebox{0.6}{$(n)$}}_{i} \bigr\rvert^{\omega +2}\bigr]\to 0,\quad M_n^{p+1}\sigma_n^{-(p+2)}\sum_{i\in I_n}\mathbb{E} \bigl[\bigl\lvert X^{\scalebox{0.6}{$(n)$}}_{i} \bigr\rvert^{p+2}\bigr]\to 0.
  \end{equation*}
  Then there is $N$ such that for all $n\ge N$ we have
  \begin{equation}\label{eq:localwp2}
    \begin{aligned}
    &\mathcal{W}_{p}(\mathcal{L}(W_n),\mathcal{N}(0,1))\\
    \leq &C_{p}\Bigl(M_n^{1+\omega }\sigma_n^{-(\omega+2)}\sum_{i\in I_n}\mathbb{E}\bigl[\bigl\lvert X^{\scalebox{0.6}{$(n)$}}_{i}\bigr\rvert^{\omega +2}\bigr] \Bigr)^{1/\omega }+C_{p}\Bigl(M_n^{p+1}\sigma_n^{-(p+2)}\sum_{i\in I_n}\mathbb{E}\bigl[\bigl\lvert X^{\scalebox{0.6}{$(n)$}}_{i}\bigr\rvert^{p+2}\bigr] \Bigr)^{1/p},
    \end{aligned}
  \end{equation}
  for some constant $C_{p}$ that only depends on $p$.
\end{theorem}
We notably remark that if the moments are nicely behaved in the sense that  
$$B_1:=\sup_{i\in I_n, n\in \mathbb{N}_{+}}\frac{\sqrt{\lvert I_{n} \rvert}\cdot \|X^{\scalebox{0.6}{$(n)$}}_{i}\|_{p+2}}{\sigma_{n}}<\infty,$$  
and that the size of the dependency neighborhood are universally bounded, i.e., $$B_2:=\sup_{i_{1:(\lceil p\rceil +1)}\in I_n,n\in\mathbb{N}_{+}} \bigl|N_n(i_{1:(\lceil p\rceil +1)})\bigr|<\infty,$$
then there is a constant $K_p$ that only depends on $B_1$, $B_2$ and $p\ge1$ such that for $n$ large enough we have
$$ \mathcal{W}_{p}(\mathcal{L}(W_n),\mathcal{N}(0,1))\leq \frac{K_{p}}{\sqrt{|I_n|}}.$$
The rate of convergence matches the known rate for independent random variables (see \cite{bobkov2018berry}).
\section{Results for m-Dependent Random Fields}\label{m_dependent}

Let $d\in \mathbb{N}_+$ be a positive integer, in this section, we study $d$-dimensional random fields.

\begin{definition}[$m$-Dependent Random Field]\label{thm:defmdepfield}
  A random field $(X_{i})_{i \in T}$ on $T \subseteq \mathbb{Z}^{d}$ is $m$-dependent if and only if for any subsets $U_{1}, U_{2} \subseteq \mathbb{Z}^{d}$, the random variables $(X_{i_{1}})_{i_{1} \in U_{1}\cap T}$ and $(X_{i_{2}})_{i_{2}\in U_{2}\cap T}$ are independent whenever $\lVert i_{1}-i_{2}\rVert >m$ for all $i_{1} \in U_{1}$ and $i_{2} \in U_{2}$.

  Here $\lVert\cdot\rVert$ denotes the maximum norm on $\mathbb{Z}^{d}$, that is $\lVert \boldsymbol{z}\rVert=\max _{1 \leq j \leq d}\lvert z_{j}\rvert$ for $\boldsymbol{z}=(z_{1}, \cdots, z_{d})$.
\end{definition}
% \textcolor{red}{i think it is fine for this setion to focus on the non triangular array case. let me know}
 Now we consider an increasing sequence $T_1\subseteq T_2\subseteq \cdots$ of finite subsets of $\mathbb{Z}^d$ that satisfy $|T_n|\xrightarrow{n\rightarrow \infty}\infty$. We have the following result as a corollary of \cref{THM:LOCALWP2}. 
\begin{corollary}\label{THM:MDEPFIELD}
  Let $p\in \mathbb{N}_+$ and $m\in\mathbb{N}_{+}$ be positives integer. 
  Suppose that $\bigl(X^{\scalebox{0.6}{$(n)$}}_{i}\bigr)$ is a triangular array where each row is an $m$-dependent random field indexed by finite subsets $T_{n}\subseteq \mathbb{Z}^{d}$ such that $|T_n|\xrightarrow{n\rightarrow \infty}\infty$. Let $\sigma_{n}^{2}:=\operatorname{Var}\Bigl(\sum_{i\in T_{n}}X_{i}^{\scalebox{0.6}{$(n)$}}\Bigr)$ and define $W_{n}:=\sigma_{n}^{-1}\sum_{i\in T_{n}}X_{i}^{\scalebox{0.6}{$(n)$}}$. Further suppose that $\mathbb{E} \bigl[X^{\scalebox{0.6}{$(n)$}}_{i}\bigr]=0$ for any $i\in T_{n}$ and that the following conditions hold:
  \begin{itemize}
    \item Moment condition: $\sigma_n^{-(p+2)}\sum_{i\in T_n}\mathbb{E}\bigl[\bigl\lvert X^{\scalebox{0.6}{$(n)$}}_{i}\bigr\rvert^{p+2}\bigr] \to 0 \quad\text{ as } n\to \infty$;
    \item Non-degeneracy condition: $\limsup_{n}\sigma_{n}^{-2}\sum_{i\in T_{n}}\mathbb{E} \bigl[\bigl\lvert X_{i}^{\scalebox{0.6}{$(n)$}}\bigr\rvert^{2}\bigr]\leq M<\infty$ for some $M\geq 1$.
  \end{itemize}
  % \begin{equation}\label{eq:coro32condition}
  %   \sigma_n^{-(p+2)}\sum_{i\in T_n}\mathbb{E}[\lvert X_{i}\rvert^{p+2}] \to 0 \quad\text{ as } n\to \infty,
  % \end{equation}
  % where $\sigma_n^{2}:=\operatorname{Var} \bigl(\sum_{i\in T_n}X_{i}\bigr)$.
  % Denoting $W_n:=\sigma_n^{-1}\sum_{i\in T_n}X_{i}$, f
  Then for $n$ large enough, we have
  \begin{equation}\label{eq:mdepfield}
    \mathcal{W}_{p}(\mathcal{L}(W_n),\mathcal{N}(0,1))\leq C_{p,d}m^{\frac{(1+\omega)d}{\omega}}M^{\frac{p-\omega}{p\omega}}\sigma_n^{-\frac{p+2}{p}}\left (\sum_{i\in T_n}\mathbb{E}\bigl[\bigl\lvert X^{\scalebox{0.6}{$(n)$}}_{i}\bigr\rvert^{p+2}\bigr] \right )^{1/p},
  \end{equation}
  where $C_{p,d}$ only depends on $p$ and $d$.

  In particular, for a triangular array of $m$-dependent \emph{stationary} random fields, suppose that we have $\sup_{n}\mathbb{E} \bigl[\bigl\lvert X_{i}^{\scalebox{0.6}{$(n)$}} \bigr\rvert^{p+2}\bigr]<\infty$, and that the non-degeneracy condition $\liminf_{n}\sigma_n^{2}/\lvert T_n \rvert>0$ holds. Then we have
  \begin{equation*}
    \mathcal{W}_{p}(\mathcal{L}(W_n),\mathcal{N}(0,1))=\mathcal{O} (\lvert T_n \rvert^{-1/2}).
  \end{equation*}
\end{corollary}

\section{Application to U-Statistics}\label{sec:ustatistic}

\begin{definition}[U-Statistic]\label{thm:defustat}
  Let $(X_{i})_{i=1}^{n}$ be a sequence of \emph{i.i.d.} random variables. Fix $m\in\mathbb{N}_{+}$ such that $m\geq 2$. Let $h:\mathbb{R}^{m}\to \mathbb{R}$ be a fixed Borel-measurable function. The Hoeffding U-statistic is defined as
  \begin{equation*}
    \sum_{1\leq i_{1}\leq \cdots\leq i_{m}\leq n}h\bigl(X_{i_{1}},\cdots,X_{i_{m}}\bigr).
  \end{equation*}
\end{definition}

\begin{corollary}\label{THM:USTATWP}
  Given $p\geq 1$, suppose that the U-statistic of an \emph{i.i.d.} sequence $(X_{i})_{i=1}^{n}$ induced by a symmetric function $h:\mathbb{R}^{m}\to\mathbb{R}$ satisfies the following conditions
  \begin{itemize}
    \item Mean zero: $\mathbb{E} \bigl[h(X_{1}, \cdots, X_{m})\bigr]=0$;
    \item Moment condition: $\mathbb{E} \bigl[\bigl\lvert h(X_{1}, \cdots, X_{m})\bigr\rvert^{p+2}\bigr]<\infty$;
    \item Non-degeneracy condition: $\mathbb{E} [g(X_{1})^{2}]>0$, where $g(x):=\mathbb{E} \bigl[h(X_{1},\cdots,X_{m})\bigm\vert X_{1}=x\bigr]$.
  \end{itemize}
  If we let
  \begin{equation*}
    W_{n}:=\frac{1}{\sigma_{n}}\sum_{1\leq i_{1}\leq \cdots\leq i_{m}\leq n}h\bigl(X_{i_{1}},\cdots,X_{i_{m}}\bigr),
  \end{equation*}
  where
  \begin{equation*}
    \sigma_{n}^{2}:=\operatorname{Var} \left(\sum_{1\leq i_{1}\leq \cdots\leq i_{m}\leq n} h\bigl(X_{i_{1}},\cdots,X_{i_{m}}\bigr)\right),
  \end{equation*}
  the following Wasserstein bound holds:
  \begin{equation*}
    \mathcal{W}_{p}(\mathcal{L}(W_{n}),\mathcal{N}(0,1))=\mathcal{O}(n^{-1/2}).
  \end{equation*}
\end{corollary}
% I've modified this section after EJP submission - Tyler
\section{Application to Non-Uniform Berry--Esseen Bounds}\label{pol_decay}
In this section, we show a specific application of our results to non-uniform Berry--Esseen bounds with polynomial decay of any order. Mirroring the classical literature, \cite{chen2004normal} established Berry--Esseen bounds for locally dependent random variables. Notably, their Theorem 2.4 showed that if the random variables $(X_i^{(n)})$ satisfy some boundedness condition on the dependency neighborhoods, then there is a constant $C>0$ such that 
\begin{equation*}
\sup_t\bigl|\mathbb{P}(W_n\geq t)-\Phi^{c}(t)\bigr|\leq C \sum_{i\in I_{n}}\bigl\lVert X_{i}^{(n)} \bigr\rVert _{3}^{3} / \sigma_{n}^{3},
\end{equation*}
where $\Phi^{c}(t):=\mathbb{P} (Z\geq t)$ with $Z\sim \mathcal{N}(0,1)$.
This extends the classical Berry--Esseen bound to locally dependent random variables, and can potentially be used to construct Kolmogorov–Smirnov tests under local dependence in nonparametric inference. However, one of the drawbacks of this inequality is that it does not depend on $t$. One would imagine that for large $t$ we could find tighter bounds for $\bigl|\mathbb{P}(W_n\geq t)-\Phi^{c}(t)\bigr|$.  Non-uniform Berry--Esseen bounds establish this. Notably \cite{chen2004normal} (Theorem 2.5) showed that under the above conditions, there exists some universal constant $C'$ such that 
\begin{equation*}
\bigl|\mathbb{P}(W_n\geq t)-\Phi^{c}(t)\bigr|\leq \frac{C'}{1+|t|^3}\sum_{i\in I_{n}}\bigl\lVert X_{i}^{(n)} \bigr\rVert _{3}^{3} / \sigma_{n}^{3},\quad \forall t\in \mathbb{R}.
\end{equation*}
This bound does decrease as $|t|$ increases and does so at a rate of $|t|^{-3}$. However one would expect that this rate could be tightened if additional assumptions were made about the moments of $\bigl(X_i^{(n)}\bigr)$. If the random variables admit some exponential moments then \cite{liu2021cram} demonstrated that locally dependent random variables admit moderate deviation inequalities. In this section, we show how $\mathcal{W}_p$ bound
 can help us obtain bounds that decrease polynomially fast in $t$ at a small price in its dependence on $|I_n|$, and do so without assuming infinite moments.%To be specific, we explain how $\mathcal{W}_p(\mathcal{L}(W_n),\mathcal{N}(0,1))$ ($p>1$)can be used to obtain tail bounds for $W_n$. 
% Indeed, letting $t\ge 0$, an important goal for statistical inference is to obtain tight upper bounds for $\mathbb{P}(W_n\ge t)$. Notably, among many other applications, such inequalities lie at the heart of decision-making in reinforcement learning \cite{mnih2008empirical,audibert2009exploration} and generalization guarantees for high dimensional statistics \cite{wainwright2019high,bartlett2002rademacher}. When the observations $(X_i)$ are \emph{i.i.d.} and bounded, general well-known concentration inequalities such as Azuma's or Bernstein's inequalities allow one to upper-bound $\mathbb{P}(W_n\ge t)$ with a sub-Gaussian decay in $t$. Notably if  $\|X_1\|_{\infty}\le R$ and $\operatorname{Var}(X_1)=1$, then Azuma's inequality states that $\mathbb{P}(W_n \ge t)\le e^{-\frac{nt^2}{2R^2}}$. However, no such concentration inequality is known to hold for strongly mixing sequences. Moreover, even when sub-Gaussian bounds hold, those are known to be significantly looser than the asymptotically valid tail bound $\Phi^{c}(t)$ \cite{austern2022efficient}, where $\Phi^{c}(\,\cdot\,)=1-\Phi(\,\cdot\,)$ and $\Phi(\,\cdot\,)$ denotes the cumulative distribution function (CDF) of the standard normal. In this section we see that our results can be used to obtain a tail bound that is asymptotically tight, decreases polynomially fast in $t$, and is valid for mixing sequences. 

\begin{theorem}\label{gillette}
    We assume that the conditions of \cref{THM:LOCALWP2} are satisfied. There is a constant $C>0$ such that for all $\beta>0$ and $t>0$ satisfying %$\sup_{i\in I}\|X_i^{(n)}/\sigma_n\|_{p+2}$
    \begin{equation*}
    (\sqrt{2\pi }p)^{\frac{1}{p+1}}\Bigl(1-\frac{\sqrt{2\beta\log t}}{t}\Bigr)t^{1-\frac{\beta}{p+1}}\ge \mathcal{W}_p(\mathcal{L}(W_n),\mathcal{N}(0,1)),
    \end{equation*} 
    we have 
    \begin{align*}
    - \frac{C}{t}\varphi\Bigl(t\bigl({1-\frac{1}{p+1}}\bigr)\Bigr)   \leq \frac{\mathbb{P}(W_n\geq t)-\Phi^c(t)}{\mathcal{W}_p(\mathcal{L}(W_n),\mathcal{N}(0,1))^{1-\frac{1}{p+1}}}\leq \frac{C}{t^{1+\beta\bigl(1-\frac{1}{p+1}\bigr)}},
    \end{align*}
    where $\varphi$ is the density function of $\mathcal{N}(0,1)$.
    % \begin{align}
    % - \frac{p^{\frac{1}{p+1}}}{t}\varphi\Big(t({1-\frac{1}{p+1}})\Big)\big(1+\frac{1}{p}\big)   \le \frac{1}{\mathcal{W}_p(\mathcal{L}(W_n),\mathcal{N}(0,1))^{1-\frac{1}{p+1}}}\Big[P(W_n\ge t)-\Phi^c(t)\Big]\le \frac{1}{t^{1+\beta\big(1-\frac{1}{p+1}\big)}\sqrt{n}^{1-\frac{1}{p+1}}}{p^{\frac{1}{p+1}}K_{p,n}^{1-\frac{1}{p+1}}}(1+\frac{1}{p})
    % \end{align}
\end{theorem}
We can see from this result that the quantity $\lvert \mathbb{P} (W_{n}\geq t)-\Phi^{c}(t) \rvert$ decays in both $t$ and $n$. Notably given any $p,r\geq 1$ assuming that the $(p+2)$-th moments of $X _{i}$'s and the dependency neighborhoods are bounded in the sense that 
$$\sup_{n\in \mathbb{N}^{+},i_{1:(p +1)}\in I_n} \bigl|N_n(i_{1:(p +1)})\bigr|<\infty,$$ 
we have $\lvert \mathbb{P} (W_{n}\geq t)-\Phi^{c}(t) \rvert=o\bigl(t^{-r}\lvert I_{n} \rvert^{-\frac{p}{2(p+1)}}\bigr)$ for $t$ and $n$ large enough.
In particular, for $p\in \mathbb{N}^{+}$ \cref{gillette,THM:LOCALWP2} imply the uniform Berry--Esseen bound by taking the supremum over $t$:
\begin{equation*}
\sup_{t\in \mathbb{R}}\lvert \mathbb{P}(W_n\geq t)- \Phi^c(t)\rvert \leq  C \biggl(\sum_{i\in I_{n}}\bigl\lVert X_{i}^{(n)} \bigr\rVert _{3}^{3} / \sigma_{n}^{3}+\sum_{i\in I_{n}}\bigl\lVert X_{i}^{(n)} \bigr\rVert _{p+2}^{1+2/p} / \sigma_{n}^{1+2/p}\biggr).
\end{equation*}
Note that it recovers the uniform Berry--Esseen bound in \cite{chen2004normal} with $p=1$.

\section{Overview of the Proofs}\label{sec:proofoutline}
The key idea of our proofs is to approximate the sum of weakly dependent random variables $\bigl(X_{i}^{\scalebox{0.6}{$(n)$}}\bigr)_{i\in I_{n}}$ by the empirical average of $q_{n}$ \emph{i.i.d.} random variables $\xi_{1}^{\scalebox{0.6}{$(n)$}},\cdots,\xi_{q_{n}}^{\scalebox{0.6}{$(n)$}}$ which we denote $V_{n}:=\frac{1}{\sqrt{q_{n}}}\sum_{i=1}^{q_{n}}\xi_{i}^{\scalebox{0.6}{$(n)$}}$. More specifically we aim for the Wasserstein-p distance between them
 $\mathcal{W}_{p}(\mathcal{L}(W_{n}),\mathcal{L}(V_{n}))$ to be as small as possible. To establish the desired result we then exploit the triangle inequality that guarantees that \begin{equation*}
  \mathcal{W}_{p}(\mathcal{L}(W_{n}),\mathcal{N}(0,1))\leq \mathcal{W}_{p}(\mathcal{L}(W_{n}),\mathcal{L}(V_{n}))+\mathcal{W}_{p}(\mathcal{L}(V_{n}),\mathcal{N}(0,1)),
\end{equation*} and we use previously known bounds for $\mathcal{W}_{p}(\mathcal{L}(V_{n}),\mathcal{N}(0,1))$ (\cref{thm:lemiidwp}).%Note that $q_{n}$ is a number that is not necessarily equal to $n$ but does diverge to infinity with $n$. 

To be able to show that such random variables $\xi_{1}^{\scalebox{0.6}{$(n)$}},\cdots,\xi_{q_{n}}^{\scalebox{0.6}{$(n)$}}$ exist, we first show (\cref{THM:EXISTENCEXI}) that as long as the third and higher-order cumulants of $W_{n}$ decay then there exist integers $(q_{n})$ and \emph{i.i.d.} random variables such that the first $k$ ($k\in \mathbb{N}_{+}$) cumulants of 
$$V_{n}:=\frac{1}{\sqrt{q_{n}}}\sum_{i=1}^{q_{n}}\xi_{i}^{\scalebox{0.6}{$(n)$}}$$
matches those of $W_{n}$ for $n$ large enough. The decay of the cumulants can be proven to hold by exploiting the local dependence assumptions (see \cref{thm:corocumubd}).

As a reminder, our goal is to establish that the Wasserstein distance $\mathcal{W}_p(\mathcal{L}(V_n),\mathcal{L}(W_n))$ is small. We relate this to the cumulants thanks to the fact that the Wasserstein-$p$ distance can be upper-bounded by integral probability metrics (\cref{thm:lemzolo}) and the well-known Stein equation. %Indeed we establish local expansions of the latter to orders that will depend on our choice of $p\geq 1$. %Indeed, while Stein directly upper-bounded $\mathbb{E} [f'(W_{n})-W_{n}f(W_{n})]$ to derive a Wasserstein-$1$ result for dependent data, we first develop local expansions of it to a higher-order that depend on our choice of $p\geq 1$. 
Indeed for \emph{i.i.d.} random variables $\bigl(\xi_{i}^{\scalebox{0.6}{$(n)$}}\bigr)_{i=1}^{q_{n}}$, \cite{barbour1986asymptotic} showed that the following approximation holds (restated in \cref{thm:barbour})
\begin{equation}\label{eq:demo1}
\begin{aligned}
  &\mathbb{E}[ h(V_{n})]-\mathcal{N}h=\mathbb{E} [f'_h(V_{n})-V_{n}f_h(V_{n})]\\
  = &\sum_{(r,s_{1:r})\in \Gamma (\lceil p\rceil-1)}(-1)^{r}\prod_{j=1}^{r}\frac{\kappa _{s _{j}+2}(V_{n})}{(s _{j}+1)!}\mathcal{N}\ \Bigl[\prod_{j=1}^{r}(\partial ^{s _{j}+1}\Theta)\ h\Bigr] +\text{ Remainders},
\end{aligned}
\end{equation}
where $f_h$ is the solution of the Stein equation \eqref{eq:stein} and $\kappa_{j}(\,\cdot\,)$ denotes the $j$-th cumulant of a random variable. (All the other notations in \eqref{eq:demo1} will be made clear in \cref{SEC:LOCALthmpf}.) We show that we can obtain similar expansions for $\mathbb{E} [f'(W_{n})-W_{n}f(W_{n})]$ (see \cref{THM:BARBOURLIKE}):
\begin{equation}\label{eq:demo3}
\begin{aligned}
  &\mathbb{E}[ h(W_{n})]-\mathcal{N}h=\mathbb{E} [f'_h(W_{n})-W_{n}f_h(W_{n})]\\
  = &\sum_{(r,s_{1:r})\in \Gamma (\lceil p\rceil-1)}(-1)^{r}\prod_{j=1}^{r}\frac{\kappa _{s _{j}+2}(W_{n})}{(s _{j}+1)!}\mathcal{N}\ \Bigl[\prod_{j=1}^{r}(\partial ^{s _{j}+1}\Theta)\ h\Bigr] +\text{ Remainders},
\end{aligned}
\end{equation}As mentioned in the previous paragraph, $q_{n}$ and $\xi_{i}^{\scalebox{0.6}{$(n)$}}$ can be chosen to be such that $\kappa_{j}(V_{n})=\kappa_{j}(W_{n})$ for $j=1,\cdots,\lceil p\rceil+1$. Thus, by taking the difference of \eqref{eq:demo1} and \eqref{eq:demo2}, we get an upper bound on $\bigl\lvert\mathbb{E}[h(W_{n})]-\mathbb{E} [h(V_{n})]\bigr\rvert$ for a large class of function $h$. As shown in \cref{thm:lemzolo}, this allows us to obtain an upper bound of the Wasserstein-$p$ distance between $\mathcal{L}(W_{n})$ and $\mathcal{L}(V_{n})$ for a general $p\geq 1$. The desired result is therefore implied by the triangle inequality of the Wasserstein-$p$ distance
\begin{equation*}
  \mathcal{W}_{p}(\mathcal{L}(W_{n}),\mathcal{N}(0,1))\leq \mathcal{W}_{p}(\mathcal{L}(W_{n}),\mathcal{L}(V_{n}))+\mathcal{W}_{p}(\mathcal{L}(V_{n}),\mathcal{N}(0,1)),
\end{equation*}
and the already known Wasserstein-$p$ bounds for \emph{i.i.d.} random variables (\cref{thm:lemiidwp}).% \textcolor{red}{should i delete this to avoid repetition?}

To be able to show that \eqref{eq:demo2} holds, we develop new techniques to obtain such expansions, which will be carefully elaborated and discussed in \cref{sec:lemma1}. 
\appendix
\section{Adapting Stein's Method for Wasserstein-$p$ Bounds}\label{SEC:LOCALthmpf}
In this section, we provide the proofs of \cref{THM:LOCALWP,THM:LOCALWP2} using Stein's method. We first introduce some background definitions and lemmas before showing the proofs of the main theorems.

\addtocontents{toc}{\SkipTocEntry}
\subsection{Preliminaries and Notations}
%Firstly, we remind the readers of the definitions of Hölder spaces, which will be heavily used throughout the rest of the paper.
\begin{definition}[Hölder Space]\label{thm:defholder}
  For any $k\in \mathbb{N}$ and real number $\omega\in (0,1]$, the Hölder space $\mathcal{C}^{k,\omega }(\mathbb{R})$ is defined as the class of $k$-times continuously differentiable functions $f: \mathbb{R} \to \mathbb{R}$ such that the $k$-times derivative of $f$ is $\omega $-Hölder continuous, i.e.,
  \begin{equation*}
    \lvert f\rvert_{k, \omega }:=\sup _{x \neq y \in \mathbb{R}} \frac{\lvert \partial^{k} f(x)-\partial^{k} f(y)\rvert}{\lvert x-y\rvert^{\omega }}<\infty,
  \end{equation*}
  where $\partial$ denotes the differential operator. Here $\omega $ is called the Hölder exponent and $\lvert f \rvert_{k,\omega }$ is called the Hölder coefficient.
\end{definition}

Using the notions of Hölder spaces, we define the Zolotarev's ideal metrics, which are related to the Wasserstein-$p$ distances via \cref{thm:lemzolo}.

\begin{definition}[Zolotarev Distance]\label{thm:defzolo}
  Suppose $\mu$ and $\nu$ are two probability distributions on $\mathbb{R}$. For any $p>0$ and $\omega :=p+1-\lceil p\rceil\in (0,1]$, the Zolotarev-$p$ distance between $\mu$ and $\nu$ is defined by
  $$
    \mathcal{Z}_{p}(\mu, \nu):=\sup _{f\in\Lambda_{p}}\biggl(\int_{\mathbb{R}} f(x) \dif \mu(x)-\int_{\mathbb{R}} f(x) \dif \nu(x)\biggr),
  $$
  where $\Lambda_{p}:=\{ f \in \mathcal{C}^{\lceil p\rceil-1,\omega }(\mathbb{R}):\lvert f \rvert_{\lceil p\rceil-1,\omega }\leq 1 \}$
\end{definition} We will see in \cref{thm:lemzolo} how the Zolotarev distance can be used to obtain $\mathcal{W}_p(\cdot,\cdot)$ rates. To bound $\mathcal{Z}_p(\,\cdot\,,\,\cdot\,)$ we rely on the Stein's method which was introduced by \cite{stein1972bound} in order to prove the central limit theorem for dependent data. It has been widely adapted to all kinds of normal approximation problems. See \cite{ross2011fundamentals}
for a detailed exposition.

\addtocontents{toc}{\SkipTocEntry}
\subsection*{Stein equation and its solutions}
Let $Z\sim \mathcal{N}(0,1)$ be a standard normal random variable. For any measurable function $h:\mathbb{R}\to \mathbb{R}$, if $h(Z)\in \mathcal{L}^{1}(\mathbb{R})$, we write $\mathcal{N} h:=\mathbb{E} [h(Z)].$ Thus, $h(Z)\in \mathcal{L}^{1}(\mathbb{R})$ if and only if $\mathcal{N}|h|<\infty$. Moreover, we define $f_h(\,\cdot\,)$ by
\begin{equation}\label{eq:steinsol}
  f_{h}(w) :=\int_{-\infty}^{w} e^{(w^{2}-t^{2})/2}(h(t)-\mathcal{N} h) \dif t
  =-\int_{w}^{\infty} e^{(w^{2}-t^{2})/2}(h(t)-\mathcal{N} h) \dif t .
\end{equation}We remark that $f_h(\cdot)$ is a \emph{solution} of the Stein equation meaning that it satisfies
\begin{equation}\label{eq:stein}
  f_h'(w)-w f_h(w)=h(w)-\mathcal{N} h,\qquad \forall w\in\mathbb{R}.
\end{equation}
Bounding $\Big|\mathbb{E}(f'_h(W_n)-W_nf_h(W_n))\Big|$ therefore allows to control $\Big|\mathbb{E}(h(W_n))-\mathcal{N}h\Big|$. If we do this for a large class of functions $h$ we can therefore upper-bound the Zolotarev distance between $\mathcal{L}(W_n)$ and a normal distribution. This is the key idea behind the Stein's method. For notational convenience, we denote by $\Theta$ the operator that maps $h$ to $f _{h}$ for any $h$ such that $\mathcal{N}\lvert h \rvert< \infty$, i.e.,
$$
  \Theta h=f _{h}.
$$
Note that $\Theta h(\,\cdot\,)$ is a function. If $h\in \Lambda_p$, then we see in \cref{thm:lemsteinsol} that $\Theta h$ can be bounded. \begin{lemma}[Part of Lemma 6 of \cite{barbour1986asymptotic}]\label{thm:lemsteinsol}
  For any $p>0$, let $h \in \Lambda_{p}$ be as defined in \cref{thm:defzolo}. Then $\Theta h=f_{h}$ in \eqref{eq:steinsol} is a solution to \eqref{eq:stein}. Moreover, $\Theta h\in \mathcal{C}^{\lceil p\rceil-1,\omega}(\mathbb{R})\cap\mathcal{C}^{\lceil p\rceil,\omega }(\mathbb{R})$ and the Hölder coefficients $\lvert \Theta h \rvert_{\lceil p\rceil-1,\omega}$ and $\lvert \Theta h \rvert_{\lceil p\rceil,\omega }$ are bounded by some constant only depending on $p$.
\end{lemma}
%\end{definition}

\addtocontents{toc}{\SkipTocEntry}
\subsection{Key Lemmas}
First, we present an important result that states that the Wasserstein-$p$ distance can be controlled in terms of the Zolotarev distance.
\begin{lemma}[Theorem $3.1$ of \cite{rio2009upper}]\label{thm:lemzolo}
  For any $p\geq 1$, there exists a positive constant $C_{p}$, such that for any pair of distributions $\mu,\nu$ on $\mathbb{R}$ with finite absolute moments of order $p$ such that
  $$
    \mathcal{W}_{p}(\mu, \nu) \leq C_{p}\bigl(\mathcal{Z}_{p}(\mu, \nu)\bigr)^{1/p}.
  $$
  In particular, $\mathcal{W}_{1}(\mu,\nu)=\mathcal{Z}_{1}(\mu,\nu)$ by Kantorovich--Rubinstein duality.
\end{lemma}

%To bound $\mathcal{Z}_p(\,\cdot\,,\,\cdot\,)$ we will use the Stein's method and exploit the fact $\Theta h$ can be controlled for any function $h\in \Lambda_p$.

In the next two lemmas, we present already-known results for the normal approximation of sums of independent random variables. Firstly  \cref{thm:barbour} provides an expansion for the difference between $\mathbb{E}[h(S_n)]$, where $S_n$ is an empirical average and $\mathcal{N}h$. This lemma will allow us to relate the Zolotarev distance to the cummulants.  %\cref{thm:lemiidwp} gives an upper bound on the Wasserstein distance between the distribution of this empirical average, $S_n$, and the standard normal distribution.

\begin{lemma}[Theorem 1 of \cite{barbour1986asymptotic}]\label{thm:barbour}
  For any $p>0$, let $h \in \Lambda_{p}$ and $S_{n}:=\sum_{i=1}^{n} X_{i}$ where $\left\{X_{1}, \cdots, X_{n}\right\}$ are independent, with $\mathbb{E} [X_{i}]=0$ and $\mathbb{E} [S_{n}^{2}]=1$. Then it follows that
  \begin{equation}
  \begin{aligned}\label{eq:barbour}
    &\mathbb{E}[ h(S_{n})]-\mathcal{N}h
    =\\ 
    &\quad \sum_{(r,s_{1:r})\in \Gamma (\lceil p\rceil-1)}(-1)^{r}\prod_{j=1}^{r}\frac{\kappa _{s _{j}+2}(S_{n})}{(s _{j}+1)!}\mathcal{N}\ \Bigl[\prod_{j=1}^{r}(\partial ^{s _{j}+1}\Theta)\ h\Bigr] +\mathcal{O}\biggl( \sum_{i=1}^{n}\mathbb{E} [\lvert X_{i} \rvert ^{p+2}]\biggr),
  \end{aligned}
\end{equation}
 where the first sum is over $\Gamma (\lceil p\rceil -1):=\bigl\{ r, s _{1:r}\in \mathbb{N}_{+}:\sum_{j=1}^{r}s _{j}\leq \lceil p\rceil-1\bigr\}$.
\end{lemma}

Note that there is a slight abuse of notation in \eqref{eq:barbour}. The last $\prod$ indicates the composition of the operators in the parentheses rather than the product. 

Secondly \cref{thm:barbour} gives an upper bound on the Wasserstein distance between the distribution of this empirical average $S_n$, and the standard normal distribution. This lemma will guarantee that if an approximation of $W_n$ by a sum of independent random variables $V_n$ can be obtained then $V_n$ is approximately normally distributed.

\begin{lemma}[Theorem 1.1 of \cite{bobkov2018berry}]\label{thm:lemiidwp}
  For any $p\geq 1$, let $S_{n}:=\sum_{i=1}^{n} X_{i}$ where $\{X_{1}, \cdots, X_{n}\}$ are independent and satisfy that $\mathbb{E} [X_{i}]=0$ and $\mathbb{E} [S_{n}^{2}]=1$. Then it follows that
  \begin{equation}\label{eq:iidwp}
    \mathcal{W}_{p}(\mathcal{L}(S_{n}), \mathcal{N}(0,1)) \leq C_{p}\biggl(\sum_{i=1}^{n} \mathbb{E}[\lvert X_{i}\rvert ^{p+2}]\biggr)^{1/p},
  \end{equation}
  where $C_{p}$ continuously depends on $p$.
\end{lemma}

We now introduce two new lemmas crucial in the proof of \cref{THM:LOCALWP}. They will be proven in \cref{sec:lemma1} and \cref{sec:lemma2}. The first lemma generalizes \cref{thm:lemiidwp} to the dependent setting.

\begin{lemma}[Local Expansion]\label{THM:BARBOURLIKE}
  Suppose that $\bigl(X^{\scalebox{0.6}{$(n)$}}_{i}\bigr)_{i\in I_n}$ is a triangular array of random variables with dependency neighborhoods satisfying the local dependence conditions \emph{[LD-$1$]} to \emph{[LD-($\lceil p\rceil$+$1$)]}. Let $W_n:=\sum_{i\in I_n}X^{\scalebox{0.6}{$(n)$}}_{i}$ with $\mathbb{E} \bigl[X_{i}^{\scalebox{0.6}{$(n)$}}\bigr]=0$, $\mathbb{E} [W_n^{2}]=1$.
  Then for any $p>0$ and $h \in \Lambda_{p}$, we have
  \begin{equation}\label{eq:barbourlike}
    \begin{aligned}
    \mathbb{E} [h(W_n)]-\mathcal{N}h
    =&
    \sum_{(r,s_{1:r})\in \Gamma (\lceil p\rceil -1)}(-1)^{r}\prod_{j=1}^{r}\frac{\kappa _{s _{j}+2}(W_n)}{(s _{j}+1)!}\mathcal{N}\ \Bigl[\prod_{j=1}^{r}(\partial ^{s _{j}+1}\Theta)\ h\Bigr]\\
    &\quad+\mathcal{O}\biggl(\sum_{j=1}^{\lceil p\rceil-1}R _{j,1,n}^{p/j}+\sum_{j=1}^{\lceil p\rceil}R _{j,\omega,n }^{p/(j+\omega -1)}\biggr),
    \end{aligned}
  \end{equation}
  where the first sum is over $\Gamma (\lceil p\rceil -1):=\bigl\{ r, s _{1:r}\in \mathbb{N}_{+}:\sum_{j=1}^{r}s _{j}\leq \lceil p\rceil-1\bigr\}$.
\end{lemma}

We can see that \cref{thm:barbour,THM:BARBOURLIKE} look quite similar to one another with the only differences being the dependence structures of $\bigl(X^{\scalebox{0.6}{$(n)$}}_{i}\bigr)$ and the remainder terms in the expansions. This similarity inspires the proof of \cref{THM:LOCALWP}. To illustrate this, imagine that there would exist some \emph{i.i.d.} random variables $\bigl(\xi_{i}^{\scalebox{0.6}{$(n)$}}\bigr)_{i=1}^{q_{n}}$ and a large sample size $q_{n}$ such that the first $\lceil p\rceil$+$1$ cumulants of $V_n:=q_{n}^{-1/2}\sum_{i=1}^{q_{n}}\xi_{i}^{\scalebox{0.6}{$(n)$}}$ match with those of $W_n$, then the expansion \eqref{eq:barbourlike} and in \eqref{eq:barbour} would be almost identical, and the difference between those would be controlled by the remainder terms $(R_{j,1,n})$ and $(R_{j,\omega,n})$. If those remainder terms are small then we could exploit the asymptotic normality of $V_n$ to obtain the asymptotic normality of $W_n$. We show that such a sequence exists when $|I_n|$ is large.

% \textcolor{blue}{As mentioned earlier, for ease of notation, in the proof of \cref{THM:BARBOURLIKE} we will drop the dependence on $n$ in our notation and write $W$, $N(\cdot)$, $\sigma$, $R_{j,1}$, $X_i$, $I$ and $R_{j,\omega}$ for respectively $W_n$, $N_n(\cdot)$, $\sigma_n$, $R_{j,1,n}$, $X^{\scalebox{0.6}{$(n)$}}_{i}$, $I_n$ and $R_{j,\omega,n}$}

% \textcolor{red}{I adapted this to the setting of $n$}
\begin{lemma}[Cumulant Matching]\label{THM:EXISTENCEXI}
  Let $p\geq 1$ and $k:=\lceil p\rceil$. If $p>1$, let $\bigl(u_{j}^{\scalebox{0.6}{$(n)$}}\bigr)_{j=1}^{k-1}$ be a sequence of real numbers. Suppose that for any $j=1,\cdots, k-1$, we have $u_{j}^{\scalebox{0.6}{$(n)$}}\to 0$ as $n\to \infty$. Then there exist constants $C_{p}, C_{p}'$ only depending on $p$ and a positive value $N>0$ (that might depend on $\bigl(u_{j}^{\scalebox{0.6}{$(n)$}}\bigr)$~) such that for any $n >N$, there exists $q_{n}\in \mathbb{N}_{+}$ and a random variable $\xi^{\scalebox{0.6}{$(n)$}}$ such that
  \begin{lemmaenum}
    \item \label{itm:match12} $\mathbb{E} [\xi^{\scalebox{0.6}{$(n)$}}]=0$,\quad $\mathbb{E} [(\xi^{\scalebox{0.6}{$(n)$}})^{2}]=1$;
    \item \label{itm:match3more} $\kappa_{j+2}(\xi^{\scalebox{0.6}{$(n)$}})=q_{n}^{j/2}u_{j}^{\scalebox{0.6}{$(n)$}}$ for $j=1,\cdots, k-1$;
    \item \label{itm:boundedaway} Either $\max_{1\leq j\leq k-1}\bigl\lvert \kappa_{j+2}(\xi^{\scalebox{0.6}{$(n)$}}) \bigr\rvert=0$ or $\max_{1\leq j\leq k-1}\bigl\lvert \kappa_{j+2}(\xi^{\scalebox{0.6}{$(n)$}}) \bigr\rvert\geq C_{p}>0$;
    \item \label{itm:momentbound} $\mathbb{E} [\lvert \xi^{\scalebox{0.6}{$(n)$}}\rvert^{p+2}]\leq C_{p}'$.
  \end{lemmaenum}
  Furthermore, $q_{n}$ can be chosen to be such that $q_{n}\to\infty$ as $\lvert I \rvert\to\infty$.
\end{lemma}
We note that the condition that $u_{j}^{\scalebox{0.6}{$(n)$}}\to 0$ as $n\to\infty$ is crucial. \cref{THM:EXISTENCEXI} is an asymptotic statement in the sense that for a given $n\leq N$, $q_{n}$ and $\xi^{\scalebox{0.6}{$(n)$}}$ might not exist.

%\textcolor{blue}{need to rephrase this parapgrah it is not gramatically correct}
Intuitively, \cref{itm:match12} and \cref{itm:match3more} determines the cumulants of $\xi^{\scalebox{0.6}{$(n)$}}$ and relates them to the cumulants of $W_n$. \cref{itm:boundedaway} requires that the maximum $$\max_{1\leq j\leq k}\bigl\lvert \kappa_{j+2}(\xi^{\scalebox{0.6}{$(n)$}}) \bigr\rvert$$ is either $0$ or bounded away from $0$ as $n$ grows. And \cref{itm:momentbound} indicates that the ($p$+$2$)-th absolute moment is upper-bounded.

\addtocontents{toc}{\SkipTocEntry}
\subsection{Proof of Theorem \ref{THM:LOCALWP}}
The proof of \cref{THM:LOCALWP} works in three stages:
\begin{enumerate}
  \item Using \cref{THM:EXISTENCEXI} we find a sequence of \emph{i.i.d.} random variables $\bigl(\xi^{\scalebox{0.6}{$(n)$}}_{\ell}\bigr)_{\ell}$ and a sample size $q_{n}$ such that the first $k$+$1$ cumulants of $W_n$ match the first $k$+$1$ cumulants of $V_n:=q_{n}^{-1/2}\sum_{i=1}^{q_{n}}\xi^{\scalebox{0.6}{$(n)$}}_i$;
  \item Using \cref{thm:lemzolo} we remark that we can bound the Wasserstein distance between the distributions of $W_n$ and an empirical average, $V_{n}$, of \emph{i.i.d.} observations in terms of $\bigl\lvert \mathbb{E} [h(W_n)]-\mathbb{E} [h(V_n)] \bigr\rvert$ for a large class of functions $h$. We do so by exploiting \cref{thm:barbour,THM:BARBOURLIKE};
  \item We remark that \cref{thm:lemiidwp} provides us with the bound on the Wasserstein distance between the distribution of $V_{n}$ and the standard normal.
\end{enumerate} Then \cref{THM:LOCALWP} follows from the triangle inequality of the Wasserstein metric: $$\mathcal{W}_p(W_n,\mathcal{N}(0,1))\le \mathcal{W}_p(\mathcal{L}(W_n),\mathcal{L}(V_n))+\mathcal{W}_p(\mathcal{L}(V_n),\mathcal{N}(0,1)).$$

\begin{proof}[Proof of \cref{THM:LOCALWP}]
  Without loss of generality, we assume $\sigma_n= 1$ and denote $W_n:=\sum_{i\in I_{n}}X^{\scalebox{0.6}{$(n)$}}_{i}$.

  Firstly, we remark that according to \cref{thm:corocumubd}, for all $1\le j\le k-1$ we have $\bigl\lvert \kappa_{j+2}(W_n) \bigr\rvert\lesssim R_{j,1,n}$. Moreover, by assumption we have  $R_{j,1,n}\to 0$ as $n\to \infty$. Therefore, $\bigl\lvert \kappa_{j+2}(W_n) \bigr\rvert \to 0$ as $n\to \infty$ and the assumptions of \cref{THM:EXISTENCEXI} hold, which implies that there exist constants $C_p$ and $C'_p$ 
  %2}(W_n)$ \rvert\lesssim R_{j,1}\to 0$ as $\lvert I \rvert\to \infty$ $\lvert \kappa_{j+2}(W_n) \rvert\lesssim R_{j,1}\to 0$ as $\lvert I \rvert\to \infty$ where $1\leq j\leq k-1$.
such  that for any $n$ large enough there are positive integers $(q_{n})$ and random variables $(\xi^{\scalebox{0.6}{$(n)$}})$ such that
  \begin{enumerate}[label=(\alph*), ref=(\alph*)]
    \item $\mathbb{E} [\xi^{\scalebox{0.6}{$(n)$}}]=0$,\quad $\mathbb{E} [(\xi^{\scalebox{0.6}{$(n)$}})^{2}]=1$;
          \item\label{bla} $\kappa_{j+2}(\xi^{\scalebox{0.6}{$(n)$}})=q_{n}^{j/2}\kappa_{j+2} (W_n)$ for $j=1,\cdots, k-1$;
    \item \label{itm:lowerbound} Either $\max_{1\leq j\leq k-1}\bigl\lvert \kappa_{j+2}(\xi^{\scalebox{0.6}{$(n)$}}) \bigr\rvert=0$ or $\max_{1\leq j\leq k-1}\bigl\lvert \kappa_{j+2}(\xi^{\scalebox{0.6}{$(n)$}}) \bigr\rvert\geq C_{p}>0$;
          \item\label{itm:uppermombound} $\mathbb{E} [\lvert \xi^{\scalebox{0.6}{$(n)$}}\rvert^{p+2}]\leq C_{p}'$.
  \end{enumerate}
  Furthermore, we know that $(q_{n})$ satisfies that $q_{n}\to\infty$ as $n\to\infty$.

  As presented in the proof sketch we will use this to bound the distance between the distribution of $W_n$ to the one of an empirical average of at least $q_{n}$ \emph{i.i.d.} random variables. Note that when $\max_{1\le j\le k-1}|\kappa_{j+2}(\xi^{\scalebox{0.6}{$(n)$}})|>0$ then we can obtain (by combining \cref{bla,itm:lowerbound} ) a lower bound on $q_{n}$ which will be crucial in our arguments as it will allow us to control the distance between this empirical average and its normal limit. When $\kappa_{3}(W_n)=\cdots=\kappa_{k+1}(W_n)=0$, such a lower bound on $q_{n}$ cannot be obtained in a similar way. Thus, we introduce an alternative sequence $(\widetilde q_{n})$ by setting $\widetilde{q}_{n}:=\lvert I_n \rvert^{2(p+1)/p}\vee q_{n}$ if $\kappa_{3}(W_n)=\cdots=\kappa_{k+1}(W_n)=0$, and $\widetilde{q}_{n}:=q_{n}$ otherwise. We remark that the sequence $(\widetilde{q}_{n})$ still respects $\widetilde{q}_{n}\to\infty$ as $n\to \infty$.

  Let $\xi_{1}^{\scalebox{0.6}{$(n)$}},\cdots,\xi_{\widetilde{q}_{n}}^{\scalebox{0.6}{$(n)$}}$ be \emph{i.i.d.} copies of $\xi^{\scalebox{0.6}{$(n)$}}$. Define $V_n:=\widetilde{q}_{n}^{-1/2}\sum_{i=1}^{\widetilde{q}_{n}}\xi_{i}^{\scalebox{0.6}{$(n)$}}$.% It will be seen later in the proof that the choice of $\widetilde{q}_{n}$ and $\xi_{1}^{\scalebox{0.6}{$(n)$}},\cdots,\xi_{\widetilde{q}_{n}}^{\scalebox{0.6}{$(n)$}}$ ensures that the specific values of $\widetilde{q}^{\scalebox{0.6}{$(I)$}}$ do not affect the rate of convergence.
  By construction, for any $j\in\mathbb{N}_{+}$ such that $j\leq k-1=\lceil p\rceil -1$ we have
  \begin{equation*}
    \kappa_{j+2}(V_n)\overset{(*)}{=}\widetilde{q}_{n}^{-(j+2)/2} \sum_{i=1}^{\widetilde{q}_{n}}\kappa_{j+2}(\xi_{i}^{\scalebox{0.6}{$(n)$}})=\widetilde{q}_{n}^{-j/2}\kappa_{j+2}(\xi^{\scalebox{0.6}{$(n)$}})=\kappa_{j+2}(W_n).
  \end{equation*}
 {Here in $(*)$ we have used the fact that cumulants are cumulative for independent random variables, which is directly implied by their definition. For more details on this, please refer to \cite{lukacs1970characteristic}.}

  Thus, by \cref{thm:barbour} and \cref{THM:BARBOURLIKE}, for any $h\in\Lambda_{p}$ we have
  \begin{equation}\label{doute}
    \bigl\lvert \mathbb{E} [h(W_n)]-\mathbb{E} [h(V_n)] \bigr\rvert \lesssim \sum_{j=1}^{k-1}R _{j,1,n}^{p/j}+\sum_{j=1}^{k}R _{j,\omega ,n}^{p/(j+\omega -1)}+\widetilde{q}_{n}^{-(p+2)/2}\sum_{i=1}^{\widetilde{q}_{n}}\mathbb{E} \bigl[\bigl\lvert \xi_{i}^{\scalebox{0.6}{$(n)$}} \bigr\rvert^{p+2}\bigr].
  \end{equation}
  To be able to have this upper bound not depend on $\xi_{i}^{\scalebox{0.6}{$(n)$}}$ we will upper-bound 
  $$\widetilde{q}_{n}^{-(p+2)/2}\sum_{i=1}^{\widetilde{q}_{n}}\mathbb{E} [\lvert \xi_{i}^{\scalebox{0.6}{$(n)$}} \rvert^{p+2}]$$ 
  in terms of the remainders $(R_{j,1,n})$ and $(R_{j,\omega,n})$. To do so we use the lower bounds on $(\widetilde{q}_{n})$ implied by the specific form we chose.

  If $\max_{1\leq j\leq k-1}\bigl\lvert \kappa_{j+2}(W_n)\bigr\rvert>0$, \cref{itm:lowerbound} implies that
  \begin{equation*}
    C_{p}\leq \max_{1\leq j\leq k-1}\bigl\lvert \kappa_{j+2}(\xi^{\scalebox{0.6}{$(n)$}}) \bigr\rvert\overset{(*)}{=}\max_{1\leq j\leq k-1}\bigl\{\widetilde{q}_{n}^{j/2}\bigl\lvert \kappa_{j+2}(W_n)\bigr\rvert\bigr\}\overset{(**)}{\lesssim} \max_{1\leq j\leq k-1}\bigl\{\widetilde{q}_{n}^{j/2}R_{j,1,n}\bigr\}.
  \end{equation*} where to get $(*)$ we used \cref{bla} and to get $(**)$ we used \cref{thm:corocumubd}.
  Thus, the following holds

  % \textcolor{red}{[Note!!! You pointed an issue here but your previous change was incorrect. Please check it now.]}
  \begin{equation*}
    \widetilde{q}_{n}^{-p/2}=(\widetilde{q}_{n}^{-j_{0}/2})^{p/j_{0}}\lesssim R_{j_{0},1,n}^{p/j_{0}}\leq \sum_{j=1}^{k-1} R_{j,1,n}^{p/j},
  \end{equation*}
  where $j_{0}$ is the integer satisfying that $\bigl\lvert \kappa_{j_{0}+2}(\xi^{\scalebox{0.6}{$(n)$}}) \bigr\rvert=\max_{1\leq j\leq k-1}\bigl\lvert \kappa_{j+2}(\xi^{\scalebox{0.6}{$(n)$}}) \bigr\rvert$.

  On the other hand, if $\kappa_{j+2}(W_n)=0$ for all $1\leq j\leq k-1$, then by definitions we have $\widetilde{q}_{n}\geq \lvert I_n \rvert^{2(p+1)/p}$, and therefore, $ \widetilde{q}_{n}^{-p/2}\le\lvert I_n \rvert^{-(p+1)}$. %Therefore, we know that $\widetilde{q}_{n}\geq \lvert I \rvert^{2(p+1)/p}$. 
  Moreover, by Hölder's inequality we know that the following holds
  \begin{equation}\label{nelson}
    \sum_{i\in I_n}\mathbb{E} \bigl[\bigl\lvert X^{\scalebox{0.6}{$(n)$}}_{i} \bigr\rvert^{2}\bigr]\leq \lvert I_n\rvert^{p/(p+2)}\Bigl(\sum_{i\in I_n}\mathbb{E} \bigl[\bigl\lvert X^{\scalebox{0.6}{$(n)$}}_{i} \bigr\rvert^{p+2}\bigr]\Bigr)^{2/(p+2)}.
  \end{equation}
  and      \begin{equation}\label{luna}
    \Bigl(\sum_{i\in I_n}X^{\scalebox{0.6}{$(n)$}}_{i}\Bigr)^{2}\leq \lvert I_n \rvert \sum_{i\in I_n}\bigl\lvert X^{\scalebox{0.6}{$(n)$}}_{i} \bigr\rvert^{2}.
  \end{equation}
  Since $\mathbb{E} \Bigl[\Bigl(\sum_{i\in I_n}X^{\scalebox{0.6}{$(n)$}}_{i}\Bigr)^{2}\Bigr]=\sigma_n^{2}=1$, we have
  \begin{align*}
    \widetilde{q}_{n}^{-p/2}
    \le &\lvert I_n \rvert^{-(p+1)}\Bigl(\mathbb{E} \Bigl[\Bigl(\sum_{i\in I}X^{\scalebox{0.6}{$(n)$}}_{i}\Bigr)^{2}\Bigr]\Bigr)^{(p+2)/2}\\
    \overset{(*)}{\leq} &\lvert I_n \rvert^{-p/2}\Bigl(\sum_{i\in I}\mathbb{E} \bigl[\bigl\lvert X^{\scalebox{0.6}{$(n)$}}_{i} \bigr\rvert^{2}\bigr]\Bigr)^{(p+2)/2}\\
    \overset{(**)}{\leq}& \sum_{i\in I_n}\mathbb{E} \bigl[\bigl\lvert X^{\scalebox{0.6}{$(n)$}}_{i} \bigr\rvert^{p+2}\bigr]\leq R_{k,\omega,n },
  \end{align*} where to obtain $(*)$ we used \eqref{luna} and to obtain $(**)$ we used \eqref{nelson}.

  Thus, using \cref{itm:uppermombound} and the fact that $\xi_{1}^{\scalebox{0.6}{$(n)$}},\cdots,\xi_{\widetilde{q}_{n}}^{\scalebox{0.6}{$(n)$}}$ are \emph{i.i.d.}, we obtain
  \begin{equation}\label{eq:compareconnect}
    \widetilde{q}_{n}^{-(p+2)/2}\sum_{i=1}^{\widetilde{q}_{n}}\mathbb{E} \bigl[\bigl\lvert \xi_{i}^{\scalebox{0.6}{$(n)$}} \bigr\rvert^{p+2}\bigr]\leq C_{p}'\widetilde{q}_{n}^{-p/2}\lesssim \sum_{j=1}^{k-1}R _{j,1,n}^{p/j}+\sum_{j=1}^{k}R _{j,\omega ,n}^{p/(j+\omega -1)}.
  \end{equation}

  Therefore, by combining this with \eqref{doute} we obtain that there is a constant $K>0$ that does not depend on $h$ such that
  \begin{equation*}
    \bigl\lvert \mathbb{E} [h(W_n)]-\mathbb{E} [h(V_n)] \bigr\rvert \le K\Bigl( \sum_{j=1}^{k}R _{j,1,n}^{p/j}+\sum_{j=1}^{k+1}R _{j,\omega ,n}^{p/(j+\omega -1)}\Bigr).
  \end{equation*}
  By taking supremum over $h\in \Lambda_{p}$ and by \cref{thm:lemzolo} we obtain that
  \begin{equation*}
    \mathcal{W}_{p}(\mathcal{L}(W_n),\mathcal{L}(V_n))\lesssim \sup_{h\in\Lambda_{p}}\bigl\lvert \mathbb{E}[h(W_n)]-\mathbb{E} [h(V_n)] \bigr\rvert^{1/p}\lesssim \sum_{j=1}^{k-1}R _{j,1,n}^{1/j}+\sum_{j=1}^{k}R _{j,\omega ,n}^{1/(j+\omega -1)}.
  \end{equation*}
  Moreover, by combining  \cref{thm:lemiidwp} and \eqref{eq:compareconnect} we have
  \begin{equation*}
    \mathcal{W}_{p}(\mathcal{L}(V_n),\mathcal{N}(0,1))\lesssim \Bigl(\widetilde{q}_{n}^{-(p+2)/2}\sum_{i=1}^{\widetilde{q}_{n}}\mathbb{E} \bigl[\bigl\lvert \xi_{i}^{\scalebox{0.6}{$(n)$}} \bigr\rvert^{p+2}\bigr]\Bigr)^{1/p}\lesssim \sum_{j=1}^{k-1}R_{j,1,n}^{1/j}+\sum_{j=1}^{k}R _{j,\omega,n }^{1/(j+\omega -1)}.
  \end{equation*}
  Therefore, as the Wasserstein distance $\mathcal{W}_p$ satisfies the triangle inequality we conclude that
  \begin{align*}
    \mathcal{W}_{p}(\mathcal{L}(W_n),\mathcal{N}(0,1))
    \leq &\mathcal{W}_{p}(\mathcal{L}(W_n),\mathcal{L}(V_n))+\mathcal{W}_{p}(\mathcal{L}(V_n),\mathcal{N}(0,1))\\
    \lesssim &\sum_{j=1}^{k-1}R _{j,1,n}^{1/j}+\sum_{j=1}^{k}R _{j,\omega ,n}^{1/(j+\omega -1)}.
  \end{align*}
\end{proof}

\section{Proof of Lemma~\ref{THM:BARBOURLIKE}}\label{sec:lemma1}
For ease of notation, when there is no ambiguity we will drop the dependence on $n$ in our notation and write $W$, $N(\cdot)$, $\sigma$, $X_i$, $I$ and $R_{j,\omega}$ for respectively $W_n$, $N_n(\cdot)$, $\sigma_n$, $X^{\scalebox{0.6}{$(n)$}}_{i}$, $I_n$ and $R_{j,\omega,n}$.

\addtocontents{toc}{\SkipTocEntry}
\subsection{Example and Roadmap}\label{rd_mp}

Given the form of expression in \cref{THM:BARBOURLIKE}, it is natural to consider performing induction on $\lceil p\rceil$. In fact, \cite{barbour1986asymptotic} used a similar induction idea to prove \cref{thm:barbour}, the analogous result to \cref{THM:BARBOURLIKE} for independent variables. As \cite{fang2019wasserstein} suggested, the key of each inductive step is the following expansion of $\mathbb{E} [Wf(W)]$.

\begin{proposition}[Expansion of {$\mathbb{E}[Wf(W)]$}]\label{THM:WFWEXPANSION}
  Denote by $\kappa_{j}(W)$ the $j$-th cumulant of $W$. Given $k\in \mathbb{N}_{+}$ and real number $\omega\in (0,1]$, for any $f\in \mathcal{C}^{k,\omega }(\mathbb{R})$, we have
  \begin{equation}\label{eq:wfwexpansion}
    \mathbb{E} [Wf(W)]=\sum_{j=1}^{k}\frac{\kappa_{j+1}(W)}{j !}\mathbb{E} [\partial^{j} f(W)]+\mathcal{O}(\lvert f \rvert_{k,\omega }R_{k,\omega }).
  \end{equation}
\end{proposition}

The case $k=\omega=1$ is a well-known result in the literature of Stein's method (for example see \cite{barbour1989central,ross2011fundamentals}). The case $k=2, \omega=1$ was first proven by \cite{fang2019wasserstein}, and they also conjectured that it was true for any positive integer $k$ with $\omega=1$. Inspired by \cite{fang2019wasserstein}'s method, we confirm that this conjecture is correct by proving \cref{THM:WFWEXPANSION}.

To help better understand the intuition behind our proof for the general settings, let's first consider the simplest case with $k=\omega=1$. Given a positive integer $m$, suppose that $(X_{i})_{i=1}^{n}$ is an $m$-dependent random sequence (the special case of $d=1$ in \cref{thm:defmdepfield}). We let $W:=\sum_{i=1}^{n}X_{i}$ and require that $\mathbb{E} [X_{1}]=0$ and $\mathbb{E} [W^{2}]=1$. For simplicity, we further assume $f\in \mathcal{C}^{2}(\mathbb{R})\cap \mathcal{C}^{1,1}(\mathbb{R})$ meaning that $f''$ is a continuous and bounded function.% We remark that \eqref{eq:wfwexpansion} reduces to the upper bound:
% \begin{equation*}
%   \bigl\lvert \mathbb{E} [Wf(W)]-\mathbb{E} [f'(W)] \bigr\rvert\lesssim \lVert f'' \rVert R_{1,1}.
% \end{equation*}

{For any indexes $i,j\in [n]$ (by convention $[n]:=\{ 1,2,\cdots,n \}$), we write $$N(i)=\{ \ell\in [n]: \lvert \ell-i \rvert\leq m \},\quad N(i,j):=\{ \ell \in [n]:\lvert \ell-i \rvert\leq m\text{ or }\lvert \ell-j \rvert\leq m \}.$$}%
Denote $W_{i,m}:=\sum_{j\notin N(i)}X_{j}$ and $W_{i,j,m}:=\sum_{\ell\notin N(i,j)}X_{\ell}$. The idea is that for each $i$, we split $W$ into two parts, $W_{i,m}$ and $W-W_{i,m}$. The former is independent of $X_{i}$ while the latter is the sum of $X_{j}$'s in the neighborhood of $X_{i}$ and will converge to $0$ when $n$ grows to $\infty$. Thus, we perform the Taylor expansion for $f(W)$ around $W_{i,m}$.

We have
\begin{align}\label{eq:split0}
  &\quad \mathbb{E}\bigl[ Wf (W)- f'(W) \bigr]\\
  = &
  \sum_{i=1}^{n}\mathbb{E} \bigl[X_{i}\bigl(f (W)-f (W_{i,m}) - f'(W_{i,m})(W-W_{i,m})\bigr)\bigr] \nonumber                            \\
    & \  + \sum_{i=1}^{n}\mathbb{E} [X_{i}f (W_{i,m})]
  +   \sum_{i=1}^{n}\mathbb{E} \bigl[X_{i}(W-W_{i,m})f'(W_{i,m})\bigr]-\mathbb{E} [f'(W)] \nonumber                                     \\
  = &
  \sum_{i=1}^{n}\mathbb{E} \bigl[X_{i}\bigl(f (W)-f (W_{i,m}) - f'(W_{i,m})(W-W_{i,m})\bigr)\bigr] \nonumber                            \\
    & \  + \sum_{i=1}^{n}\ \mathbb{E} [X_{i}]\ \mathbb{E} [f (W_{i,m})]
  +   \sum_{i=1}^{n}\sum_{j\in N(i)}\mathbb{E} \bigl[X_{i}X_{j}f'(W_{i,m})\bigr]-\mathbb{E} [f'(W)] \nonumber                           \\
  = &
  \sum_{i=1}^{n}\mathbb{E} \bigl[X_{i}\bigl(f (W)-f (W_{i,m}) - f'(W_{i,m})(W-W_{i,m})\bigr)\bigr] \nonumber                            \\
    & \    +   \Bigl(\sum_{i=1}^{n}\sum_{j\in N(i)}\mathbb{E} \bigl[X_{i}X_{j}f'(W_{i,m})\bigr]-\mathbb{E} [f'(W)] \Bigr)=:E_{1}+E_{2}.\nonumber
\end{align}

By assumption, $\lVert f''\rVert$ is bounded and we have
\begin{align}\label{eq:bounde1}
  \lvert E_{1} \rvert= & \biggl\lvert
  \sum_{i=1}^{n}\mathbb{E} \bigl[X_{i}\bigl(f (W)-f (W_{i,m})
  - f'(W_{i,m})(W-W_{i,m})\bigr)\bigr]\biggr\rvert       \\
  \leq                 & \frac{\lVert f''\rVert}{2}\sum_{i=1}^{n}
  \mathbb{E}\bigl[\bigl\lvert X_{i}(W-W_{i,m})^{2}\bigr\rvert\bigr]
  =     \frac{\lVert f''\rVert}{2}\sum_{i=1}^{n}
  \mathbb{E}\Bigl[\lvert X_{i}\rvert\ \Bigl(
  \sum_{j\in N(i)}X_{j}\Bigr)^{2} \Bigr]\nonumber                 \\
  =                    & \frac{\lVert f''\rVert}{2}\sum_{i=1}^{n}
  \sum_{j\in N(i)}\sum_{\ell\in N(i)}
  \mathbb{E}[\lvert X_{i}X_{j}X_{\ell}\rvert]\leq \frac{\lVert f''\rVert}{2}\sum_{i=1}^{n}
  \sum_{j\in N(i)}\sum_{\ell\in N(i,j)}
  \mathbb{E}[\lvert X_{i}X_{j}X_{\ell}\rvert].\nonumber
\end{align}

Now we bound $E_{2}$.
\begin{align}\label{eq:split20}
  E_{2} =             & \sum_{i=1}^{n}\sum_{j\in N(i)}\mathbb{E} \bigl[X_{i}X_{j}f'(W_{i,m})\bigr]-\mathbb{E} [f'(W)]                                                                                                     \\
  =                   & \sum_{i=1}^{n}\sum_{j\in N(i)}\mathbb{E}\bigl[X_{i}X_{j}\bigl(f'(W_{i,m})-f'(W_{i,j,m})\bigr)\bigr]\nonumber\\*
  &\ +\sum_{i=1}^{n}\sum_{j\in N(i)} \mathbb{E}\bigl[X_{i}X_{j}f'(W_{i,j,m})\bigr] -\mathbb{E} [f'(W)]\nonumber      \\
  \overset{(*)} {=}   & \sum_{i=1}^{n}\sum_{j\in N(i)} \mathbb{E}\bigl[X_{i}X_{j}\bigl(f'(W_{i,m})-f'(W_{i,j,m})\bigr)\bigr]\nonumber\\*
  &\ +\sum_{i=1}^{n}\sum_{j\in N(i)}\mathbb{E}[X_{i}X_{j}]\ \mathbb{E} [f'(W_{i,j,m})] -\mathbb{E} [f'(W)]\nonumber \\
  \overset{(**)}{  =} & \sum_{i=1}^{n}\sum_{j\in N(i)}\mathbb{E}\bigl[X_{i}X_{j}\bigl(f'(W_{i,m})-f'(W_{i,j,m})\bigr)\bigr]\nonumber\\*
  &\ + \sum_{i=1}^{n}\sum_{j\in N(i)}\mathbb{E}[X_{i}X_{j}]\ \mathbb{E} \bigl[f'(W_{i,j,m})-f'(W)\bigr] \nonumber             \\
  =                 & (t_1)+(t_2),\nonumber
\end{align}
where to obtain $(*)$ we have used the fact that $W_{i,j,m}$ is independent of $(X_{i},X_{j})$ in the second equation and to obtain $(**)$ we have assumed hat $\mathbb{E}(W^2)=1$.

The first term $(t_1)$, can be upper-bounded by the mean value theorem as
\begin{align*}
       & \biggl\lvert\sum_{i=1}^{n}\sum_{j\in N(i)}\mathbb{E}\bigl[X_{i}X_{j}\bigl(f'(W_{i,m})-f'(W_{i,j,m})\bigr)\bigr]\biggr\rvert   \\
  \leq & \sum_{i=1}^{n}\sum_{j\in N(i)}\lVert f'' \rVert \ \mathbb{E} \bigl[\bigl\lvert X_{i}X_{j}(W_{i,m}-W_{i,j,m})\bigr\rvert\bigr] \\
  \leq & \lVert f'' \rVert\sum_{i=1}^{n}
  \sum_{j\in N(i)}\sum_{\ell\in N(i,j)}
  \mathbb{E}[\lvert X_{i}X_{j}X_{\ell}\rvert].
\end{align*}

By another application of the mean-value theorem, we remark that the second term $(t_2)$, is controlled by
\begin{align*}
       & \biggl\lvert\sum_{i=1}^{n}\sum_{j\in N(i)}\mathbb{E}[X_{i}X_{j}]\ \mathbb{E} \bigl[f'(W_{i,j,m})-f'(W)\bigr]\biggr\rvert                        \\
  \leq & \sum_{i=1}^{n}\sum_{j\in N(i)}\lVert f'' \rVert\ \mathbb{E}[\lvert X_{i}X_{j}\rvert]\ \mathbb{E} \bigl[\bigl\lvert W_{i,j,m}-W\bigr\rvert\bigr] \\
  \leq & \lVert f'' \rVert\sum_{i=1}^{n}\sum_{j\in N(i)}\sum_{\ell\in N(i,j)}\mathbb{E}[\lvert X_{i}X_{j}\rvert]\ \mathbb{E} [\lvert X_{\ell} \rvert].
\end{align*}

Thus,
\begin{equation*}
  \begin{aligned}
  &\bigl\lvert \mathbb{E} \bigl[Wf(W)-f'(W)\bigr] \bigr\rvert\\
  \leq &\lVert f'' \rVert \sum_{i=1}^{n}\sum_{j\in N(i)}\sum_{\ell\in N(i,j)}\biggl(\frac{3}{2}\mathbb{E} [\lvert X_{i}X_{j}X_{\ell} \rvert]+\mathbb{E} [\lvert X_{i}X_{j} \rvert]\ \mathbb{E} [\lvert X_{\ell} \rvert]\biggr)\\
  \leq &\frac{3\lVert f'' \rVert}{2}R_{1,1}.
  \end{aligned}
\end{equation*}
This gives us a bound that matches with \eqref{eq:wfwexpansion}.

For $k\geq 2$, we would like to carry out the expansion in the same spirit. However, it would be too tedious to write out every sum in the process. Thus, in \cref{sec:pflocalnotation}, we introduce the terms called $\mathcal{S}$-sums, $\mathcal{T}$-sums, and $\mathcal{R}$-sums, which serve as useful tools in tracking different quantities when we approximate $\mathbb{E} [f'(W)-Wf(W)]$ with respect to locally dependent random variables. Instead of performing the expansion to get \eqref{eq:wfwexpansion} for $\mathbb{E} [Wf(W)]$, we first do it for any $\mathcal{T}$-sum and use induction to prove a more general result for the existence of such expansions (see \cref{thm:grandexpand}). In the general situation of $\mathcal{T}$-sums, the cumulants are replaced by other constants that only depend on the specific $\mathcal{T}$-sum in consideration and the joint distribution of $(X_{i})_{i\in I}$. Finally, we prove that in particular, those constants for $\mathbb{E} [Wf(W)]$ are precisely the cumulants of $W$. This will be done by direct calculation when $f$ is a polynomial and then extended to more general $f$'s by applying \cref{thm:uniqueexp}.

\addtocontents{toc}{\SkipTocEntry}
\subsection{Notations and Definitions}\label{sec:pflocalnotation}

% Given integers $m\geq 2$ and $\ell\geq 1$, let $\eta_{1},\eta_{2},\cdots,\eta_{\ell}$ be a composition of the integer $m$, i.e., $\eta_{1:\ell}$ is a sequence of positive integers such that $\eta_{1}+\eta_{2}+\cdots+\eta_{\ell}=m$. For any random variables $Y_{1},\cdots,Y_{m}$, we define an order-$m$ \textbf{compositional expectation} with respect to $\eta_{1:\ell}$ as
% \begin{equation*}
%   [\eta_{1},\cdots,\eta_{\ell}]\triangleright (Y_{1},\cdots,Y_{m}):=\mathbb{E} [Y_{1}\cdots Y_{\eta_{1}}]\ \mathbb{E} [Y_{\eta_{1}+1}\cdots Y_{\eta_{1}+\eta_{2}}]\ \cdots \ \mathbb{E} [Y_{\eta_{1}+\cdots+\eta_{\ell-1}+1}\cdots Y_{m}].
% \end{equation*}
% Intuitively, $\eta_{1:\ell}$ determines how it decomposes into different expectations. 

As in \cref{SEC:LOCAL}, given an integer $k\geq 1$, suppose $(X_{i})_{i\in I}$ is a class of mean zero random variables indexed by $I$ that satisfy the local dependence assumptions {[LD-$1$]} to {[LD-$k$]}. Without loss of generality, we always assume that $\sigma^{2}:=\operatorname{Var}\left(\sum_{i\in I}X_{i}\right)=1$. We denote $W:=\sigma^{-1}\sum_{i\in I}X_{i}=\sum_{i\in I}X_{i}$.

%\textcolor{blue}{improve writing on first point}
% \begin{itemize}
%   \item $N_{i.}[0]:=\emptyset$, and shorthand $N(i_{1},\cdots, i_{j})$ as $N_{i.}[j]$ (for a definition see \textbf{[LD-$1$]} to \textbf{[LD-$k$]});
%   \item $W_{i.}[0]:=W:=\sigma^{-1}\sum_{i\in I}X_{i}=\sum_{i\in I}X_{i}$;
%   \item $W_{i.}[j]:=\sigma^{-1}\sum_{i\in I\setminus N_{i.}[j]}X_{i}=\sum_{i\in I\setminus N_{i.}[j]}X_{i}$.% for $1\leq j\leq k$.
% \end{itemize}
\subsubsection*{$\mathcal{S}$-sums}

  Fix $k\in\mathbb{N}_{+}$ and $t_{1},\cdots,t_{k}\in \mathbb{Z}$ be integers such that $\lvert t_{j}\rvert \leq j-1$ for any $j\in [k]$. We set  $t_{1}=0$. Let $z=\bigl\lvert\{ j:t_{j}>0 \}\bigr\rvert$ be the number of indexes $j$ for which $t_j$ is strictly positive.  If $z\geq 1$, we write $\{ j:t_{j}>0 \}=\{ q_{1},\cdots,q_{z} \}$. Without loss of generality, we suppose that the sequence $2\leq q_{1}<\cdots<q_{z}\leq k$ is increasing. We further let $q_{0}:=1$ and $q_{z+1}:=k+1$. We define an order-$k$ \textbf{$\mathcal{S}$-sum} with respect to the sequence $t_{1:k}$ as
\begin{align}\label{eq:defcomp1}
    &\mathcal{S} [t_{1},\cdots,t_{k}]
  :=  \sum_{i_{1}\in N_{1}}\sum_{i_{2}\in N_{2}}\cdots\sum_{i_{k}\in N_{k}}\ [q_{1}-q_{0},\cdots,q_{z+1}-q_{z}]\triangleright \bigl(X_{i_{1}},\cdots,X_{i_{k}}\bigr)                                                                                                                                  \\
  &\ = \sum_{i_{1}\in N_{1}}\sum_{i_{2}\in N_{2}}\cdots\sum_{i_{k}\in N_{k}}\ \mathbb{E} \bigl[X_{i_{q_{0}}}\cdots X_{i_{q_{1}-1}}\bigr]\ \mathbb{E} \bigl[X_{i_{q_{1}}}\cdots X_{i_{q_{2}-1}}\bigr]\ \cdots\ \mathbb{E} \bigl[X_{i_{q_{z}}}\cdots X_{i_{q_{z+1}-1}}\bigr],\nonumber 
\end{align}
where $N_{1}:=I$, and for $j\in\mathbb{N}_{+}$ such that $j\geq 2$, we let
$$
N_{j}:=\begin{cases} N (i_{1:\lvert t_{j} \rvert})=N(i_{1},\cdots,i_{\lvert t_{j}\rvert})& \text{ if }t_{j}\neq 0\\ 
  \emptyset &\text{ if }t_{j}=0
\end{cases}.
$$
Note that $N_{j}$ depends on $t_{j}$ and the sequence $i_{1:(j-1)}$. For ease of notation, we do not explicitly write out the dependencies on $i_{1:(j-1)}$ when there is no ambiguity. Further note that if any $t_j$, that is not $t_1$, is null then $N_j=\emptyset$ therefore, the $\mathcal{S}$-sum $\mathcal{S}[t_1,\cdots,t_k]=0$.

By definition all $\mathcal{S}$-sums are deterministic quantities, the value of which only depends on $t_{1:k}$, and the joint distribution of $(X_{i})_{i\in I}$. We also remark that the signs of $t_{j}$'s determine how an $\mathcal{S}$-sum factorizes into different expectations.  Notably if $z=0$ (meaning that all the $t_j$ are negative) then the {$\mathcal{T}$-sum} is 
$$
\mathcal{T}_{f,s} [t_{1},\cdots,t_{k}]=\sum_{i_{1}\in N_{1}}\sum_{i_{2}\in N_{2}}\cdots\sum_{i_{k}\in N_{k}}\mathbb{E} \Bigl[X_{i_{1}}\cdots X_{i_{k}}\partial^{k-1} f\bigl(W_{i.}[k-s]\bigr)\Bigr].
$$  
Since by assumption, $X_{i}$'s are centered random variables, the $\mathcal{S}$-sum vanishes if $q_{j+1}=q_{j}+1$ for some $0\leq j\le z$: 
\begin{equation}
\begin{aligned}\label{j}
    &\mathcal{S} [t_{1},\cdots,t_{k}]
=\sum_{i_{1}\in N_{1}}\sum_{i_{2}\in N_{2}}\cdots\sum_{i_{k}\in N_{k}}\ \mathbb{E} \bigl[X_{i_{q_{0}}}\cdots X_{i_{q_{1}-1}}\bigr]\ \cdot\\*
&\qquad\mathbb{E} \bigl[X_{i_{q_{1}}}\cdots X_{i_{q_{2}-1}}\bigr]\ \cdots \mathbb{E}[X_{i_{q_j}}]\ \cdots\ \mathbb{E} \bigl[X_{i_{q_{z}}}\cdots X_{i_{q_{z+1}-1}}\bigr]=0.
\end{aligned}
\end{equation}
Furthermore, the {absolute} value of $t_{j}$'s influences the range of running indexes. The bigger $\lvert t_{j}\rvert$ is the larger the set $N_{j}$ is. The largest possible index set for $i_{j+1}$ is $N(i_{1:(j-1)})$, which corresponds to the case $\lvert t_j\rvert=j-1$. On the other hand, if $t_{j}=0$, the sum is over an empty set and vanishes. 
In particular, if we require that the $\mathcal{S}$-sum is not always zero, then $t_{2}$ is always taken to be $-1$ and $i_{2}\in N(i_{1})$.

\subsubsection*{$\mathcal{T}$-sums}

For any function $f\in \mathcal{C}^{k-1}(\mathbb{R})$ and integer $s\in \mathbb{N}$ such that $s\leq k$, the order-$k$ \textbf{$\mathcal{T}$-sum}, with respect to the sequence $t_{1:k}$, is defined as 
  \begin{equation}\label{eq:defcomp2}
    \begin{aligned}
   & \mathcal{T}_{f,s} [t_{1},\cdots,t_{k}]
   := \\
   &\!\!\!\!\!\!  \sum_{i_{1}\in N_{1}}\sum_{i_{2}\in N_{2}}\cdots \sum_{i_{k}\in N_{k}}\ [q_{1}-q_{0},\cdots,q_{z+1}-q_{z}]\triangleright \bigl(X_{i_{1}},\cdots,X_{i_{k-1}},X_{i_{k}}\partial^{k-1}f\bigl(W_{i.}[k-s]\bigr)\bigr)\\
   = &
     \left\{\begin{aligned}
    &\sum_{i_{1}\in N_{1}}\sum_{i_{2}\in N_{2}}\cdots\sum_{i_{k}\in N_{k}}\mathbb{E} \Bigl[X_{i_{1}}\cdots X_{i_{k}}\partial^{k-1} f\bigl(W_{i.}[k-s]\bigr)\Bigr]&&\text{ if }z=0\\
    &\begin{aligned}
   &\sum_{i_{1}\in N_{1}}\sum_{i_{2}\in N_{2}}\cdots\sum_{i_{k}\in N_{k}}\mathbb{E} \bigl[X_{i_{q_{0}}}\cdots X_{i_{q_{1}-1}}\bigr]\ \cdots\ \mathbb{E} \bigl[X_{i_{q_{(z-1)}}}\cdots X_{i_{q_{z}-1}}\bigr]\ \cdot\\
   &\qquad\quad \mathbb{E} \bigl[X_{i_{q_{z}}}\cdots X_{i_{k}}\partial^{k-1}f\bigl(W_{i.}[k-s]\bigr)\bigr]
    \end{aligned}&&\text{ if }z\geq 1,
   \end{aligned}\right.
 \end{aligned}
\end{equation} 
%  and if $z>0$  as  \begin{align}\label{eq:defcomp2}
%     \mathcal{T}_{f,s} [t_{1},\cdots,t_{k}]
%   := &\sum_{i_{1}\in N_{1}}\sum_{i_{2}\in N_{2}}\cdots \sum_{i_{k}\in N_{k}}\ [q_{1}-q_{0},\cdots,q_{z+1}-q_{z}]\triangleright \bigl(X_{i_{1}},\cdots,X_{i_{k-1}},X_{i_{k}}\partial^{k-1}f\bigl(W_{i.}[k-s]\bigr)\bigr)\nonumber\\
%   = &
%   \sum_{i_{1}\in N_{1}}\sum_{i_{2}\in N_{2}}\cdots\sum_{i_{k}\in N_{k}}\mathbb{E} \bigl[X_{i_{q_{0}}}\cdots X_{i_{q_{1}-1}}\bigr]\ \cdots\ \mathbb{E} \bigl[X_{i_{q_{(z-1)}}}\cdots X_{i_{q_{z}-1}}\bigr]\ \mathbb{E} \bigl[X_{i_{q_{z}}}\cdots X_{i_{k}}\partial^{k-1}f\bigl(W_{i.}[k-s]\bigr)\bigr];
%  \end{align} 
 where $N_{1:k},z,q_{0:(z+1)}$ are defined as in the definition of $\mathcal{S}$-sums and $W_{i.}[j]$ is defined as
 \begin{equation*}
  W_{i.}[j]:=\begin{cases}
    W& \text{ if }j=0\\
  \sum_{i\in I\backslash N(i_{1:j})}X_{i}&\text{ if }1\leq j\leq k
  \end{cases}.
 \end{equation*}

 Note that the bigger $s$ is, the larger the set $I\backslash N(i_{1:(k-s)})$ is, which means that $W_{i.}[k-s]$ is the sum of more $X_{i}$'s. Again we remark that the values of {$\mathcal{T}$-sums} can depend on the values of $s$ and the sequences $t_{1:k}$. In particular, if $s=0$, then we have $W_{i.}[k-s]=W_{i.}[k]=\sum_{i\in I \backslash N (i_{1:k})}X_{i}$, which implies that $W_{i.}[k-s]$ is independent of $X_{i_{1}},\cdots,X_{i_{k}}$ by the assumption \textbf{[LD-$k$]}. Thus, we have
\begin{equation*}
  \mathbb{E} \bigl[X_{i_{q_{z}}}\cdots X_{i_{k}}\partial^{k-1}f\bigl(W_{i.}[k-s]\bigr)\bigr]=\mathbb{E}[X_{i_{q_{z}}}\cdots X_{i_{k}}]\ \mathbb{E} \bigl[\partial^{k-1}f\bigl(W_{i.}[k-s]\bigr)\bigr].
\end{equation*}
By definitions \eqref{eq:defcomp1} and \eqref{eq:defcomp2} we get
  \begin{equation}\label{eq:indepcomp}
    \mathcal{T}_{f,0} [t_{1},\cdots,t_{k}]=\mathcal{S}[t_{1},\cdots,t_{k}]\ \mathbb{E} [\partial^{k-1}f(W_{i.}[k])].
  \end{equation}
 This equation will be useful in our discussion later. In general if $z>0$ then
  \begin{align*}
    \mathcal{T}_{f,s} [t_{1},\cdots,t_{k}]
= &\mathcal{S}[t_{1},\cdots,t_{q_{z}-1}]\ \sum_{i_{q_{z}}\in N_{q_{z}}}\sum_{i_{q_{z}+1}\in N_{q_{z}+1}}\cdots\sum_{i_{k}\in N_{k}}\\*
   &\quad \mathbb{E} \bigl[X_{i_{q_{z}}}\cdots X_{i_{k}}\partial^{k-1}f\bigl(W_{i.}[k-s]\bigr)\bigr].
 \end{align*}

\subsubsection*{$\mathcal{R}$-sums}

For $k\geq 2$ and given a real number $\omega\in (0, 1]$, we further define an order-$k$ \textbf{$\mathcal{R}$-sum} with respect to the sequence $t_{1:k}$ as 
\begin{equation}\label{eq:defcomp3}
  \begin{aligned}
    &\mathcal{R}_{\omega} [t_{1},\cdots, t_{k}]:= \\
  &\!\!\!\!\!\!\sum_{i_{1}\in N_{1}}\sum_{i_{2}\in N_{2}}\cdots\sum_{i_{k-1}\in N_{k-1}}[q_{1}-q_{0},\cdots,q_{z+1}-q_{z}]\triangleright \Bigl(\lvert X_{i_{1}}\rvert,\cdots,\lvert X_{i_{k-1}}\rvert,\bigl(\sum_{i_{k}\in N_{k}}\lvert X_{i_{k}}\rvert\bigr)^{\omega }\Bigr)\\
  = & \left\{\begin{aligned}
    &\sum_{i_{1}\in N_{1}}\sum_{i_{2}\in N_{2}}\cdots\sum_{i_{k}\in N_{k}}\mathbb{E} \Bigl[X_{i_{1}}\cdots X_{i_{k-1}}\bigl(\sum_{i_{k}\in N_{k}}\lvert X_{i_{k}}\rvert\bigr)^{\omega }\Bigr]&&\text{ if }z=0\\
    &\begin{aligned}
   &\sum_{i_{1}\in N_{1}}\sum_{i_{2}\in N_{2}}\cdots\sum_{i_{k}\in N_{k}}\mathbb{E} \bigl[X_{i_{q_{0}}}\cdots X_{i_{q_{1}-1}}\bigr]\ \cdots\ \mathbb{E} \bigl[X_{i_{q_{(z-1)}}}\cdots X_{i_{q_{z}-1}}\bigr]\ \cdot\\
   &\qquad\mathbb{E} \biggl[X_{i_{q_{z}}}\cdots X_{i_{k-1}}\bigl(\sum_{i_{k}\in N_{k}}\lvert X_{i_{k}}\rvert\bigr)^{\omega }\biggr]
    \end{aligned}&&\text{ if }z\geq 1.
   \end{aligned}\right.
  \end{aligned}
\end{equation}
We again remark that if $z\ge 1$ then   
\begin{align*}
    \mathcal{R}_{\omega} [t_{1},\cdots, t_{k}]
=  & \mathcal{R}_{1}[t_{1},\cdots,t_{q_{z}-1}]\ \sum_{i_{q_{z}}\in N_{q_{z}}}\sum_{i_{q_{z}+1}\in N_{q_{z}+1}}\cdots\sum_{i_{k}\in N_{k}}\\*
 &\quad \mathbb{E} \biggl[X_{i_{q_{z}}}\cdots X_{i_{k-1}}\Bigl(\sum_{i_{k}\in N_{k}}\lvert X_{i_{k}}\rvert\Bigr)^{\omega }\biggr].
  \end{align*}

We call $\omega$ the exponent of the $\mathcal{R}$-sum. If $\omega=1$, the only difference between an $\mathcal{R}$-sum and an $\mathcal{S}$-sum is that the $X_{i_{j}}$'s in \eqref{eq:defcomp1} are replaced by $\lvert X_{i_{j}} \rvert$'s in \eqref{eq:defcomp3}. Thus, an $\mathcal{S}$-sum is always upper-bounded by the corresponding compositional $1$-sum, i.e.,

  \begin{equation}\label{eq:type13}
    \bigl\lvert \mathcal{S} [t_{1},\cdots,t_{k}]\bigr\rvert\leq \mathcal{R}_{1} [t_{1},\cdots,t_{k}].
  \end{equation} 

Another important observation is that we can compare the values of $\mathcal{R}$-sums with respect to two different sequences $t_{1},\cdots,t_{k}$ and $t_{1}',\cdots,t_{k}'$ in certain situations. In specific, if for any $j\in [k]$ we have that if $t_j$ and $t'_j$ are of the same sign and $ |t_{j}|\leq  |t_{j}'| $, then
      \begin{equation}\label{eq:comparecomp1}
        \mathcal{R}_{\omega}[t_{1},\cdots,t_{k}]\leq \mathcal{R}_{\omega}[t_{1}',\cdots,t_{k}'].
      \end{equation}
  In fact, the sequences $(t_{j})$ and $(t_{j}')$ having the same sign indicates that $\{ j:t_{j}>0 \}=\{ j:t_{j}'>0 \}$.
  Thus, we can write

  \begin{equation}\label{eq:sthneedstocompare}
  \begin{aligned}
     \mathcal{R}_{\omega} [t_{1}',\cdots, t_{k}']
  = &\sum_{i_{1}\in N_{1}'}\sum_{i_{2}\in N_{2}'}\cdots\sum_{i_{k-1}\in N_{k-1}'}[q_{1}-q_{0},\cdots,q_{z+1}-q_{z}]\triangleright\\ 
  &\quad\biggl(\lvert X_{i_{1}}\rvert,\cdots,\lvert X_{i_{k-1}}\rvert,\Bigl(\sum_{i_{k}\in N_{k}'}\lvert X_{i_{k}}\rvert\Bigr)^{\omega }\biggr),
  \end{aligned}
\end{equation}
  where we note that $N_{1}'=I=N_{1}$ and for $j=2,\cdots,k$ we have
  \begin{equation*}
    N_{j}'=N(i_{1},\cdots,i_{\lvert t_{j}' \rvert})\supseteq N(i_{1},\cdots,i_{\lvert t_{j} \rvert})=N_{j}.
  \end{equation*}
  By comparing \eqref{eq:defcomp3} with \eqref{eq:sthneedstocompare}, we obtain \eqref{eq:comparecomp1}.

\subsubsection*{Re-expression of the remainder terms $R_{k,\omega}$}

Using the notion of $\mathcal{R}$-sums, we rewrite the $R_{k,\omega}$ in \cref{SEC:LOCAL} as
  \begin{equation}\label{eq:rkalpha}
    \begin{aligned}
    R_{k,\omega}:= & \sum_{(\ell,\eta_{1:\ell})\in C^{*}(k+2)}\sum_{i_{1}\in N_{1}'}\sum_{i_{2}\in N_{2}'}\cdots\sum_{i_{k+1}\in N_{k+1}'}\\
    &\quad [\eta_{1},\cdots,\eta_{\ell}]\triangleright \biggl(\lvert X_{i_{1}}\rvert,\cdots,\lvert X_{i_{k+1}}\rvert,\Bigl(\sum_{i_{k+2}\in N_{k+2}'}\lvert X_{i_{k+2}}\rvert\Bigr)^{\omega }\biggr) \\
     = &\sum_{t_{1:(k+2)}\in \mathcal{M}_{1,k+2}}\ \mathcal{R}_{\omega}[t_{1},t_{2},\cdots,t_{k+2}].
  \end{aligned}
\end{equation}
where $N_{1}':=I$ and $N_{j}':=N(i_{1:(j-1)})$ for $j\geq 2$. $C^{*}(k+2)$ and $\mathcal{M}_{1,k+2}$ are given by
$$
C^{*}(k+2)=\bigl\{\ell,\eta_{1:\ell}\in\mathbb{N}_{+}: \eta_{j}\geq 2\ \forall  j\in [\ell-1], \ \sum_{j=1}^{\ell}\eta_{j}=k+2\bigr\},
$$ and 
\begin{equation*}
  \mathcal{M}_{1,k+2}:=\Bigl\{t_{1:(k+2)}:~ t_{j+1}=\pm j\ \ \&\ \  t_{j}\wedge t_{j+1}< 0\ \ \forall 1\leq j\leq k+1\Bigr\}.
\end{equation*} 
Note that $t_{j}\wedge t_{j+1}<0$ for any $j\in [k+1]$ means that there is at least one $-1$ in any two consecutive signs, which corresponds to the requirement that $\eta_{j}\geq 2$ for $j\in [\ell-1]$ in \eqref{eq:rkalpha}.

\addtocontents{toc}{\SkipTocEntry}
\subsection{Proofs of Proposition \ref{THM:WFWEXPANSION} and Lemma \ref{THM:BARBOURLIKE}}\label{sec:pflocalkeylemma}

In this section, we carry out the local expansion technique and prove \cref{THM:WFWEXPANSION,THM:BARBOURLIKE}.

  Firstly, we establish the following lemma, which will be crucial in the inductive step of proving the main theorem.

\begin{lemma}\label{thm:proptaylor}
  Fix $k\in \mathbb{N}_{+}$. For any $s\in [k]\cup \{ 0 \}$ and $f\in \mathcal{C}^{k,\omega}(\mathbb{R})$, we have
  \begin{equation}
  \begin{aligned}\label{eq:taylorcomp4}
    &\Bigl\lvert \mathcal{T}_{f,s}[t_{1},\cdots,t_{k+1}]-\mathcal{S}[t_{1},\cdots,t_{k+1}]\ \mathbb{E} [\partial^{k}f(W)]\Bigr\rvert\\
    \leq & \lvert f \rvert_{k,\omega}\bigl(-(\mathbb{I}(t_{k+1}<0)\cdot \mathcal{R}_{\omega}[t_{1},\cdots,t_{k+1},k+1]
    + \mathbb{I}(s\ge 1)\mathcal{R}_{\omega}[t_{1},\cdots,t_{k+1},-(k+1)]\bigr).
  \end{aligned}
\end{equation}
  Given any $\ell\in [k]$ and $s\in [\ell]\cup \{ 0 \}$, we further have
    \begin{align}\label{eq:taylorcomp3}
    &\biggl\lvert \mathcal{T}_{f,s}[t_{1},\cdots,t_{\ell}]-\mathcal{S}[t_{1},\cdots,t_{\ell}]\ \mathbb{E} [\partial^{\ell-1}f(W)]\\*
    &\ -\mathbb{I}(s\ge 1)\cdot\sum_{j=1}^{k-\ell+1}\sum_{h=0}^{j}(-1)^{h}\frac{1}{h !(j-h)!}
    \mathcal{T}_{f,j}[t_{1},\cdots,t_{\ell},\underbrace{s - \ell,\cdots,s  - \ell}_{h\text{ times}},\underbrace{-\ell,\cdots,-\ell}_{(j-h)\text{ times}}] \nonumber\\*
    &\ -(\mathbb{I}(t_{\ell}<0)\sum_{j=1}^{k-\ell+1}\frac{1}{j!}
    \mathcal{T}_{f,j}[t_{1},\cdots,t_{\ell},\ell,\underbrace{-\ell,\cdots,-\ell}_{(j-1)\text{ times}}] \biggr\rvert\nonumber\\
    \leq & \frac{\lvert f \rvert_{k,\omega}}{(k-\ell+1)!}\bigl(-\mathbb{I}(t_{\ell}<0)\cdot\mathcal{R}_{\omega}[t_{1},\cdots,t_{\ell},\ell,\underbrace{-\ell,\cdots,-\ell}_{(k-\ell+1)\text{ times}}]\nonumber\\*
    &\quad+ \mathbb{I}(s\ge 1)\cdot\mathcal{R}_{\omega}[t_{1},\cdots,t_{\ell},\underbrace{-\ell,\cdots,-\ell}_{(k-\ell+2)\text{ times}}]\bigr).\nonumber
    \end{align}
\end{lemma}
%

%\textcolor{cyan}{[Note: the proof is rewritten]}

\begin{proof}
  Firstly, we remark that the definition of Hölder continuity implies that
  \begin{equation}\label{eq:somesomeholder}
    \bigl\lvert\partial^{k}f(y)-\partial^{k}f (x)\bigr\rvert\leq \lvert f \rvert_{k,\omega }\lvert y-x \rvert^{\omega },
  \end{equation}
  where $\omega$ is the Hölder exponent of $f$ and $\lvert f \rvert_{k,\omega}$ is the Hölder constant (see \cref{thm:defholder}).
  Let $z=\bigl\lvert\{ j\in [k+1]: t_{j}>0 \}\bigr\rvert$ be the number of positive indexes $(t_j)$. If $z\geq 1$, we write $\{ j\in [k+1]:t_{j}>0 \}=\{ q_{1},\cdots,q_{z} \}$. Without loss of generality, we suppose that the sequence $2\leq q_{1}<\cdots<q_{z}\leq k+1$ is increasing. We further let $q_{0}:=1$ and $q_{z+1}:=k+2$. 
  Applying \eqref{eq:somesomeholder} we have
  \begin{align}\label{eq:someannoyingthing}
    &\Bigl\lvert\mathbb{E} \bigl[X_{i_{q_{z}}}\cdots X_{i_{k+1}}\partial^{k}f\bigl(W_{i.}[k+1-s]\bigr)\bigr]-\mathbb{E} \bigl[X_{i_{q_{z}}}\cdots X_{i_{k+1}}\partial^{k}f\bigl(W_{i.}[k+1]\bigr)\bigr]\Bigr\rvert\\
    \leq & \lvert f \rvert_{k,\omega}\mathbb{E} \bigl[\bigl\lvert X_{i_{q_{z}}}\cdots X_{i_{k+1}} \bigr\rvert\cdot  \bigl\lvert W_{i.}[k+1-s]-W_{i.}[k+1] \bigr\rvert^{\omega}\bigr]\nonumber\\
    \leq &\lvert f \rvert_{k,\omega}\mathbb{E} \biggl[\bigl\lvert X_{i_{q_{z}}}\cdots X_{i_{k+1}} \bigr\rvert\cdot  \Bigl\lvert \sum_{i\in N(i_{1:(k+1)})\backslash N (i_{1:(k+1-s)})} X_{i}\Bigr\rvert^{\omega}\biggr]\nonumber\\
    \leq &\lvert f \rvert_{k,\omega}\mathbb{E} \biggl[\bigl\lvert X_{i_{q_{z}}}\cdots X_{i_{k+1}} \bigr\rvert\cdot  \Bigl\lvert \sum_{i\in N(i_{1:(k+1)})} X_{i}\Bigr\rvert^{\omega}\biggr],\nonumber
  \end{align}
where in the last inequality we have used the fact that $N(i_{1:(k+1)})\backslash N (i_{1:(k+1-s)})\subseteq N(i_{1:(k+1)})$.
If $z=0$, this directly implies that
\begin{align}\label{eq:firstwannaprove}
  & \Bigl\lvert \mathcal{T}_{f,s}[t_{1},\cdots,t_{k+1}]-\mathcal{T}_{f,0}[t_{1},\cdots,t_{k+1}] \Bigr\rvert
  \leq \mathbb{I}(s\ge 1)\cdot \lvert f \rvert_{k,\omega}\mathcal{R}_{\omega}[t_{1},\cdots,t_{k+1},-(k+1)].
\end{align}
If $z\geq 1$, by definition \eqref{eq:defcomp2} we have for $s\geq 1$
  \begin{align*}
    & \Bigl\lvert \mathcal{T}_{f,s}[t_{1},\cdots,t_{k+1}]-\mathcal{T}_{f,0}[t_{1},\cdots,t_{k+1}] \Bigr\rvert\\
    =& \biggl\lvert \sum_{i_{1}\in N_{1}}\sum_{i_{2}\in N_{2}}\cdots\sum_{i_{k+1}\in N_{k+1}}\mathbb{E} \bigl[X_{i_{q_{0}}}\cdots X_{i_{q_{1}-1}}\bigr]\ \cdots\ \mathbb{E} \bigl[X_{i_{q_{z-1}}}\cdots X_{i_{q_{z}-1}}\bigr]\cdot\\*
    &\  \mathbb{E} \bigl[X_{i_{q_{z}}}\cdots X_{i_{k+1}}\bigl(\partial^{k}f(W_{i.}[k+1-s])-\partial^{k}f(W_{i.}[k+1])\bigr)\bigr] \biggr\rvert\\
    \leq & \sum_{i_{1}\in N_{1}}\sum_{i_{2}\in N_{2}}\cdots\sum_{i_{k+1}\in N_{k+1}}\mathbb{E} \bigl[\bigl\lvert X_{i_{q_{0}}}\cdots X_{i_{q_{1}-1}}\bigr\rvert \bigr]\ \cdots\ \mathbb{E} \bigl[\bigl\lvert X_{i_{q_{z-1}}}\cdots X_{i_{q_{z}-1}}\bigr\rvert \bigr]\cdot\\*
    &\ \Bigl\lvert \mathbb{E} \bigl[ X_{i_{q_{z}}}\cdots X_{i_{k+1}}\textbf{}\partial^{k}f(W_{i.}[k+1-s])-\partial^{k}f(W_{i.}[k+1])\bigr]\Bigr\rvert\\
    \overset{\eqref{eq:someannoyingthing}}{\leq}& \lvert f \rvert_{k,\omega}\sum_{i_{1}\in N_{1}}\sum_{i_{2}\in N_{2}}\cdots\sum_{i_{k+1}\in N_{k+1}}\mathbb{E} \bigl[\bigl\lvert X_{i_{q_{0}}}\cdots X_{i_{q_{1}-1}}\bigr\rvert \bigr]\ \cdots\ \mathbb{E} \bigl[\bigl\lvert X_{i_{q_{z-1}}}\cdots X_{i_{q_{z}-1}}\bigr\rvert \bigr]\ \cdot\\*
    &\qquad \mathbb{E} \biggl[\bigl\lvert X_{i_{q_{z}}}\cdots X_{i_{k+1}} \bigr\rvert\cdot  \Bigl\lvert \sum_{i\in N(i_{1:(k+1)})} X_{i}\Bigr\rvert^{\omega}\biggr]\\
    =&\lvert f \rvert_{k,\omega}\mathcal{R}_{\omega}[t_{1},\cdots,t_{k+1},-(k+1)].
  \end{align*}
  Here the last equality is due to the definition \eqref{eq:defcomp3}. 
  Thus, \eqref{eq:firstwannaprove} is proven for both $z=0$ and $z\geq 1$. Next we show that
  \begin{equation}\label{eq:secondwannaprove}
    \begin{aligned}
    &\Bigl\lvert \mathcal{S}[t_{1},\cdots,t_{k+1}] \bigl(\mathbb{E} [\partial^{k}f(W)]-\mathbb{E} [\partial^{k}f(W_{i.}[k+1])]\bigr)\Bigr\rvert\\
    \leq &-\mathbb{I}(t_{k+1}<0)\lvert f \rvert_{k,\omega}\mathcal{R}_{\omega}[t_{1},\cdots,t_{k+1},k+1].
    \end{aligned}
  \end{equation}
In this goal, we first note that if $t_{k+1}\geq 0$, by definition \eqref{eq:defcomp1} we know that $q_z=k+1$ and therefore, according to \eqref{j} we know that $$\mathcal{S}[t_{1},\cdots,t_{k+1}]=0,$$and so \eqref{eq:secondwannaprove} holds. Otherwise, we note that we have
\begin{equation}
\begin{aligned}\label{eq:someannoyingthing2}
  &\Bigl\lvert\mathbb{E} \bigl[\partial^{k}f\bigl(W\bigr)\bigr]-\mathbb{E} \bigl[\partial^{k}f\bigl(W_{i.}[k+1]\bigr)\bigr]\Bigr\rvert
  \leq  \lvert f \rvert_{k,\omega}\mathbb{E} \bigl[ \bigl\lvert W_{i.}[k+1-s]-W_{i.}[k+1] \bigr\rvert^{\omega}\bigr]\\
  \leq &\lvert f \rvert_{k,\omega}\mathbb{E} \biggl[ \Bigl\lvert \sum_{i\in N(i_{1:(k+1)})\backslash N (i_{1:(k+1-s)})} X_{i}\Bigr\rvert^{\omega}\biggr]
  \leq \lvert f \rvert_{k,\omega}\mathbb{E} \biggl[ \Bigl\lvert \sum_{i\in N(i_{1:(k+1)})} X_{i}\Bigr\rvert^{\omega}\biggr].
\end{aligned}
\end{equation}
This implies that
\begin{align*}
  &\Bigl\lvert \mathcal{S}[t_{1},\cdots,t_{k+1}] \bigl(\mathbb{E} [\partial^{k}f(W)]-\mathbb{E} [\partial^{k}f(W_{i.}[k+1])]\bigr)\Bigr\rvert\\
  \leq &\Bigl\lvert \mathcal{S}[t_{1},\cdots,t_{k+1}] \Bigr\rvert\cdot \Bigl\lvert\mathbb{E} \bigl[\partial^{k}f\bigl(W\bigr)\bigr]-\mathbb{E} \bigl[\partial^{k}f\bigl(W_{i.}[k+1]\bigr)\bigr]\Bigr\rvert\\
  \overset{(*)}{\leq}&\lvert f \rvert_{k,\omega}\mathcal{R}_{1}[t_{1},\cdots,t_{k+1}]\ \mathbb{E} \biggl[ \Bigl\lvert \sum_{i\in N(i_{1:(k+1)})} X_{i}\Bigr\rvert^{\omega}\biggr]\\
  =& \lvert f \rvert_{k,\omega}\sum_{i_{1}\in N_{1}}\sum_{i_{2}\in N_{2}}\cdots\sum_{i_{k+1}\in N_{k+1}}[q_{1}-q_{0},\cdots,q_{z+1}-q_{z}]\triangleright \\*
  &\quad\bigl(\lvert X_{i_{1}}\rvert,\cdots,\lvert X_{i_{k+1}}\rvert\bigr)\cdot \mathbb{E} \biggl[ \Bigl\lvert \sum_{i\in N(i_{1:(k+1)})} X_{i}\Bigr\rvert^{\omega}\biggr]\\
  =& \lvert f \rvert_{k,\omega}\sum_{i_{1}\in N_{1}}\sum_{i_{2}\in N_{2}}\cdots\sum_{i_{k+1}\in N_{k+1}}[q_{1}-q_{0},\cdots,q_{z+1}-q_{z},1]\triangleright\\* 
  &\quad\Bigl(\lvert X_{i_{1}}\rvert,\cdots,\lvert X_{i_{k+1}}\rvert,\bigl(\sum_{i_{k+2}\in N(i_{1:(k+1)})}\lvert X_{i_{k+2}}\rvert\bigr)^{\omega }\Bigr)\\
  =& \lvert f \rvert_{k,\omega}\mathcal{R}_{\omega}[t_{1},\cdots,t_{k+1},k+1],
\end{align*}
where $(*)$ is due to \eqref{eq:type13} and \eqref{eq:someannoyingthing2}.
Taking the difference of \eqref{eq:firstwannaprove} and \eqref{eq:secondwannaprove}, we obtain \eqref{eq:taylorcomp4} by applying the equation \eqref{eq:indepcomp}.

  For $\ell\leq k$, we apply the Taylor expansion with remainders taking the integral form and obtain that
  \begin{align}\label{eq:someholder}
                       & \partial^{\ell-1} f(y)-\partial^{\ell-1} f(x)
    =  \sum_{j=1}^{m-\ell}\frac{1}{j!}(y-x)^{j}\partial^{\ell-1+j}f(x)                                                                                                          \\
                       & +\frac{1}{(k-\ell+1)!}(y-x)^{k-\ell+1}\int_{0}^{1}(k-\ell+1)v^{k-\ell}\partial^{k} f(v x+(1-v)y)\dif v  \nonumber                                                           \\
    \overset{(*)}{  =} & \sum_{j=1}^{k-\ell+1}\frac{1}{j!}(y-x)^{j}\partial^{\ell-1+j}f(x)                               \nonumber                                                          \\
                       & +\frac{1}{(k-\ell+1)!}(y-x)^{k-\ell+1}\int_{0}^{1}(k-\ell+1) v ^{k-\ell}\bigl(\partial^{k} f( v  x+(1- v )y)-\partial^{k} f(x)\bigr)\dif  v,\nonumber
  \end{align}
  where to obtain $(*)$ we added and subtracted $\frac{(y-x)^{k-\ell+1}}{(k-\ell+1)!}\partial^k f(x).$  Moreover, using the fact that $\partial^k f(\cdot)$ is assumed to be Hölder  continuous we obtain that
  \begin{equation}\label{cabris}
    \bigl\lvert \partial^{k} f( v  x+(1- v )y)-\partial^{k} f(x)\bigr\rvert\leq \lvert f \rvert_{k,\omega}(1-v)^{\omega}\ \lvert y-x \rvert^{\omega}\leq \lvert f \rvert_{k,\omega}\lvert y-x \rvert^{\omega}.
  \end{equation}
  Therefore, as $\int_{0}^{1}(k-\ell+1)v^{k-\ell}\dif v =1$, by combining \eqref{cabris} with \eqref{eq:someholder} we get that
  \begin{equation}\label{eq:whattaylor}
    \biggl\lvert \partial^{\ell-1} f(y)-\partial^{\ell-1} f(x)- \sum_{j=1}^{k-\ell+1}\frac{1}{j!}(y-x)^{j}\partial^{\ell-1+j}f(x)\biggr\rvert     \leq  \frac{\lvert f \rvert_{k,\omega }}{(k-\ell+1)!}\lvert y-x \rvert^{k-\ell+1+\omega }.
  \end{equation}

  We prove that the following inequality holds:
    \begin{align}
    &\begin{aligned}\label{eq:taylorcomp1}
      &\biggl\lvert \mathcal{T}_{f,s}[t_{1},\cdots,t_{\ell}]-\mathcal{T}_{f,0}[t_{1},\cdots,t_{\ell}]\\*
      &\ -\mathbb{I}(s\geq 1)\cdot\sum_{j=1}^{k-\ell+1}\sum_{h=0}^{j}(-1)^{h}\frac{1}{h !(j-h)!}
      \mathcal{T}_{f,j}[t_{1},\cdots,t_{\ell},\underbrace{s - \ell,\cdots,s  - \ell}_{h\text{ times}},\underbrace{-\ell,\cdots,-\ell}_{(j-h)\text{ times}}] \biggr\rvert\\
      \leq & \frac{\mathbb{I}(s\ge 1)\cdot\lvert f \rvert_{k,\omega}}{(k-\ell+1)!}\mathcal{R}_{\omega}[t_{1},\cdots,t_{\ell},\underbrace{-\ell,\cdots,-\ell}_{(k-\ell+2)\text{ times}}],
    \end{aligned}
    \end{align}

  First, let's establish \eqref{eq:taylorcomp1}. Let $z=\bigl\lvert\{ j\in [\ell]: t_{j}>0 \}\bigr\rvert$. If $z\geq 1$, we write $\{ j\in [\ell]:t_{j}>0 \}=\{ q_{1},\cdots,q_{z} \}$. Without loss of generality, we suppose that the sequence $2\leq q_{1}<\cdots<q_{z}\leq \ell$ is increasing. We further let $q_{0}:=1$ and $q_{z+1}:=\ell+1$. 
  Applying \eqref{eq:whattaylor} we have
  \begin{equation}
  \begin{aligned}\label{eq:someannoyingthing3}
    &\biggl\lvert\mathbb{E} \bigl[X_{i_{q_{z}}}\cdots X_{i_{\ell}}\partial^{\ell-1}f\bigl(W_{i.}[\ell-s]\bigr)\bigr]-\mathbb{E} \bigl[X_{i_{q_{z}}}\cdots X_{i_{\ell}}\partial^{\ell-1}f\bigl(W_{i.}[\ell]\bigr)\bigr]\\
    &\ -\sum_{j=1}^{k-\ell+1}\frac{1}{j!}\mathbb{E} \bigl[ X_{i_{q_{z}}}\cdots X_{i_{\ell}}    (W_{i.}[\ell-s]-W_{i.}[\ell] )^{j}\partial^{\ell-1+j}f(W_{i.}[\ell])\bigr]
    \biggr\rvert\\
    \leq & \frac{\lvert f \rvert_{k,\omega}}{(k-\ell+1)!}\mathbb{E} \bigl[\bigl\lvert X_{i_{q_{z}}}\cdots X_{i_{\ell}} \bigr\rvert\cdot  \bigl\lvert W_{i.}[\ell-s]-W_{i.}[\ell] \bigr\rvert^{k-\ell+1+\omega}\bigr].
  \end{aligned}
\end{equation}

  For convenience let
  \begin{align*}
    &E_{1}:=\sum_{i_{q_{z}}\in N_{q_{z}}}\cdots\sum_{i_{\ell}\in N_{\ell}}\mathbb{E} \bigl[X_{i_{q_{z}}}\cdots X_{i_{\ell}}\partial^{\ell-1}f\bigl(W_{i.}[\ell-s]\bigr)\bigr]-\mathbb{E} \bigl[X_{i_{q_{z}}}\cdots X_{i_{\ell}}\partial^{\ell-1}f\bigl(W_{i.}[\ell]\bigr)\bigr],\\
    &E_{2,j}:=\sum_{i_{q_{z}}\in N_{q_{z}}}\cdots\sum_{i_{\ell}\in N_{\ell}}\mathbb{E} \bigl[ X_{i_{q_{z}}}\cdots X_{i_{\ell}}    (W_{i.}[\ell-s]-W_{i.}[\ell] )^{j}\partial^{\ell-1+j}f(W_{i.}[\ell])\bigr],\\
    &E_{3}:=\sum_{i_{q_{z}}\in N_{q_{z}}}\cdots\sum_{i_{\ell}\in N_{\ell}}\mathbb{E} \bigl[\bigl\lvert X_{i_{q_{z}}}\cdots X_{i_{\ell}} \bigr\rvert\cdot  \bigl\lvert W_{i.}[\ell-s]-W_{i.}[\ell] \bigr\rvert^{k-\ell+1+\omega}\bigr].
  \end{align*}
  Then \eqref{eq:someannoyingthing3} reduces to
  \begin{equation}\label{eq:derivedfromtaylor}
  \textstyle\bigl\lvert E_{1}-\sum_{j=1}^{k-\ell+1}E_{2,j}/j! \bigr\rvert\leq \lvert f \rvert_{k,\omega}E_{3}/(k-\ell+1)!.
  \end{equation}
  Then we observe that by definition of $W_{i.}[\cdot]$ we have
  \begin{equation}
  \begin{aligned}\label{eq:annoyance1}
    &\mathbb{E} \bigl[ X_{i_{q_{z}}}\cdots X_{i_{\ell}}    (W_{i.}[\ell-s]-W_{i.}[\ell] )^{j}\partial^{\ell-1+j}f(W_{i.}[\ell])\bigr]\\
    =& \mathbb{E} \biggl[ X_{i_{q_{z}}}\cdots X_{i_{\ell}}    \Bigl(\sum_{i\in N(i_{1:\ell})}X_{i}-\sum_{i\in N(i_{1:\ell-s})}X_{i} \Bigr)^{j}\partial^{\ell-1+j}f(W_{i.}[\ell])\biggr]\\
    =& \sum_{h=0}^{j}(-1)^{h}\binom{j}{h}\ \mathbb{E} \biggl[ X_{i_{q_{z}}}\cdots X_{i_{\ell}}  \Bigl(\sum_{i\in N(i_{1:\ell-s})}X_{i} \Bigr)^{h}  \Bigl(\sum_{i\in N(i_{1:\ell})}X_{i}\Bigr)^{j-h}\partial^{\ell-1+j}f(W_{i.}[\ell])\biggr],
  \end{aligned}
\end{equation}
and that
\begin{align}\label{eq:annoyance2}
  &\mathbb{E} \bigl[\bigl\lvert X_{i_{q_{z}}}\cdots X_{i_{\ell}} \bigr\rvert\cdot  \bigl\lvert W_{i.}[\ell-s]-W_{i.}[\ell] \bigr\rvert^{k-\ell+1+\omega}\bigr]\\\nonumber
  \leq &\mathbb{E} \biggl[\bigl\lvert X_{i_{q_{z}}}\cdots X_{i_{k+1}} \bigr\rvert\cdot  \Bigl\lvert \sum_{i\in N(i_{1:(k+1)})\backslash N (i_{1:(k+1-s)})} X_{i}\Bigr\rvert^{k-\ell+1+\omega}\biggr]\\\nonumber
  \leq &\mathbb{E} \biggl[\bigl\lvert X_{i_{q_{z}}}\cdots X_{i_{k+1}} \bigr\rvert\cdot  \Bigl\lvert \sum_{i\in N(i_{1:(k+1)})} X_{i}\Bigr\rvert^{k-\ell+1+\omega}\biggr]\\\nonumber
  \leq &\mathbb{E} \biggl[\bigl\lvert X_{i_{q_{z}}}\cdots X_{i_{k+1}} \bigr\rvert\cdot  \Bigl( \sum_{i\in N(i_{1:(k+1)})} \lvert X_{i}\rvert\Bigr)^{k-\ell+1}\cdot \Bigl\lvert \sum_{i\in N(i_{1:(k+1)})} X_{i}\Bigr\rvert^{\omega}\biggr].\nonumber
\end{align}
  If $z=0$, we take the sum of \eqref{eq:annoyance1} or \eqref{eq:annoyance2} over $i_{q_{z}}\in N_{q_{z}},\cdots,i_{\ell}\in N_{\ell}$. By definition \eqref{eq:defcomp2} and \eqref{eq:defcomp3} we have
  \begin{equation}\label{eq:annoyance4}
    \begin{aligned}
      & E_{1}=\mathcal{T}_{f,s}[t_{1},\cdots,t_{\ell}]-\mathcal{T}_{f,0}[t_{1},\cdots,t_{\ell}],\\
      &E_{2,j}=\sum_{h=0}^{j}(-1)^{h}\binom{j}{h}\ \mathcal{T}_{f,j}[t_{1},\cdots,t_{\ell},\underbrace{s - \ell,\cdots,s  - \ell}_{h\text{ times}},\underbrace{-\ell,\cdots,-\ell}_{(j-h)\text{ times}}],\\
      &E_{3}\leq \mathcal{R}_{\omega}[t_{1},\cdots,t_{\ell},\ell,\underbrace{-\ell,\cdots,-\ell}_{(k-\ell+1)\text{ times}}].
      \end{aligned}
  \end{equation}
  Combining \eqref{eq:annoyance4} and \eqref{eq:derivedfromtaylor}, we have for $s\geq 1$
  \begin{align*}
    &\biggl\lvert \mathcal{T}_{f,s}[t_{1},\cdots,t_{\ell}]-\mathcal{T}_{f,0}[t_{1},\cdots,t_{\ell}]\\*
    &\ -\sum_{j=1}^{k-\ell+1}\sum_{h=0}^{j}(-1)^{h}\frac{1}{h !(j-h)!}
    \mathcal{T}_{f,j}[t_{1},\cdots,t_{\ell},\underbrace{s - \ell,\cdots,s  - \ell}_{h\text{ times}},\underbrace{-\ell,\cdots,-\ell}_{(j-h)\text{ times}}] \biggr\rvert\\
    \overset{\eqref{eq:annoyance4}}{=}&\bigl\lvert E_{1}-\sum_{j=1}^{k-\ell+1}E_{2,j}/j! \bigr\rvert\overset{\eqref{eq:derivedfromtaylor}}{\leq} \lvert f \rvert_{k,\omega}E_{3}/(k-\ell+1)!\\
    \overset{\eqref{eq:annoyance4}}{\leq} & \frac{\lvert f \rvert_{k,\omega}}{(k-\ell+1)!}\mathcal{R}_{\omega}[t_{1},\cdots,t_{\ell},\underbrace{-\ell,\cdots,-\ell}_{(k-\ell+2)\text{ times}}].
  \end{align*}
  Thus, \eqref{eq:taylorcomp1} holds for $z=0$.
  
  If $z\geq 1$, similar to \eqref{eq:annoyance2} we have
  \begin{equation}\label{eq:annoyance3}
    \begin{aligned}
    & \mathcal{S}[t_{1},\cdots,t_{q_{z}-1}]
    \cdot E_{1}=\mathcal{T}_{f,s}[t_{1},\cdots,t_{\ell}]-\mathcal{T}_{f,0}[t_{1},\cdots,t_{\ell}],\\
    &\mathcal{S}[t_{1},\cdots,t_{q_{z}-1}]
    \cdot E_{2,j}
  =\sum_{h=0}^{j}(-1)^{h}\binom{j}{h}\ \mathcal{T}_{f,j}[t_{1},\cdots,t_{\ell},\underbrace{s - \ell,\cdots,s  - \ell}_{h\text{ times}},\underbrace{-\ell,\cdots,-\ell}_{(j-h)\text{ times}}],\\
    &\mathcal{R}_{1}[t_{1},\cdots,t_{q_{z}-1}]
    \cdot E_{3}\leq \mathcal{R}_{\omega}[t_{1},\cdots,t_{\ell},\underbrace{-\ell,\cdots,-\ell}_{(k-\ell+2)\text{ times}}].
    \end{aligned}
  \end{equation}
  Combining \eqref{eq:annoyance3} and \eqref{eq:derivedfromtaylor} we get for $s\geq 1$
  \begin{align*}
    &\biggl\lvert \mathcal{T}_{f,s}[t_{1},\cdots,t_{\ell}]-\mathcal{T}_{f,0}[t_{1},\cdots,t_{\ell}]\\*
      &\ -\sum_{j=1}^{k-\ell+1}\sum_{h=0}^{j}(-1)^{h}\frac{1}{h !(j-h)!}
      \mathcal{T}_{f,j}[t_{1},\cdots,t_{\ell},\underbrace{s - \ell,\cdots,s  - \ell}_{h\text{ times}},\underbrace{-\ell,\cdots,-\ell}_{(j-h)\text{ times}}] \biggr\rvert\\
    \overset{\eqref{eq:annoyance3}}{=}& \bigl\lvert \mathcal{S}[t_{1},\cdots,t_{q_{z}-1}] \bigr\rvert\cdot \bigl\lvert E_{1}-\sum_{j=1}^{k-\ell+1}E_{2,j}/j! \bigr\rvert
    \overset{\eqref{eq:type13}}{\leq}  \mathcal{R}_{1}[t_{1},\cdots,t_{q_{z}-1}]\cdot \bigl\lvert E_{1}-\sum_{j=1}^{k-\ell+1}E_{2,j}/j! \bigr\rvert\\
    \overset{\eqref{eq:derivedfromtaylor}}{\leq} & \mathcal{R}_{1}[t_{1},\cdots,t_{q_{z}-1}]\cdot \lvert f \rvert_{k,\omega}E_{3}/(k-\ell+1)!
    \overset{\eqref{eq:annoyance3}}{\leq}  \frac{\lvert f \rvert_{k,\omega}}{(k-\ell+1)!}\mathcal{R}_{\omega}[t_{1},\cdots,t_{\ell},\underbrace{-\ell,\cdots,-\ell}_{(k-\ell+2)\text{ times}}].
  \end{align*}
Thus, we have shown \eqref{eq:taylorcomp1} for both $z=0$ and $z\geq 1$.

Next we prove that the following inequality holds:
    \begin{align}
    \label{eq:taylorcomp2}
    &\begin{aligned}
      &\biggl\lvert \mathcal{S}[t_{1},\cdots,t_{\ell}]\bigl(\mathbb{E} [\partial^{\ell-1}f(W)]-\mathbb{E} [\partial^{\ell-1}f(W_{i.}[\ell])]\bigr)\\*
      &\ -\mathbb{I}(t_{\ell}<0)\cdot\sum_{j=1}^{k-\ell+1}\frac{1}{j!}
      \mathcal{T}_{f,j}[t_{1},\cdots,t_{\ell},\ell,\underbrace{-\ell,\cdots,-\ell}_{(j-1)\text{ times}}] \biggr\rvert\\
      \leq & \frac{\mathbb{I}(t_{\ell}<0)\cdot\lvert f \rvert_{k,\omega}}{(k-\ell+1)!}\mathcal{R}_{\omega}[t_{1},\cdots,t_{\ell},\ell,\underbrace{-\ell,\cdots,-\ell}_{(k-\ell+1)\text{ times}}].
    \end{aligned}
    \end{align}
For \eqref{eq:taylorcomp2}, we apply \eqref{eq:whattaylor} again and get that
\begin{equation}\label{eq:someannoyingthing4}
  \begin{aligned}
  &\biggl\lvert\mathbb{E} \bigl[\partial^{k}f(W)\bigr]-\mathbb{E} \bigl[\partial^{k}f(W_{i.}[\ell])\bigr]\\
  &\ -\sum_{j=1}^{k-\ell+1}\frac{1}{j!}\mathbb{E} \bigl[     (W_{i.}-W_{i.}[\ell] )^{j}\partial^{\ell-1+j}f(W_{i.}[\ell])\bigr]
  \biggr\rvert
  \leq \frac{\lvert f \rvert_{k,\omega}}{(k-\ell+1)!}\mathbb{E} \bigl[  \bigl\lvert W-W_{i.}[\ell] \bigr\rvert^{k-\ell+1+\omega}\bigr].
  \end{aligned}
\end{equation}
For convenience let
\begin{align*}
  &E_{4}:=\mathbb{E} [\partial^{\ell-1}f(W)]-\mathbb{E} [\partial^{\ell-1}f(W_{i.}[\ell])],\\
  &E_{5,j}:=\mathbb{E} \bigl[     (W_{i.}-W_{i.}[\ell] )^{j}\partial^{\ell-1+j}f(W_{i.}[\ell])\bigr],\\
  &E_{6}:=\mathbb{E} \bigl[  \bigl\lvert W-W_{i.}[\ell] \bigr\rvert^{k-\ell+1+\omega}\bigr].
\end{align*}
Then \eqref{eq:someannoyingthing4} reduces to
\begin{equation}\label{eq:someannoyingthing5}
  \textstyle\bigl\lvert E_{4}-\sum_{j=1}^{k-\ell+1}E_{5,j}/j! \bigr\rvert\leq \lvert f \rvert_{k,\omega}E_{6}/(k-\ell+1)!.
\end{equation} We first note that if $t_{\ell}\ge 0$ then $\mathcal{S}[t_1,\cdots,t_{\ell}]=0$ therefore, \eqref{eq:taylorcomp2} holds. 
Moreover, similar to \eqref{eq:annoyance3}, we have for $t_{\ell}<0$
\begin{equation}\label{eq:annoyance5}
  \begin{aligned}
    & \mathcal{S}[t_{1},\cdots,t_{\ell}]\cdot E_{4}=\mathcal{S}[t_{1},\cdots,t_{\ell}]\bigl(\mathbb{E} [\partial^{\ell-1}f(W)]-\mathbb{E} [\partial^{\ell-1}f(W_{i.}[\ell])]\bigr),\\
    & \mathcal{S}[t_{1},\cdots,t_{\ell}]\cdot E_{5,j}=\mathcal{T}_{f,j}[t_{1},\cdots,t_{\ell},\ell,\underbrace{-\ell,\cdots,-\ell}_{(j-1)\text{ times}}],\\
    & \mathcal{R}_{1}[t_{1},\cdots,t_{\ell}]\cdot E_{6}\leq \mathcal{R}_{\omega}[t_{1},\cdots,t_{\ell},\ell,\underbrace{-\ell,\cdots,-\ell}_{(k-\ell+1)\text{ times}}].
  \end{aligned}
\end{equation}
Combining \eqref{eq:annoyance5} and \eqref{eq:someannoyingthing5}, we have
\begin{align*}
  &\biggl\lvert \mathcal{S}[t_{1},\cdots,t_{\ell}]\bigl(\mathbb{E} [\partial^{\ell-1}f(W)]-\mathbb{E} [\partial^{\ell-1}f(W_{i.}[\ell])]\bigr)
  -\sum_{j=1}^{k-\ell+1}\frac{1}{j!}
  \mathcal{T}_{f,j}[t_{1},\cdots,t_{\ell},\ell,\underbrace{-\ell,\cdots,-\ell}_{(j-1)\text{ times}}] \biggr\rvert\\
  &\overset{\eqref{eq:annoyance5}}{=} \bigl\lvert \mathcal{S}[t_{1},\cdots,t_{\ell}] \bigr\rvert\cdot \bigl\lvert E_{4}-\sum_{j=1}^{k-\ell+1}E_{5,j}/j! \bigr\rvert
  \overset{\eqref{eq:type13}}{\leq}  \mathcal{R}_{1}[t_{1},\cdots,t_{\ell}]\cdot \bigl\lvert E_{4}-\sum_{j=1}^{k-\ell+1}E_{5,j}/j! \bigr\rvert\\
  &\overset{\eqref{eq:someannoyingthing5}}{\leq} \mathcal{R}_{1}[t_{1},\cdots,t_{\ell}]\cdot \lvert f \rvert_{k,\omega}E_{6}/(k-\ell+1)!
  \overset{\eqref{eq:annoyance5}}{\leq}  \frac{\lvert f \rvert_{k,\omega}}{(k-\ell+1)!}\mathcal{R}_{\omega}[t_{1},\cdots,t_{\ell},\ell,\underbrace{-\ell,\cdots,-\ell}_{(k-\ell+1)\text{ times}}].
\end{align*}
Therefore, we have established both \eqref{eq:taylorcomp1} and \eqref{eq:taylorcomp2}. Taking their difference and applying \eqref{eq:indepcomp}, we obtain \eqref{eq:taylorcomp3}.
 
\end{proof}

Equipped with the tools in \cref{thm:proptaylor}, we approximate any $\mathcal{T}$-sum $\mathcal{T}_{f,s}[t_{1},\cdots,t_{\ell}]$ by order-$j$ $\mathcal{S}$-sums ($j=\ell,\cdots,k+1$) with remainder terms being order-($k+2$) $\mathcal{R}$-sums.

% \textcolor{cyan}{[todo: check this theorem carefully]}
\begin{theorem}\label{thm:grandexpand}
Fix $k\in\mathbb{N}_{+}$. For any $\ell\in [k+1], s\in [\ell]\cup \{ 0 \}$, and $t_{1},\cdots,t_{\ell}\in \mathbb{Z}$ such that $\lvert t_{j}\rvert\leq j-1$ for any $j\in [\ell]$, there exist $Q_{\ell},\cdots,Q_{k+1}$ (which depend on $s$ and $t_{1:\ell}$ and the joint distribution of $(X_{i})_{i\in I}$) and a constant $C_{k,\ell}$ ($C_{k,\ell}\leq 4^{k-\ell+1}$) such that for any $f\in \mathcal{C}^{k,\omega }(\mathbb{R})$, we have
  \begin{equation}\label{eq:grandexpand}
    \biggl\lvert \mathcal{T}_{f,s} [t_{1},\cdots, t_{\ell}]-\sum_{j=\ell}^{k+1}Q_{j}\mathbb{E} [\partial^{j-1} f(W)]\biggr\rvert\leq C_{k,\ell}\lvert f \rvert_{k,\omega }R_{k,\omega }.
  \end{equation}
  Note that by \eqref{eq:rkalpha} $R_{k,\omega }$ is given as
  \begin{equation*}
    R_{k,\omega }= \sum_{t_{1:(k+2)}\in \mathcal{M}_{1,k+2}}\ \mathcal{R}_{\omega}[t_{1},t_{2},\cdots,t_{k+2}],
  \end{equation*}
  where 
  $$
  \mathcal{M}_{1,k+2}:=\Bigl\{t_{1:(k+2)}:~ t_{j+1}=\pm j\ \ \&\ \  t_{j}\wedge t_{j+1}< 0\ \ \forall 1\leq j\leq k+1\Bigr\}.
  $$
\end{theorem}

\begin{proof}
  If there exists an integer $2\le j\le \ell$ such that $t_{j}=0$ or there exists $j\in [\ell-1]$ such that $t_{j}$ and $t_{j+1}$ are both positive, then $\mathcal{T}_{f,s}[t_{1},\cdots,t_{\ell}]=0$ by definition and the theorem already holds by setting $Q_{j}=\cdots=Q_{k+1}=0$.

  Otherwise, we claim:
  \begin{claim}
    Let $\mathcal{T}_{f,s} [t_{1},\cdots,t_{\ell}]$ be a $\mathcal{T}$-sum. For any $j=\ell+1,\cdots,k+1$, let 
    $$
    \mathcal{E}_{\ell+1,j}:=\bigl\{ t_{(\ell+1):j}:~ \lvert t_{h+1} \rvert\leq h\ \ \&\ \  t_{h}\wedge t_{h+1}\ \ \forall  \ell\leq h\leq j-1\bigr\}.
    $$
   For all $j=\ell+1,\cdots, k+1$, $\nu\in [j]\cup \{ 0 \}$, and  $(t_{\ell+1},\cdots,t_{j})\in \mathcal{E}_{\ell,j}$, there are coefficients $a_{j,\nu,t_{(\ell+1):j}}$ (additionally depending on $s$) such that if we write 
   \begin{equation}\label{eq:defqj}
   Q_j=\sum_{t_{(l+1):j}\in \mathcal{E}_{\ell,j}}\sum_{\nu=0}^j a_{j,\nu,t_{(\ell+1):j}}\mathcal{T}_{f,\nu}[t_1,\cdots,t_{\ell},t_{\ell+1},\cdots,t_j],
   \end{equation}
then the following holds
    \begin{equation}\label{eq:tediousexpr}
      \begin{aligned}
      &\biggl\lvert \mathcal{T}_{f,s} [t_{1},\cdots, t_{\ell}]-\sum_{j=\ell}^{k+1}Q_{j}\mathbb{E} [\partial^{j-1} f(W)]\biggr\rvert\\
      \leq & 4^{k-\ell+1}\lvert f \rvert_{k,\omega }\sum_{\substack{t_{(\ell+1):(k+2)}\in \mathcal{M}_{\ell,k+1}}}%:\\t_{j}=\pm (j-1)\,\forall \ell+1\leq j\leq k+2,\\ t_{j}\wedge t_{j+1}< 0\,\forall \ell\leq j\leq k+1}}\ 
      \mathcal{R}_{\omega}\bigl[t_{1},\cdots,t_{\ell}, \cdots ,t_{k+2}\bigr],
      \end{aligned}
    \end{equation} 
    where 
    $$
    \mathcal{M}_{\ell+1,k+2}:=\Bigl\{t_{(\ell+1):(k+2)}:~ t_{j+1}=\pm j\ \ \&\ \  t_{j}\wedge t_{j+1}< 0\quad \forall \ell\leq j\leq k+1\Bigr\}.
    $$
  \end{claim}

We establish this claim by performing induction on $\ell$ with $\ell$ taking the value $k+1,k,\cdots, 1$ in turn.

  For $\ell=k+1$, by \eqref{eq:taylorcomp4} we have
  \begin{align*}
    &\Bigl\lvert \mathcal{T}_{f,s}[t_{1},\cdots,t_{k+1}]-\mathcal{S}[t_{1},\cdots,t_{k+1}]\ \mathbb{E} [\partial^{k}f(W)]\Bigr\rvert\\
    \leq & \lvert f \rvert_{k,\omega}\bigl(\mathbb{I}(t_{k+1}<0)\cdot\mathcal{R}_{\omega}[t_{1},\cdots,t_{k+1},k+1]
    + \mathbb{I}(s\ge 1)\cdot\mathcal{R}_{\omega}[t_{1},\cdots,t_{k+1},-(k+1)]\bigr).
  \end{align*}
  If there exists $j\in [k]$ such that $t_{j}$ and $t_{j+1}$ are both positive, then $\mathcal{T}_{f,s}[t_{1},\cdots,t_{k+1}]=0$ and the claim holds with all $a_{j,\nu,t_{\ell:(k+1)}}=0$. Otherwise, for all $j\le k$ either $t_{j}$ is negative or $t_{j+1}$ is negative for $j\in [k]$. If $t_{k+1}<0$, then we have 
  \begin{align*}
    & \mathbb{I}(t_{k+1}<0)\cdot\mathcal{R}_{\omega}[t_{1},\cdots,t_{k+1},k+1]
    + \mathbb{I}(s\geq 1)\cdot\mathcal{R}_{\omega}[t_{1},\cdots,t_{k+1},-(k+1)]\\
    = & \mathcal{R}_{\omega}[t_{1},\cdots,t_{k+1},k+1]
    + \mathbb{I}(s\geq 1)\cdot\mathcal{R}_{\omega}[t_{1},\cdots,t_{k+1},-(k+1)]\\
  \overset{(*)}{ \leq} & \mathcal{R}_{\omega}[0, \operatorname{sgn}(t_{2}),2\operatorname{sgn}(t_{3}),\cdots,k\cdot\operatorname{sgn}(t_{k+1}),k+1]\\
    &\ +\mathcal{R}_{\omega}[0, \operatorname{sgn}(t_{2}),2\operatorname{sgn}(t_{3}),\cdots,k\cdot\operatorname{sgn}(t_{k+1}),-(k+1)] \\
    \leq & \sum_{\substack{t_{k+2}=\pm (k+1):\\ t_{k+1}\wedge t_{k+2}< 0}}\ \mathcal{R}_{\omega}\bigl[t_{1},\cdots,t_{k+1},t_{k+2}\bigr],
  \end{align*}
  where $(*)$ is a consequence of \eqref{eq:comparecomp1} and $\operatorname{sgn}(x)=0,1,\text{ or }-1$ denotes the sign of a real number $x$.
  
  Further note that if $t_{k+1}>0$, then $\mathbb{I}(t_{k+1}<0)=0$ and we get
  \begin{align*}
    & \mathbb{I}(t_{k+1}<0)\cdot\mathcal{R}_{\omega}[t_{1},\cdots,t_{k+1},k+1]
    + \mathbb{I}(s\geq 1)\cdot\mathcal{R}_{\omega}[t_{1},\cdots,t_{k+1},-(k+1)]\\
    = & \mathbb{I}(s\geq 1)\cdot\mathcal{R}_{\omega}[t_{1},\cdots,t_{k+1},-(k+1)]\\
   \overset{(*)}{ \leq} & \mathcal{R}_{\omega}[0, \operatorname{sgn}(t_{2}),2\operatorname{sgn}(t_{3}),\cdots,k\cdot\operatorname{sgn}(t_{k+1}),-(k+1)]\\
    \leq & \sum_{\substack{t_{k+2}=\pm (k+1):\\ t_{k+1}\wedge t_{k+2}< 0}}\ \mathcal{R}_{\omega}\bigl[t_{1},\cdots,t_{k+1}, t_{k+2}\bigr],
  \end{align*} where $(*)$ is a consequence of \eqref{eq:comparecomp1}.  Thus, we have shown that
  \begin{align*}
    &\Bigl\lvert \mathcal{T}_{f,s}[t_{1},\cdots,t_{k+1}]-\mathcal{S}[t_{1},\cdots,t_{k+1}]\ \mathbb{E} [\partial^{k}f(W)]\Bigr\rvert\\
    \leq & \lvert f \rvert_{k,\omega}\bigl(\mathbb{I}(t_{k+1}<0)\cdot\mathcal{R}_{\omega}[t_{1},\cdots,t_{k+1},k+1]
    + \mathbb{I}(s\geq 1)\cdot\mathcal{R}_{\omega}[t_{1},\cdots,t_{k+1},-(k+1)]\bigr)\\
    \leq & \lvert f \rvert_{k,\omega}\sum_{\substack{t_{k+2}=\pm (k+1):\\ t_{k+1}\wedge t_{k+2}< 0}}\ \mathcal{R}_{\omega}\bigl[t_{1},\cdots,t_{k+1},t_{k+2}\bigr].
  \end{align*}
Now suppose the claim holds for $\ell+1$ and consider the case of $\ell$. By \eqref{eq:taylorcomp3} we have
\begin{align*}
  &\biggl\lvert \mathcal{T}_{f,s}[t_{1},\cdots,t_{\ell}]-\mathcal{S}[t_{1},\cdots,t_{\ell}]\ \mathbb{E} [\partial^{\ell-1}f(W)]\\*
  &\ -\mathbb{I}(s\geq 1)\cdot\sum_{j=1}^{k-\ell+1}\sum_{h=0}^{j}(-1)^{h}\frac{1}{h !(j-h)!}
  \mathcal{T}_{f,j}[t_{1},\cdots,t_{\ell},\underbrace{s - \ell,\cdots,s  - \ell}_{h\text{ times}},\underbrace{-\ell,\cdots,-\ell}_{(j-h)\text{ times}}] \\*
  &\ +\mathbb{I}(t_{\ell}<0)\sum_{j=1}^{k-\ell+1}\frac{1}{j!}
  \mathcal{T}_{f,j}[t_{1},\cdots,t_{\ell},\ell,\underbrace{-\ell,\cdots,-\ell}_{(j-1)\text{ times}}] \biggr\rvert\\
  \leq & \frac{\lvert f \rvert_{k,\omega}}{(k-\ell+1)!}\bigl(\mathbb{I}(t_{\ell}<0)\cdot\mathcal{R}_{\omega}[t_{1},\cdots,t_{\ell},\ell,\underbrace{-\ell,\cdots,-\ell}_{(k-\ell+1)\text{ times}}]\\*
  &\ + \mathbb{I}(s\geq 1)\cdot\mathcal{R}_{\omega}[t_{1},\cdots,t_{\ell},\underbrace{-\ell,\cdots,-\ell}_{(k-\ell+2)\text{ times}}]\bigr).
\end{align*} 
Note that $\mathcal{T}_{f,j}[t_{1},\cdots,t_{\ell},\underbrace{s - \ell,\cdots,s  - \ell}_{h\text{ times}},\underbrace{-\ell,\cdots,-\ell}_{(j-h)\text{ times}}]$ and $\mathcal{T}_{f,j}[t_{1},\cdots,t_{\ell},\ell,\underbrace{-\ell,\cdots,-\ell}_{(j-1)\text{ times}}]$ are  $\mathcal{T}$-sums of order at least $\ell+j$ ($j\geq 1$). Therefore, we can apply inductive hypothesis on them. In specific, the remainder term ($\mathcal{R}$-sums) in the expansion of 
$$
\mathcal{T}_{f,j}[t_{1},\cdots,t_{\ell},\underbrace{s - \ell,\cdots,s  - \ell}_{h\text{ times}},\underbrace{-\ell,\cdots,-\ell}_{(j-h)\text{ times}}]
$$
is given by this
\begin{align*}
  &\ \ 4^{k-\ell-j+1}\lvert f \rvert_{k,\omega }\!\!\!\!\!\!\!\!\!\!\!\!\!\!\!\!\!\!\!\!\sum_{\substack{t_{(\ell+j+1):(k+2)}\in \mathcal{M}_{\ell+j+1,k+2}}}\!\!\!\!\!\!\!\!\!\!\!\!\!\! \!\!\!\!\!\!\mathcal{R}_{\omega}\bigl[t_{1},\cdots,t_{\ell}, \underbrace{s - \ell,\cdots,s  - \ell}_{h\text{ times}},\underbrace{-\ell,\cdots,-\ell}_{(j-h)\text{ times}},t_{\ell+j+1},\cdots ,t_{k+2}\bigr]\\
  &\overset{\eqref{eq:comparecomp1}}{\leq} 4^{k-\ell-j+1}\lvert f \rvert_{k,\omega }\!\!\!\!\!\!\!\!\!\!\!\!\!\!\!\!\!\!\!\!\!\!\!\!\!\!\sum_{\substack{t_{(\ell+j+1):(k+2)}\in \mathcal{M}_{\ell+j+1,k+2}}}\!\!\!\!\!\!\!\!\!\!\!\!\!\!\!\!\!\!\!\!\!\!\!\!\!\! \mathcal{R}_{\omega}\bigl[t_{1},\cdots,t_{\ell}, -\ell,-(\ell+1),\cdots,-(\ell+j-1),t_{\ell+j+1},\cdots ,t_{k+2}\bigr]\\
  &\leq  4^{k-\ell-j+1}\lvert f \rvert_{k,\omega }\!\!\!\!\!\!\!\!\!\!\!\!\!\!\!\!\sum_{\substack{t_{(\ell+2):(k+2)}\in  \mathcal{M}_{\ell+2,k+2}}}\!\!\!\!\!\!\!\!\!\!\!\!\!\!\!\! \mathcal{R}_{\omega}\bigl[t_{1},\cdots,t_{\ell}, -\ell,t_{\ell+2},\cdots ,t_{k+2}\bigr]=:4^{k-\ell-j+1}\lvert f \rvert_{k,\omega}\cdot U_{1}.
\end{align*}
Similarly, the remainder term in the expansion of $\mathcal{T}_{f,j}[t_{1},\cdots,t_{\ell},\ell,\underbrace{-\ell,\cdots,-\ell}_{(j-1)\text{ times}}]$ is given by
\begin{align*}
  &4^{k-\ell-j+1}\lvert f \rvert_{k,\omega }\!\!\!\!\!\!\sum_{\substack{t_{(\ell+j+1):(k+2)}\in \mathcal{M}_{\ell+j+1,k+2}}}\!\!\!\!\!\! \mathcal{R}_{\omega}\bigl[t_{1},\cdots,t_{\ell}, \ell,\underbrace{-\ell,\cdots,-\ell}_{(j-1)\text{ times}},t_{\ell+j+1},\cdots ,t_{k+2}\bigr]\\
  \leq & 4^{k-\ell-j+1}\lvert f \rvert_{k,\omega }\!\!\!\!\sum_{\substack{t_{(\ell+2):(k+2)}\in \mathcal{M}_{\ell+2,k+2}}}\!\!\!\! \mathcal{R}_{\omega}\bigl[t_{1},\cdots,t_{\ell}, \ell,t_{\ell+2},\cdots ,t_{k+2}\bigr]=:4^{k-\ell-j+1}\lvert f \rvert_{k,\omega}\cdot U_{2}.
\end{align*}
Note that $U_{1}+\mathbb{I}(t_{\ell}<0)\cdot U_{2}$ is controlled by
\begin{align}\label{eq:notethatu1u2}
  &U_{1}+\mathbb{I}(t_{\ell}<0)\cdot U_{2}\\
  =&\sum_{t_{(\ell+2):(k+2)}\in \mathcal{M}_{\ell+2,k+2}}\ \mathcal{R}_{\omega}\bigl[t_{1},\cdots,t_{\ell}, -\ell,t_{\ell+2},\cdots ,t_{k+2}\bigr]\nonumber\\
  &\ +\mathbb{I}(t_{\ell}<0)\cdot\sum_{t_{(\ell+2):(k+2)}\in \mathcal{M}_{\ell+2,k+2}}\ \mathcal{R}_{\omega}\bigl[t_{1},\cdots,t_{\ell}, \ell,t_{\ell+2},\cdots ,t_{k+2}\bigr]\nonumber\\
  \leq &\sum_{t_{(\ell+1):(k+2)}\in \mathcal{M}_{\ell+1,k+2}}\ \mathcal{R}_{\omega}\bigl[t_{1},\cdots,t_{\ell},t_{\ell+1},\cdots ,t_{k+2}\bigr].\nonumber
\end{align}

As we mentioned above, by inductive hypothesis we have that there exist coefficients $Q_{j}$ satisfying \eqref{eq:defqj} such that 
\begin{align*}
  &\biggl\lvert \mathcal{T}_{f,s} [t_{1},\cdots, t_{\ell}]-\sum_{j=\ell}^{k+1}Q_{j}\ \mathbb{E} [\partial^{j-1} f(W)]\biggr\rvert\\
  \leq &\sum_{j=1}^{k-\ell+1}\sum_{h=0}^{j}\frac{1}{h!(j-h)!}4^{k-\ell-j+1}\lvert f \rvert_{k,\omega}\cdot U_{1}+\mathbb{I}(t_{\ell}<0)\sum_{j=1}^{k-\ell+1}\frac{1}{j!}4^{k-\ell-j+1}\lvert f \rvert_{k,\omega}\cdot U_{2}\\
  &\ +\frac{\lvert f \rvert_{k,\omega}}{(k-\ell+1)!}\bigl(\mathbb{I}(t_{\ell}<0)\cdot\mathcal{R}_{\omega}[t_{1},\cdots,t_{\ell},\ell,\underbrace{-\ell,\cdots,-\ell}_{(k-\ell+1)\text{ times}}]
  + \mathcal{R}_{\omega}[t_{1},\cdots,t_{\ell},\underbrace{-\ell,\cdots,-\ell}_{(k-\ell+2)\text{ times}}]\bigr).
\end{align*}
Noting that $\sum_{h=0}^{j}1/(h!(j-h)!)=2^{j}/j!$, we have
\begin{align*}
  &\biggl\lvert \mathcal{T}_{f,s} [t_{1},\cdots, t_{\ell}]-\sum_{j=\ell}^{k+1}Q_{j}\ \mathbb{E} [\partial^{j-1} f(W)]\biggr\rvert\\
  \leq \ \ & \sum_{j=1}^{k-\ell+1}\frac{2^{j}\cdot 4^{k-\ell-j+1}}{j!}\lvert f \rvert_{k,\omega}\cdot \bigl( U_{1}+\mathbb{I}(t_{\ell}<0)\cdot U_{2}\bigr)\\
  &\ +\frac{\lvert f \rvert_{k,\omega}}{(k-\ell+1)!}\bigl(\mathbb{I}(t_{\ell}<0)\cdot U_{2}
  + U_{1}\bigr)\\
  \leq \ \ & \bigl(1+{\textstyle\sum_{j=1}^{k-\ell+1}} 2^{2k-2\ell-j+2}\bigr)\lvert f \rvert_{k,\omega}\bigl( U_{1}+\mathbb{I}(t_{\ell}<0)\cdot U_{2}\bigr)\\
  \leq \ \ & 4^{k-\ell+1}\lvert f \rvert_{k,\omega}\bigl( U_{1}+\mathbb{I}(t_{\ell}<0)\cdot U_{2}\bigr)\\
  \overset{\eqref{eq:notethatu1u2}}{\leq} & 4^{k-\ell+1}\lvert f \rvert_{k,\omega}\sum_{t_{(\ell+1):(k+2)}\in \mathcal{M}_{\ell+1,k+2}}\ \mathcal{R}_{\omega}\bigl[t_{1},\cdots,t_{\ell},t_{\ell+1},\cdots ,t_{k+2}\bigr].
\end{align*}
Thus, we have shown \eqref{eq:tediousexpr}.

Finally, we note that for all $t_{1:\ell}\in \mathcal{M}_{1,\ell}$ and then by \eqref{eq:comparecomp1} we have
  \begin{align*}
    &\sum_{\substack{t_{(\ell+1):(k+1)}\in \mathcal{M}_{\ell+1,k+2}}}\ \mathcal{R}_{\omega}\bigl[t_{1},\cdots,t_{\ell}, \cdots ,t_{k+2}\bigr]\\
    \leq &\sum_{\substack{t_{(\ell+1):(k+2)}\in  \mathcal{M}_{\ell+1,k+2}}}\ \mathcal{R}_{\omega}\bigl[0,\operatorname{sgn}(t_{2}),2\operatorname{sgn}(t_{3}),\cdots,(\ell-1)\operatorname{sgn}(t_{\ell}), t_{\ell+1}\cdots ,t_{k+2}\bigr]\\
    \leq &\sum_{\substack{t_{1:(k+2)}\in \mathcal{M}_{1,k+2}}}\ \mathcal{R}_{\omega}[t_{1},t_{2},\cdots,t_{k+2}]=R_{k,\omega}.
  \end{align*}

\end{proof}

  We remark that if $f$ is a polynomial of degree at most $k$, then the Hölder constant $\lvert f \rvert_{k,\omega }=0$ and hence the remainder $C_{k,\ell}\lvert f \rvert_{k,\omega }R_{k,\omega}$ vanishes.

For any $\mathcal{T}$-sum, we have established the existence of expansions in \cref{thm:grandexpand}. Next we show the uniqueness of such expansions.

\begin{lemma}[Uniqueness]\label{thm:uniqueexp}
  Under the same settings as \cref{thm:grandexpand}, suppose that there exist two sets of coefficients $Q_{\ell},\cdots,Q_{k+1}$ and $Q_{\ell}',\cdots,Q_{k+1}'$ only depending on $s$ and $t_{1:\ell}$, and the joint distribution of $(X_{i})_{i\in I}$ such that for any polynomial $f$ of degree at most $\ell$, we have
  \begin{align*}
    \mathcal{T}_{f,s} [t_{1},\cdots, t_{\ell}]= & Q_{\ell}\mathbb{E} [\partial^{\ell-1} f(W)]+\cdots+Q_{k+1}\mathbb{E} [\partial^{k} f(W)]    \\
    =                                                               & Q_{\ell}'\mathbb{E} [\partial^{\ell-1} f(W)]+\cdots+Q_{k+1}'\mathbb{E} [\partial^{k} f(W)],
  \end{align*}
  Then $Q_{j}= Q_{j}'$ for any $j=\ell,\cdots, k+1$.
\end{lemma}

\begin{proof}
  We prove this lemma by contradiction.

  Let $j$ be the smallest number such that $Q_{j}\neq Q_{j}'$. Since the coefficients $Q_{\ell},\cdots,Q_{k+1}$ do not depend on $f$, we can choose $f(x)=c x^{j-1}$ such that $\partial^{j-1} f(x)= c(j-1)!\neq 0$. But $Q_{j+1}\mathbb{E} [\partial^{j}f(W)]=\cdots=Q_{k+1}\mathbb{E} [\partial^{k} f(W)]=0$, which implies $cQ_{j}=cQ_{j}'$. This is a contradiction. Therefore, $Q_{j}= Q_{j}'$ for any $j=\ell,\cdots, k+1$.
\end{proof}

\begin{proof}[Proof of \cref{THM:WFWEXPANSION}]
  Applying \cref{thm:grandexpand} with $\ell=1$, and $s=t_{1}=t_{2}=0$, we have for any $f\in \mathcal{C}^{k,\omega}(\mathbb{R})$,
  \begin{equation*}
    \mathbb{E} [Wf(W)]=\sum_{i_{1}\in I}\mathbb{E} [X_{i_{1}}f(W)]=\mathcal{T}_{f,0}[0]=\sum_{j=1}^{k+1}Q_{j}\mathbb{E} [\partial^{j-1} f(W)]+\mathcal{O}(\lvert f \rvert_{k,\omega }R_{k,\omega }),
  \end{equation*}
  for some $Q_{1},\cdots, Q_{k+1}$ that only depend on the distribution of $(X_{i})_{i\in I}$ and where $R_{k,\omega}$ is defined in \eqref{eq:rkalpha}. Suppose that $f$ is a polynomial of degree at most $k$, then we observe that $f\in \mathcal{C}^{k,\omega}(\mathbb{R})$ and  $\lvert f \rvert_{k,\omega}=0$. Thus, this implies that
  \begin{equation}\label{eq:bracketf}
    \mathcal{T}_{f,0}[0]=\mathbb{E} [Wf(W)]=\sum_{j=1}^{k+1}Q_{j}\mathbb{E} [\partial^{j-1} f(W)].
  \end{equation}
  On the other hand, for any random variable, the moments $(\mu_{j})_{j\geq 0}$ and cumulants $(\kappa_{j})_{j\geq 0}$, provided that they exist, are connected through the following relations \cite{smith1995recursive}:
  \begin{equation}\label{eq:lemcumueq}
    \mu_{n}=\sum_{j=1}^{n}\binom{n-1}{j-1}\kappa_{j}\mu_{n-j}.
  \end{equation}
  %We remark that a similar expansion can be obtained using the cumulants $(\kappa_j)$. 
  Using this we will obtain a similar expansion to \eqref{eq:bracketf} by using the cumulants $(\kappa_j)$. In this goal, we first remark that if $f(x)= x^{j}$ where $j\leq k$, then by using \eqref{eq:lemcumueq} we obtain that
  \begin{align*}
      & \mathbb{E} [Wf(W)]=\mu_{j+1}(W)
    =\sum_{h=1}^{j+1}\binom{j}{h-1}\kappa_{h}(W)\mu_{j+1-h}(W) \\
    = & \sum_{h=0}^{j}\binom{j}{h}\kappa_{h+1}(W)\mu_{j-h}(W)
    =\sum_{h=0}^{k}\frac{\kappa_{h+1}(W)}{h !}\mathbb{E} [\partial^{h} f(W)].
  \end{align*}
  Moreover, we remark that this can be generalized to arbitrary polynomials $f$ of degree $k$. Indeed, any polynomial $f$ of degree $k$ can be written as $f(x)=\sum_{j=0}^{k}a_{j}x^{j}$ for certain coefficients $(a_j)$. By the linearity of expectations, we know that
  \begin{equation*}
    \mathbb{E} [Wf(W)]=\sum_{j=0}^{k}\frac{\kappa_{j+1}(W)}{j !}\mathbb{E} [\partial^{j} f(W)].
  \end{equation*}
  Compare this to \eqref{eq:bracketf} and apply \cref{thm:uniqueexp}. We conclude that $Q_{j}=\kappa_{j}(W)/(j-1)!$ for any $j\in [k+1]$. In particular, $Q_{1}=0=\kappa_{1}(W)$.
\end{proof}

Next we upper-bound the cumulants of $W$ using $R_{k,1}$.

\begin{corollary}[Bounds for Cumulants]\label{thm:corocumubd}
  For any $k\in\mathbb{N}_{+}$, there exists a constant $C_{k}$ that only depends on $k$ such that $\bigl\lvert \kappa_{k+2}(W) \bigr\rvert\leq C_{k}R_{k,1}$.
\end{corollary}

\begin{proof}
  Let $f(x)=x^{k+1}/(k+1)!$. We remark that $f\in \Lambda_{k+1}$ where $\Lambda_{k+1}:=\{ f\in \mathcal{C}^{k,1}(\mathbb{R}):\lvert f \rvert_{k,1}\leq 1 \}$. Moreover, by using \cref{THM:WFWEXPANSION} we have
  \begin{equation*}
    \mathbb{E} [Wf(W)]=\sum_{j=1}^{k}\frac{\kappa_{j+1}(W)}{j !}\mathbb{E} [\partial^{j} f(W)]+\mathcal{O}(R_{k,1}).
  \end{equation*}
  Here the constant dropped from the big $\mathcal{O}$ analysis is controlled by $4^{k}$.
  On the other hand, by \eqref{eq:lemcumueq} we have
  \begin{align*}
    &\mathbb{E} [Wf(W)]=\frac{1}{(k+1)!}\mu_{k+2}(W)\\
    =&\sum_{j=1}^{k+1}\binom{k+1}{j}\kappa_{j+1}(W)\mu_{k+1-j}(W)\\
    = &\sum_{j=1}^{k}\frac{\kappa_{j+1}(W)}{j !}\mathbb{E} [\partial^{j} f(W)]+\frac{\kappa_{k+2}(W)}{(k+1)!}.
  \end{align*}
  Thus, there exists $C_{k}$ such that $\bigl\lvert \kappa_{k+2}(W) \bigr\rvert\leq C_{k}R_{k,1}$.
\end{proof}

Finally, we are able to prove \cref{THM:BARBOURLIKE} based on \cref{THM:WFWEXPANSION} and \cref{thm:corocumubd}.

\begin{proof}[Proof of \cref{THM:BARBOURLIKE}]
  We perform induction on $k:=\lceil p\rceil$. We start with $k=1$. In this goal, we first remark that by \cref{thm:lemsteinsol}, we have $f = \Theta h\in \mathcal{C}^{1,\omega }(\mathbb{R})$ and that $\lvert f \rvert_{1,\omega }$ is bounded by a constant. Moreover, as $f=\Theta h$ is the solution to the Stein equation \eqref{eq:stein}. By \cref{THM:WFWEXPANSION} we obtain that
  \begin{equation*}
    \mathbb{E} [h(W)]-\mathcal{N}h= \mathbb{E} [f'(W)]-\mathbb{E} [W f(W)] = \mathcal{O}(R_{1,\omega }).
  \end{equation*}
  Therefore, the desired result is established for $1$.
  Suppose that the proposition holds for $1,\cdots,k-1$, we want to prove that it will also hold for $k$. Let $f=\Theta h$, then by \cref{thm:lemsteinsol} we know that $f\in  \mathcal{C}^{k,\omega }(\mathbb{R})$ and that $\lvert f \rvert_{k,\omega}$ is bounded by some constant that only depends on $k,\omega$. Thus, by \cref{THM:WFWEXPANSION}, we have
  \begin{equation*}
    \mathbb{E} [Wf(W)]=\sum_{j=1}^{k}\frac{\kappa_{j+1}(W)}{j !}\mathbb{E} [\partial^{j} f(W)]+\mathcal{O}(R_{k,\omega }).
  \end{equation*}
  Hence we have the following expansion of the Stein equation
  \begin{align}\label{eq:lemma57start}
    &\mathbb{E} [h(W)]-\mathcal{N}h= \mathbb{E} [ f'(W)]-\mathbb{E} [W f(W)]\\
    =&-\sum_{j=2}^{k}\frac{\kappa_{j+1}(W)}{j!}\mathbb{E} [\partial^{j}f(W)]+\mathcal{O} (R_{k,\omega }) \nonumber \\
    =                               & -\sum_{j=1}^{k-1}\frac{\kappa_{j+2}(W)}{(j+1)!}\mathbb{E} [\partial^{j+1}\Theta h(W)]+\mathcal{O} (R_{k,\omega }).\nonumber
  \end{align}
  Noting that
  $\partial^{j+1}\Theta h\in \mathcal{C}^{k-j-1,\omega }(\mathbb{R})$ and $\lvert \partial^{j+1}\Theta h \rvert_{k-j-1,\omega }$ is bounded by a constant only depending on $k,\omega$, then by inductive hypothesis we obtain that 
  \begin{align}\label{eq:lemma57ind}
    &\mathbb{E} [\partial^{j+1}\Theta h(W)]-\mathcal{N}[\partial^{j+1}\Theta h]\\
    = & \sum_{(r,s_{1:r})\in \Gamma(k-j-1)}(-1)^{r}\prod_{\ell=1}^{r}\frac{\kappa _{s _{\ell}+2}(W)}{(s _{\ell}+1)!}\mathcal{N}\ \Bigl[\prod_{\ell=1}^{r}(\partial ^{s _{\ell}+1}\Theta)\comp \partial^{j+1}\Theta \ h\Bigr] \nonumber \\*
     & \ +\mathcal{O}\biggl(\sum_{\ell=1}^{k-j-1}R _{\ell,1}^{(k-j-1+\omega )/\ell}+\sum_{\ell=1}^{k-j}R _{\ell,\omega }^{(k-j-1+\omega )/(\ell+\omega -1)}\biggr),\nonumber
  \end{align}
  where we denoted $\Gamma(k-j-1):=\bigl\{ r,s_{1:r}\in \mathbb{N}_{+}: \sum_{\ell=1}^{r}s_{\ell}\leq k-j-1 \bigr\}$.

  By \cref{thm:corocumubd} and Young's inequality, we have
  \begin{equation}\label{eq:lemma57ga}
    \begin{gathered}
      \lvert \kappa_{j+2}(W) R_{\ell,\omega }^{\frac{k-j+\omega -1}{\ell+\omega -1}}\rvert\lesssim R_{j,1}R_{\ell,\omega }^{\frac{k-j+\omega -1}{\ell+\omega -1}}\leq \frac{j}{k+\omega-1 }R_{j,1}^{\frac{k+\omega -1}{j}}+\frac{k-j+\omega -1}{k+\omega -1}R_{\ell,\omega }^{\frac{k+\omega -1}{\ell+\omega -1}},\\
      \lvert \kappa_{j+2}(W) R_{\ell,1}^{\frac{k-j+\omega -1}{\ell}}\rvert\lesssim R_{j,1}R_{\ell,1 }^{\frac{k-j+\omega -1}{\ell}}\leq \frac{j}{k+\omega-1 }R_{j,1}^{\frac{k+\omega-1}{j}}+\frac{k-j+\omega -1}{k+\omega-1}R_{\ell,1 }^{\frac{k+\omega -1}{\ell}}.
    \end{gathered}
  \end{equation}
  Thus, we derive that
  \begin{align*}
                                         & \quad\mathbb{E} [h(W)]-\mathcal{N}h                                                                                                                                                                                                                                                                                                                        \\
    \overset{\eqref{eq:lemma57start}}{=} & -\sum_{j=1}^{k-1}\frac{\kappa_{j+2}(W)}{(j+1)!}\mathbb{E} [\partial^{j+1}\Theta h(W)]+\mathcal{O} (R_{k,\omega })                                                                                                                                                                                                                                     \\
    \overset{\eqref{eq:lemma57ind}}{=}   & -\sum_{j=1}^{k-1}\frac{\kappa_{j+2}(W)}{(j+1)!}\mathcal{N} [\partial^{j+1}\Theta h] +\sum_{j=1}^{k-1}\frac{\kappa_{j+2}(W)}{(j+1)!}\cdot\\*
    &\sum_{(r,s_{1:r})\in\Gamma(k-j-1)}(-1)^{r}\prod_{\ell=1}^{r}\frac{\kappa _{s _{\ell}+2}(W)}{(s _{\ell}+1)!}\mathcal{N}\ \Bigl[\prod_{\ell=1}^{r}(\partial ^{s _{\ell}+1}\Theta)\comp \partial^{j+1}\Theta \ h\Bigr] \\*
                                         & \  +\mathcal{O}\biggl(R_{k,\omega }+\sum_{j=1}^{k-1}\lvert \kappa_{j+2}(W) \rvert\sum_{\ell=1}^{k-j-1}R _{\ell,1}^{(k+\omega -j-1)/\ell}\\*
                                         &\qquad +\sum_{j=1}^{k-1}\lvert \kappa_{j+2}(W) \rvert\sum_{\ell=1}^{k-j}R _{\ell,\omega }^{(k+\omega -j-1)/(\ell+\omega -1)}\biggr)                                                                                  \\
    \overset{\eqref{eq:lemma57ga}}{=}    & -\sum_{j=1}^{k-1}\frac{\kappa_{j+2}(W)}{(j+1)!}\mathcal{N} [\partial^{j+1}\Theta h]+\sum_{j=1}^{k-1}\frac{\kappa_{j+2}(W)}{(j+1)!}\cdot\\*
    &\sum_{(r,s_{1:r})\in\Gamma(k-j-1)}(-1)^{r}\prod_{\ell=1}^{r}\frac{\kappa _{s _{\ell}+2}(W)}{(s _{\ell}+1)!}\mathcal{N}\ \Bigl[\prod_{\ell=1}^{r}(\partial ^{s _{\ell}+1}\Theta)\comp \partial^{j+1}\Theta \ h\Bigr] \\*
                                         & \  +\mathcal{O}\biggl(R_{k,\omega }+\sum_{j=1}^{k-1}R_{j,1}^{(k+\omega-1 )/j}+\sum_{j=1}^{k-1}\sum_{\ell=1}^{k-j-1}R _{\ell,1}^{(k+\omega -1)/\ell}\\*
                                         &\qquad+\sum_{j=1}^{k-1}\sum_{\ell=1}^{k-j}R _{\ell,\omega }^{(k+\omega -1)/(\ell+\omega -1)}\biggr)                                                                                                      \\
    =\ \                                 & \sum_{(r,s_{1:r})\in\Gamma(k-1)}(-1)^{r}\prod_{\ell=1}^{r}\frac{\kappa _{s _{\ell}+2}(W)}{(s _{\ell}+1)!}\mathcal{N}\ \Bigl[\prod_{\ell=1}^{r}(\partial ^{s _{\ell}+1}\Theta)\ h\Bigr]\\*
    &\ +\mathcal{O}\biggl(\sum_{\ell=1}^{k-1}R _{\ell,1}^{(k+\omega -1)/\ell}+\sum_{\ell=1}^{k}R _{\ell,\omega }^{(k+\omega -1)/(\ell+\omega -1)}\biggr).
  \end{align*}
  %  By induction, this is true for any non-negative integer $k$.
  Therefore, the desired property was established by induction.
\end{proof}

\section{Proof of Lemma~\ref{THM:EXISTENCEXI}}\label{sec:lemma2}
In \cref{THM:EXISTENCEXI}, we would like to find a random variable with a given sequence of real numbers as its cumulants. Constructing a random variable from its cumulants can be difficult in practice. However, there is a rich literature on establishing the existence of a random variable given the moment sequence. And it is well-known that the moments can be recovered from the cumulants, and vice versa. The explicit expression between moments $\mu_{n}$ and cumulants $\kappa_{n}$ is achieved by using the Bell polynomials, i.e.,
\begin{gather}
  \mu_{n} =B_{n}(\kappa_{1},\cdots,\kappa_{n})=\sum_{j=1}^{n}B_{n,j}(\kappa_{1},\cdots,\kappa_{n-j+1}),\\
  \kappa_{n} =\sum_{j=1}^{n}(-1)^{j-1}(j-1)!B_{n,j}(\mu_{1},\cdots,\mu_{n-j+1}),\label{eq:cumufrommom}
\end{gather}
where $B_{n}$ and $B_{n,j}$ are the exponential Bell polynomial defined by
\begin{equation}\label{eq:defbell}
  \begin{gathered}
    B_{n}(x_{1},\cdots,x_{n}):=\sum_{j=1}^{n}B_{n,j}(x_{1},x_{2},\cdots,x_{n-j+1}),\\
    B_{n,j}(x_{1},x_{2},\cdots,x_{n-j+1}):=\sum \frac{n!}{i_{1}!i_{2}!\cdots i_{n-j+1}!}\Bigl(\frac{x_{1}}{1!}\Bigr)^{i_{1}}\Bigl(\frac{x_{2}}{2!}\Bigr)^{i_{2}}\cdots\Bigl(\frac{x_{n-j+1}}{(n-j+1)!}\Bigr)^{i_{n-j+1}}.
  \end{gathered}
\end{equation}
The sum here is taken over all sequences $i_{1},\cdots,i_{n-j+1}$ of non-negative integers such that the following two conditions are satisfied:
\begin{gather*}
  i_{1}+i_{2}+\cdots+i_{n-j+1}=j,\\
  i_{1}+2 i_{2}+\cdots +(n-j+1)i_{n-j+1}=n.
\end{gather*}

In mathematics, the classical \textit{moment problem} is formulated as follows: Given a sequence $(\mu_{i})_{i\geq 0}$, does there exist a random variable defined on a given interval such that $\mu_{j}=\mathbb{E} [X^{j}]$ for any non-negative integer $j$? There are three essentially different types of (closed) intervals. Either two end-points are finite, one end-point is finite, or no end-points are finite, which corresponds to the \emph{Hamburger}, \emph{Hausdorff}, and \emph{Stieltjes} moment problem respectively. See \cite{akhiezer2020classical,berg1995indeterminate} or \cite{tamarkin1943problem} for a detailed discussion. For our purpose, there is no restriction on the support of random variables. Thus, the following lemma for the Hamburger moment problem is all we need.

\begin{lemma}\label{thm:hamburger}
  The Hamburger moment problem is solvable, i.e., $(\mu_{j})_{j\geq 0}$ is a sequence of moments if and only if $\mu_{0}=1$ and the corresponding Hankel kernel
  \begin{equation}\label{eq:hankel}
    H=\left(\begin{array}{cccc}
        \mu_{0} & \mu_{1} & \mu_{2} & \cdots \\
        \mu_{1} & \mu_{2} & \mu_{3} & \cdots \\
        \mu_{2} & \mu_{3} & \mu_{4} & \cdots \\
        \vdots  & \vdots  & \vdots  & \ddots
      \end{array}\right)
  \end{equation}
  is positive definite, i.e.,
  $$
    \sum_{j, k \geq 0} \mu_{j+k} c_{j} c_{k} \geq 0
  $$
  for every real sequence $(c_{j})_{j \geq 0}$ with finite support, i.e., $c_{j}=0$ except for finitely many $j$'s.
\end{lemma}

If we define the ($j$+$1$)-th upper-left determinant of a Hankel matrix by
\begin{equation}\label{eq:defhankel}
  H_{j}(x_{0},x_{1},\cdots,x_{2j}):=\left\lvert\begin{array}{cccc}
    x_{0}  & x_{1}   & \cdots & x_{j}   \\
    x_{1}  & x_{2}   & \cdots & x_{j+1} \\
    \vdots & \vdots  & \ddots & \vdots  \\
    x_{j}  & x_{j+1} & \cdots & x_{2j}
  \end{array}\right\rvert,
\end{equation}
by Sylvester's criterion in linear algebra \cite{gilbert1991positive}, the positive-definite condition above is equivalent to $H_{j}(\mu_{0},\cdots,\mu_{2j})>0$ for any $j\in \mathbb{N}_{+}$.

In order to prove \cref{THM:EXISTENCEXI}, we construct a Hankel matrix from given values of cumulants and ensure that the upper-left determinants of \eqref{eq:hankel} are all positive. Then by \cref{thm:hamburger}, there exists a random variable that has matched moments with the ones in \eqref{eq:hankel} and hence it also has the required cumulants by \eqref{eq:cumufrommom}.

For convenience, we write
\begin{equation*}
  L_{j}(x_{1},\cdots,x_{2j}):=H_{j}(1,B_{1}(x_{1}),B_{2}(x_{1},x_{2}),\cdots,B_{2j}(x_{1},\cdots,x_{2j})).
\end{equation*}
Taking $x_{1}=0$, from the definitions \eqref{eq:defbell} and \eqref{eq:defhankel}, there is an expansion
\begin{align}\label{sh}
  \begin{aligned}
  L_{j}(0,x_{2},\cdots,x_{2j})
  =&H_{j}(1,0,B_{2}(0,x_{2}),\cdots,B_{2j}(0,x_{2},\cdots,x_{2j}))\\
  =&\sum a_{t_{2},\cdots , t_{2j}}^{(j)}x_{2}^{t_{2}}\cdots x_{2j}^{t_{2j}},
  \end{aligned}
\end{align}
where the sum is taken over
\begin{equation}\label{eq:illustratemore}
\begin{gathered}
  t_{2}+t_{3}+\cdots+t_{2j}\geq j,\\
  2t_{2}+3t_{3}+\cdots+(2j)t_{2j}=j(j+1).
\end{gathered}
\end{equation}

To illustrate, set $B_{0}=1$. By \eqref{eq:defhankel} each term in $L_{j}(x_{1},\cdots,x_{2j})=H_{j}(B_{0},B_{1},\cdots,B_{2j})$ is a product of the form $\prod_{s=1}^{j+1}B_{\ell_{s}}$ (coefficients omitted) such that $\sum_{s=1}^{j+1}\ell_{s}=j(j+1)$. Note that by \eqref{eq:defbell} each monomial $M$ in $B_{\ell}(x_{1},\cdots,x_{\ell})$ satisfies that
\begin{gather*}
  \deg_{x_{1}}(M)+\deg_{x_{2}}(M)+\cdots+\deg_{x_{2j}}(M)\geq 1\quad \text{for }\ell\geq 1,\\
  \deg_{x_{1}}(M)+2\deg_{x_{2}}(M)+\cdots+(2j)\deg_{x_{2j}}(M)=\ell.
\end{gather*}
Thus, each monomial $a_{t_{1},\cdots,t_{2j}}^{(j)}x_{1}^{t_{1}}\cdots x_{2j}^{t_{2j}}$ in $L_{j}(x_{1},\cdots,x_{2j})$ satisfies that
\begin{gather*}
  \textstyle t_{1}+t_{2}+\cdots+t_{2j}\geq \bigl\lvert \{ s\in [j+1]:\ell_{s}\geq 1 \} \bigr\rvert\geq j,\\
  \textstyle t_{1}+2 t_{2}+\cdots+(2j)t_{2j}=\sum_{s=1}^{j+1}\ell_{s}=j(j+1).
\end{gather*}
Now if we take $x_{1}=0$, then $t_{1}=0$ holds for all remaining terms, and \eqref{eq:illustratemore} follows.

We further define in the following way a sequence of univariate polynomials which will be essential in our construction in \cref{THM:EXISTENCEXI}, by setting
\begin{equation*}
  P_{j}(x):=L_{j}(0,1,x,x^{2},{x^{3}},\cdots,x^{2j-2}).
\end{equation*}

Firstly, we present a lemma on the properties of $P_{j}(x)$.
\begin{lemma}\label{thm:constantterm}
  $P_{j}(x)$ is a polynomial of degree at most $j(j-1)$ with only even-degree terms and if we write
  \begin{equation*}
    P_{j}(x)=\sum_{\ell=0}^{j(j-1)/2}b_{2\ell}^{(j)}x^{2\ell},
  \end{equation*}
  we have $b_{0}^{(j)}=a_{j(j+1)/2,0,\cdots,0}^{(j)}\geq 2\ $ for any $j\geq 2,\,j\in \mathbb{N}_{+}$.
\end{lemma}

\begin{proof}
  Note that by applying \eqref{sh} we obtain that
  \begin{equation}\label{eq:pjlj}
    P_{j}(x)=L_{j}(0,1,x,\cdots,x^{2j-2})=\sum a_{t_{2},\cdots,t_{2j}}^{(j)}x^{t_{3}+2t_{4}+\cdots+(2j-2)t_{2j}},
  \end{equation}
  where the sum is taken over
  \begin{gather*}
    t_{2}+t_{3}+\cdots+t_{2j}\geq j,\\
    2t_{2}+3t_{3}+\cdots+(2j)t_{2j}=j(j+1).
  \end{gather*}
  The degree of each term in \eqref{eq:pjlj} is
  \begin{align*}
    &t_{3}+2t_{4}+\cdots+(2j-2)t_{2j}\\
    =&(2t_{2}+3t_{3}+\cdots+(2j)t_{2j})-2 (t_{2}+t_{3}+\cdots+t_{2j})\\
    = &j(j+1)-2 (t_{2}+t_{3}+\cdots+t_{2j}).
  \end{align*}
  This is even and no greater than $j(j-1)$ since $t_{2}+t_{3}+\cdots+t_{2j}\geq j$.

  Then we show the constant term $b_{0}^{(j)}\geq 2$. Consider a standard normal random variable $\xi\sim \mathcal{N}(0,1)$. Then $\kappa_{j}(\xi)=0$ for all $j\geq 1,j\neq 2$, and $\kappa_{2}(\xi)=1$, which is straightforward by checking that the moment generating function of $\xi$ is $\exp (t^{2}/2)$. By \cref{thm:hamburger}, we have
  \begin{align*}
    &b_{0}^{(j)}=P_{j}(0)=L_{j}(0,1,0,\cdots,0)\\
    =&L_{j}(\kappa_{1}(\xi),\kappa_{2}(\xi),\cdots,\kappa_{2j}(\xi))\\
    =&H_{j}(\mu_{0}(\xi),\mu_{1}(\xi),\cdots,\mu_{2j}(\xi))>0.
  \end{align*}
  Since $\mu_{2\ell}(\xi)=(2\ell-1)!!$ and $\mu_{2\ell-1}(\xi)=0$ are integers for $\ell\in\mathbb{N}_{+}$, $b_{0}^{(j)}$ is also an integer. Checking Leibniz formula of the determinant for the Hankel matrix $H_{j}$ \cite{lang2012introduction}, we observe that there is an even number of terms and that each term is odd. In specific, the determinant for the Hankel matrix is given by
    \begin{align*}
      b_{0}^{(j)}=H_{j}(\mu_{0}(\xi),\mu_{1}(\xi),\cdots,\mu_{2j}(\xi))=\sum_{\tau\in S_{j}}\operatorname{sgn}(\tau)\prod_{i=1}^{j}\mu_{\tau(i)+i-2}(\xi),
    \end{align*}
    where by abuse of notation $\operatorname {sgn}$ is the sign function of permutations in the $j$-th permutation group $S_{j}$, which returns $+1$ and $-1$ for even and odd permutations, respectively. Since $\mu_{2\ell}(\xi)=(2\ell-1)!!$ and $\mu_{2\ell-1}(\xi)=0$ for all $\ell\in\mathbb{N}_{+}$, we have
    \begin{equation*}
      \operatorname{sgn}(\tau)\prod_{i=1}^{j}\mu_{\tau(i)+i-2}(\xi)
      \begin{cases}
        \text{is odd } & \text{ if }\tau (i)+i\text{ is even }\forall i=1,\cdots,j \\
        =0             & \text{ otherwise }
      \end{cases}.
    \end{equation*}
    Noting that the number of permutations  $\tau$ that satisfies $\tau (i)+i$ is even for all $i=1,\cdots,j$ is $(j!)^{2}$, which is even when $j\geq 2$, we conclude that $b_{0}^{(j)}$ is even, and thus, it should be at least $2$.
\end{proof}

As we have explained at the beginning of this section, we would like to construct a `moment' sequence such that the corresponding Hankel kernel is positive definite. The following lemma offers one single step in the construction.
\begin{lemma}\label{thm:choosenewmu}
  Suppose there is some constant $C$ such that $\lvert \mu_{\ell}\rvert\leq C$ for $\ell=1,\cdots, 2j+1$ and $H_{j}(\mu_{0},\cdots,\mu_{2j})\geq 1$. Then there exists $C'$ only depending on $j$ and $C$ such that 
  $$H_{j+1}(\mu_{0},\cdots,\mu_{2j},\mu_{2j+1},C')\geq 1.$$
\end{lemma}

\begin{proof}
  Let $C'=(j+1) (j+1)!C^{j+2}+1$. Then by the Laplace expansion \cite{lang2012introduction} of the determinant, we have
  \begin{align*}
    H_{j+1}(\mu_{0},\cdots,\mu_{2j},\mu_{2j+1},C')= & C'H_{j}(\mu_{0},\cdots,\mu_{2j})+\sum_{\ell=0}^{j}(-1)^{j+1+\ell}\mu_{j+1+\ell}A_{j+2,\ell+1} \\
    \geq                                            & C'-(j+1)C\cdot (j+1)!C^{j+1}\geq 1,
  \end{align*}
  where $A_{j+2,\ell+1}$ is the determinant of the ($j$+$1$)$\times$ ($j$+$1$) submatrix obtained by deleting the ($j$+$2$)-th row and ($\ell$+$1$)-th column of
  \begin{equation*}
    A=\left(\begin{array}{cccc}
      \mu_{0}   & \mu_{1}   & \cdots & \mu_{j+1} \\
      \mu_{1}   & \mu_{2}   & \cdots & \mu_{j+2} \\
      \vdots    & \vdots    & \ddots & \vdots    \\
      \mu_{j+1} & \mu_{j+2} & \cdots & C'
    \end{array}\right).
  \end{equation*}
\end{proof}

Now we prove \cref{THM:EXISTENCEXI}.

\begin{proof}[Proof of \cref{THM:EXISTENCEXI}]
  The key of the proof will be to use \cref{thm:hamburger}. To do so we need to postulate an infinite sequence that will be our candidates for of potential moments and check that the conditions of \cref{thm:hamburger} hold. We remark that as we already know what we want the first $k$+$1$ cumulants to be, we already know what the candidates are for the first $k$+$1$ moments; and we only to find adequate proposal for the ($k$+$2$)-th moment onward. We will do so by iteratively using \cref{thm:choosenewmu}.

  In this goal, we remark that since by \cref{thm:constantterm} we know that $b_{0}^{(j)}\geq 2$. Therefore, we can choose a small enough constant $0<C_{p}<1$ only depending on $k=\lceil p\rceil$ such that
  \begin{align}\label{boire}
    b_{0}^{(j)}-\sum_{\ell=1}^{j(j-1)/2}\sum_{2t_{2}+2t_{3}+\cdots+2t_{2j}=j(j+1)-2\ell\atop 2t_{2}+3t_{3}+\cdots+2jt_{2j}=j(j+1)}\lvert a_{t_{2},\cdots,t_{2j}}^{(j)} \rvert C_{p}^{2\ell}\geq 1,
  \end{align}
  for any integer $j=1,\cdots, \lceil k/2\rceil$. Given an index set $I_n$, if $u_{j}^{\scalebox{0.6}{$(n)$}}= 0$ for all $j=1,\cdots, k-1$, let $\xi^{\scalebox{0.6}{$(n)$}}\sim \mathcal{N}(0,1)$ and $q_{n}\gg \lvert I_n \rvert$. Then $q_{n}$ and $\xi^{\scalebox{0.6}{$(n)$}}$ satisfy all the requirements since $\kappa_{j}(\xi^{\scalebox{0.6}{$(n)$}})=0$ for all $j\in \mathbb{N}_{+},j\neq 2$ and $\kappa_{2}(\xi^{\scalebox{0.6}{$(n)$}})=1$, which is straightforward by checking that the momemt generating function of $\xi^{\scalebox{0.6}{$(n)$}}$ is $\exp (t^{2}/2)$.

  Otherwise, let
  \begin{equation}\label{eq:chooseq}
    q_{n}:=\Bigl\lfloor\min_{1\leq j\leq k-1, u_{j}^{\scalebox{0.4}{$(n)$}}\neq 0}\bigl\{  C_{p}^{2}\lvert u_{j}^{\scalebox{0.6}{$(n)$}}\rvert^{-2/j} \bigr\}\Bigr\rfloor,
  \end{equation}
  where $\lfloor x\rfloor$ denotes the largest integer not exceeding $x$. Since by assumption, for any $j=1,\cdots, k-1$, $u_{j}^{\scalebox{0.6}{$(n)$}}\to 0$ as $n\to \infty$, then we know that there exists $N>0$ such that (i) $q_{n}\geq 1$ for any $n>N$ and (ii) $q_{n}\to\infty$ as $n\to \infty$.
  We note that by definition 
  $$
  \min_{1\leq j\leq k-1, u_{j}^{\scalebox{0.4}{$(n)$}}\neq 0}\bigl\{  C_{p}^{2}\lvert u_{j}^{\scalebox{0.6}{$(n)$}}\rvert^{-2/j} \bigr\}< q_{n}+1,
  $$ 
  which implies
  \begin{equation}\label{eq:csatisfied}
    \max_{1\leq j\leq k-1}\bigl\{ q_{n}^{j/2}\lvert u_{j}^{\scalebox{0.6}{$(n)$}} \rvert \bigr\}>C_{p}^{j}\bigl(q_{n}/(q_{n}+1)\bigr)^{j/2}>C_{p}^{p}/2^{p/2}.
  \end{equation}
  On the other hand, \eqref{eq:chooseq} also implies that $C_{p}^{2}\lvert u_{j}^{\scalebox{0.6}{$(n)$}}\rvert^{-2/j}\geq q_{n}$. Thus, $q_{n}^{j/2}\lvert u_{j}^{\scalebox{0.6}{$(n)$}}\rvert \leq C_{p}^{j}$. Now let $\widetilde{\kappa}_{j+2}:=q_{n}^{j/2}u_{j}^{\scalebox{0.6}{$(n)$}}.$ We remark that $\lvert\widetilde {\kappa}_{j+2}\rvert\le C_p^j$ and $\lvert\widetilde \kappa_{j+2}\rvert\ge C_p^p/2^{p/2}$. We write $\widetilde{\mu}_{j+2}:=B_{j+2}(0,\widetilde{\kappa}_{2},\cdots,\widetilde{\kappa}_{j+2})$ for $j=1,\cdots,k-1$. Those will be our candidates for the first $k$+$1$ moments. Moreover, if $k$ is odd, we also propose a candidate for ($k$+$2$)-th moment by setting $\widetilde{\mu}_{k+2}:=0$.

  For $j=1,\cdots,\lceil k/2\rceil$ 
  % \textcolor{blue}{[i imagine that you want here to let $j$ be smaller than $(k+2)/2$ moment is k is even and $(k+1)/2$ is $k $ is odd?} \textcolor{cyan}{no. $k/2$ when $k$ is even, $k+2$ moment cannot be trivially set to be zero when $k$ is even; $(k+1)/2$ when $k$ is odd. In general we want to control absolute $k+2$ moment, for so for odd $k$ we set $k+2$ moment to be zero and control $k+3$ moment]}, 
  by \eqref{sh} we have
  \begin{align*}
                         & H_{j}(1,0,\widetilde{\mu}_{2},\widetilde{\mu}_{3},\cdots,\widetilde{\mu}_{2j})=L_{j}(0,\widetilde{\kappa}_{2},\widetilde{\kappa}_{3},\cdots,\widetilde{\kappa}_{2j})                                                                                                                                                                                                   \\
    =                    & \sum_{2t_{2}+3t_{3}+\cdots+2jt_{2j}=j(j+1)}a_{t_{2},\cdots , t_{2j}}^{(j)}\widetilde{\kappa}_{2}^{t_{2}}\cdots \widetilde{\kappa}_{2j}^{t_{2j}}\\
    =&\sum_{\ell=0}^{j(j-1)/2}\sum_{2t_{2}+2t_{3}+\cdots+2t_{2j}=j(j+1)-2\ell\atop 2t_{2}+3t_{3}+\cdots+2jt_{2j}=j(j+1)}a_{t_{2},\cdots , t_{2j}}^{(j)}\widetilde{\kappa}_{2}^{t_{2}}\cdots \widetilde{\kappa}_{2j}^{t_{2j}} \\
    \overset{(a)}{\geq}  & b_{0}^{(j)}-\sum_{\ell=1}^{j(j-1)/2}\sum_{2t_{2}+2t_{3}+\cdots+2t_{2j}=j(j+1)-2\ell\atop 2t_{2}+3t_{3}+\cdots+2jt_{2j}=j(j+1)}\bigl\lvert a_{t_{2},\cdots , t_{2j}}^{(j)}\widetilde{\kappa}_{2}^{t_{2}}\cdots \widetilde{\kappa}_{2j}^{t_{2j}}\bigr\rvert                                                                                                              \\
    \overset{(b)}{ \geq} & b_{0}^{(j)}-\sum_{\ell=1}^{j(j-1)/2}\sum_{2t_{2}+2t_{3}+\cdots+2t_{2j}=j(j+1)-2\ell\atop 2t_{2}+3t_{3}+\cdots+2jt_{2j}=j(j+1)}\lvert a_{t_{2},\cdots,t_{2j}}^{(j)} \rvert C_{p}^{2\ell}              \quad                                                                                                                                                                 
    \overset{(c)}{\geq}  1.
  \end{align*}
  where to get $(a)$ we used the definition of $b_{0}^{(j)}$, and where to obtain $(b)$ we used the fact that  $\lvert \widetilde{\kappa}_{j+2} \rvert\leq C_{p}^{j}$, and where to get $(c)$ we used \eqref{boire}.
  Moreover, as $|\widetilde{\kappa}_{j+2}|\le C_p^j$, then we know that there exists some constant $C_{p}'$ such that $\lvert \widetilde{\mu}_{j+2}\rvert=\lvert B_{j+2}(0,\widetilde{\kappa}_{2},\cdots,\widetilde{\kappa}_{j+2}) \rvert \leq C_{p}'$ for any integer $j=1,\cdots, 2\lceil k/2\rceil-1$.
  Therefore, by \cref{thm:choosenewmu}, there exists $C_{p}''$ depending on $k=\lceil p\rceil$ and $C_{p}'$ such that
  \begin{equation*}
    H_{\lceil k/2\rceil+1}(1,0,\widetilde{\mu}_{2},\cdots,\widetilde{\mu}_{2\lceil k/2\rceil+1},C_{p}'')\geq 1.
  \end{equation*}
  Let $\widetilde{\mu}_{2\lceil k/2\rceil+2}:=C_{p}''$. Applying \cref{thm:choosenewmu} repeatedly, we get a sequence $(\widetilde{\mu}_{j})_{j\geq 1}$ such that $\widetilde{\mu}_{0}=1$ and $H_{j}(\widetilde{\mu}_{0},\widetilde{\mu}_{1},\cdots,\widetilde{\mu}_{2j})\geq 1>0$ for any $j\in \mathbb{N}_{+}$. The sequence $(\widetilde \mu_j)$ is then our candidate for the moments and we remark that they satisfy the conditions of \cref{thm:hamburger}. Therefore, by \cref{thm:hamburger}, we conclude that there exists $\xi^{\scalebox{0.6}{$(n)$}}$ such that $\mu_{j}(\xi^{\scalebox{0.6}{$(n)$}})=\widetilde{\mu}_{j}$ for any $j\in \mathbb{N}_{+}$. As the first $k$+$1$ moments uniquely define the first $k$+$1$ cumulants of a random variable we have  $\kappa_{j+2}(\xi^{\scalebox{0.6}{$(n)$}})=\widetilde{\kappa}_{j+2}=q_{n}^{j/2}u_{j}^{\scalebox{0.6}{$(n)$}}$ for all $j=1,\cdots, k-1$. Thus, the $q_{n}$ and $\xi^{\scalebox{0.6}{$(n)$}}$ that we have constructed meet the requirements of \cref{itm:match12,itm:match3more}. Moreover, \eqref{eq:csatisfied} implies that \cref{itm:boundedaway} is also satisfied. Lastly, to show \cref{itm:momentbound} we note that
    $$
    \begin{aligned}
      &\mathbb{E} [\lvert \xi^{\scalebox{0.6}{$(n)$}} \rvert^{p+2}]=\lVert \xi^{\scalebox{0.6}{$(n)$}} \rVert_{p+2}^{p+2}\overset{(*)}{\leq} \lVert \xi^{\scalebox{0.6}{$(n)$}} \rVert_{2\lceil k/2\rceil +2}^{p+2}\\
      = &\bigl(\mu_{2\lceil k/2\rceil +2}(\xi^{\scalebox{0.6}{$(n)$}})\bigr)^{(p+2)/(2\lceil k/2\rceil +2)}\\
      \leq &(C_{p}'')^{(p+2)/(2\lceil k/2\rceil +2)}.
    \end{aligned}
    $$
    Here $(*)$ is due to the fact that $k=\lceil p\rceil\geq p$.
\end{proof}

\section{Proofs of Other Results}\label{sec:lemma3}
In this section, we provide the proofs of all the other results in the main text. 

\addtocontents{toc}{\SkipTocEntry}
\subsection{Proof of Proposition \ref{THM:LEMMACONTROLBRACKET} and Theorem \ref{THM:LOCALWP2}}

For ease of notation, in this section we will drop the dependence on $n$ in our notation and write $W$, $N(\,\cdot\,)$, $\sigma$, $X_i$, $I$ and $R_{j,\omega}$ for respectively $W_n$, $N_n(\,\cdot\,)$, $\sigma_n$, $X^{\scalebox{0.6}{$(n)$}}_{i}$, $I_n$ and $R_{j,\omega,n}$.

Before we prove the bounds for $R_{k,\omega}$, we note that $R_{k,\omega}$ can be defined without assuming local dependence \textup{[LD*]}. Thus, we first aim to generalize this concept, which makes the result derived in \cref{thm:controlbracketnew} also applicable in general dependent situations. Let $(X_{i})_{i\in I}$ be a class of mean zero random variables indexed by $I$. For any graph $G$ (not necessarily the dependency graph) with the vertex set $I$ and a subset $J\subseteq I$, we define $N(J)$ to be vertex set of the neighborhood of $J$. As in \cref{sec:lemma1}, we assume $\operatorname{Var}\left(\sum_{i\in I}X_{i}\right)=1$, without loss of generality. Let $W=\sum_{i\in I}X_{i}$.

We extend the notation of $\mathcal{R}$-sums defined in \eqref{eq:defcomp3} to this general setting. Given an integer $k\in\mathbb{N}_{+}$ such that $k\geq 2$, for any $t_{1:k}\in \mathbb{Z}$ such that $\lvert t_{j}\rvert \leq j-1$ for any $j\in [k]$, let $z=\bigl\lvert\{ j:t_{j}>0 \}\bigr\rvert$. If $z\geq 1$, we write $\{ j:t_{j}>0 \}=\{ q_{1},\cdots,q_{z} \}$, where the sequence $2\leq q_{1}<\cdots<q_{z}\leq k$ is taken to be increasing. We further let $q_{0}:=1$ and $q_{z+1}:=k+1$. Then we could still define the $\mathcal{R}$-sums by
\begin{align*}
   & \mathcal{R}_{\omega}[t_{1},t_{2},\cdots,t_{k}] : =\\
   &\ \sum_{i_{1}\in N_{1}}\sum_{i_{2}\in N_{2}}\cdots\sum_{i_{k-1}\in N_{k-1}}[q_{1}-q_{0},\cdots,q_{z+1}-q_{z}]\triangleright \Bigl(\lvert X_{i_{1}}\rvert,\cdots,\lvert X_{i_{k-1}}\rvert,\bigl(\sum_{i_{k}\in N_{k}}\lvert X_{i_{k}}\rvert\bigr)^{\omega }\Bigr),
\end{align*}
where $N_{1}:=I$, and for $2\leq j\leq k$
$$
N_{j}:=\begin{cases} N (i_{1:\lvert t_{j} \rvert})=N(i_{1},\cdots,i_{\lvert t_{j}\rvert})& \text{ if }t_{j}\neq 0\\ 
  \emptyset &\text{ if }t_{j}=0
\end{cases}.
$$
Now the remainder term $R_{k,\omega}$ is defined as
\begin{align}\label{eq:rkalphanew}
  R_{k,\omega}:= & \sum_{(\ell,\eta_{1:\ell})\in C^{*}(k+2)}\sum_{i_{1}\in N_{1}'}\sum_{i_{2}\in N_{2}'}\cdots\sum_{i_{k+1}\in N_{k+1}'}\\*
  &\quad [\eta_{1},\cdots,\eta_{\ell}]\triangleright \biggl(\lvert X_{i_{1}}\rvert,\cdots,\lvert X_{i_{k+1}}\rvert,\Bigl(\sum_{i_{k+2}\in N_{k+2}'}\lvert X_{i_{k+2}}\rvert\Bigr)^{\omega }\biggr)\nonumber \\
   = &\sum_{t_{1:(k+2)}\in \mathcal{M}_{1,k+2}}\ \mathcal{R}_{\omega}[t_{1},t_{2},\cdots,t_{k+2}].\nonumber
\end{align}
where $N_{1}':=I$ and $N_{j}':=N(i_{1:(j-1)})$ for $j\geq 2$. $C^{*}(k+2)$ and $\mathcal{M}_{1,k+2}$ are given by
$$
C^{*}(k+2)=\bigl\{\ell,\eta_{1:\ell}\in\mathbb{N}_{+}: \eta_{j}\geq 2\ \forall  j\in [\ell-1], \ \sum_{j=1}^{\ell}\eta_{j}=k+2\bigr\},
$$ and 
\begin{equation*}
\mathcal{M}_{1,k+2}:=\Bigl\{t_{1:(k+2)}:~ t_{j+1}=\pm j\ \ \&\ \  t_{j}\wedge t_{j+1}< 0\ \ \forall 1\leq j\leq k+1\Bigr\}.
\end{equation*} 

Note that the expressions of $\mathcal{R}$-sums and $R_{k,\omega}$ have the same forms as those in \cref{sec:pflocalnotation}, but here we do not impose the assumption of the local dependence of $(X_{i})_{i\in I}$ anymore as $N(i_{1:q})$'s are defined directly from the graph structure we constructed on $I$. The main goal of this section is to prove the following proposition.

\begin{proposition}\label{thm:controlbracketnew}
  Fix $k\in\mathbb{N}_{+}$ such that $k\geq 2$ and real number $\omega\in (0,1]$. Let $N(J)$ be defined as above and suppose the cardinality of $N(J)$ is upper-bounded by $M$ for any $\lvert J \rvert\leq k$. Then there exists a constant $C_{k+\omega }$ only depending on $k+\omega$ such that
  \begin{equation*}
    \mathcal{R}_{\omega} [t_{1},t_{2},\cdots,t_{k}]\leq C_{k+\omega }  M^{k-2+\omega }\sum_{i\in I}\mathbb{E} [\lvert X_{i} \rvert^{k-1+\omega }].
  \end{equation*}
\end{proposition}

Before proving \cref{thm:controlbracketnew}, we need the following two lemmas. \cref{thm:lemmacontrolset} helps us change the order of summation in $\mathcal{R}_{\omega}[t_{1},\cdots,t_{k}]$ and \cref{thm:lemmayoung} is a generalized version of Young's inequality, which allows us to bound the expectations of products by sums of moments.

\begin{lemma}\label{thm:lemmacontrolset}
  Fix $k\in\mathbb{N}_{+}$ such that $k\geq 2$. For any $J\subseteq I$, let $N(J)$ be defined as above. Suppose $(i_{1},\cdots,i_{k})$ is a tuple such that $i_{1}\in I$, $i_{2}\in N(i_{1})$, $\cdots$, $i_{k}\in N(i_{1:(k-1)})$. Then for any $1\leq h\leq k$, there exists a permutation $\pi$ on $[k]$ such that $\pi (1)=h$, $i_{\pi(1)}\in I$, $i_{\pi(2)}\in N\bigl(i_{\pi(1)}\bigr)$, $\cdots$, $i_{\pi(k)}\in N\bigl(i_{\pi(1)},\cdots,i_{\pi(k-1)}\bigr)$.
\end{lemma}

\begin{proof}
  We perform induction on $k$.

  Firstly, suppose that $k=2$, then we remark that $i_{2}\in N(i_{1})\Leftrightarrow i_{1}\in N(i_{2})$. For $h=1$, we can choose $\pi$ to be the identity and the desired identity holds. For $h=2$, we let $\pi(1):=2$ and $\pi(2):=1$ and remark than once again the desired result holds.

  Suppose that the proposition is true for $2,\cdots,k-1$. We want to prove that it holds for $k$. If $h<k$, consider the tuple $(i_{1},\cdots, i_{h})$. By inductive hypothesis, there is a permutation $\widetilde{\pi}$ on $\{ 1,2,\cdots,h \}$ such that $\widetilde{\pi}(1)=h$, $i_{\widetilde{\pi}(2)}\in N\bigl(i_{\widetilde{\pi}(1)}\bigr)$, $\cdots$, $i_{\widetilde{\pi}(h)}\in N\bigl(i_{\widetilde{\pi}(1)},\cdots,i_{\widetilde{\pi}(q-1)}\bigr)$. Define
  \begin{equation*}
    \pi(j):=\left\{
    \begin{aligned}
       & \widetilde{\pi}(j) &  & \text{ if }1\leq j\leq h   \\
       & j                  &  & \text{ if }h+1\leq j\leq k
    \end{aligned}
    \right..
  \end{equation*}
  Then $\pi$ satisfies the requirements in the lemma.

  Now suppose $h=k$. $i_{k}\in N(i_{1:(k-1)})$ indicates that $i_{k}$ is a neighbor of $\{ i_{1},\cdots,i_{k-1} \}$. Then there exists $1\leq \ell\leq k-1$ such that there is an edge between $i_{k}$ and $i_{\ell}$ in the graph $G=(I,E)$. Thus, $i_{h}\in N(i_{\ell})$.

  By inductive hypothesis, there is a permutation $\widetilde{\pi}$ on $[\ell]$ such that $\widetilde{\pi}(1)=\ell$, $i_{\widetilde{\pi}(2)}\in N\bigl(i_{\widetilde{\pi}(1)}\bigr)$, $\cdots$, $i_{\widetilde{\pi}(\ell)}\in N\bigl(i_{\widetilde{\pi}(1)},\cdots,i_{\widetilde{\pi}(\ell-1)}\bigr)$.

  Define
  \begin{equation*}
    \pi(j):=\begin{cases}
        k                      & \text{ if }j=1                 \\
        \widetilde{\pi}(j-1)   & \text{ if }2\leq j\leq \ell+1  \\
        j-1                    & \text{ if }\ell+2\leq j\leq  k
    \end{cases}.
  \end{equation*}
  Then $\pi(1)=h=k$. Moreover, we have $i_{\pi({2})}=i_{\ell}\in N(i_{k})=N\bigl(i_{\pi(1)}\bigr)$, and note that for all $j=3,\cdots,\ell$ we have $i_{\pi(j+1)}=i_{\widetilde{\pi}(j)}\in N\bigl(i_{\widetilde{\pi}(1)},\cdots,i_{\widetilde{\pi}(j-1)}\bigr)=N\bigl(i_{\pi(1)},\cdots,i_{\pi(j)}\bigr)$. Finally, for all $j\ge \ell+1$ we have $i_{\pi(j+1)}=i_{j}\in N(i_{1:(j-1)})\subseteq N\bigl(i_{1},\cdots,i_{j-1},i_{k}\bigr)=N\bigl(i_{\pi(1)},\cdots,i_{\pi(j)}\bigr)$. Thus, the lemma holds for $k$ as well. By induction, the proof is complete.
\end{proof}

Also, we need a generalization of Young's inequality.

\begin{lemma}\label{thm:lemmayoung}
  Given $t\in\mathbb{N}_{+}$, let $Y_{1},\cdots,Y_{t}$ be a sequence of random variables, and real numbers $p_{1},\cdots, p_{t}>1$ satisfy that $1/p_{1}+\cdots+1/p_{t}=1$. Then for any $(\ell, \eta_{1:\ell})\in C(t):=\{ \ell,\eta_{1:\ell}\in\mathbb{N}_{+}:\sum_{j=1}^{\ell}\eta_{j}=t\}$, we have that
  \begin{equation}\label{eq:lemmayoung3}
    [\eta_{1},\cdots, \eta_{\ell}]\triangleright (\lvert Y_{1} \rvert,\cdots,\lvert Y_{t} \rvert)\leq
    \frac{1}{p_{1}}\mathbb{E} [\lvert Y_{1} \rvert^{p_{1}}]+\cdots+ \frac{1}{p_{t}}\mathbb{E} [\lvert Y_{t} \rvert^{p_{t}}].
  \end{equation}
\end{lemma}

\begin{proof}
  First, we prove
  \begin{align}
     & \mathbb{E} [\lvert Y_{1}\cdots Y_{t} \rvert]\leq \frac{1}{p_{1}}\mathbb{E} [\lvert Y_{1} \rvert^{p_{1}}]+\cdots+ \frac{1}{p_{t}}\mathbb{E} [\lvert Y_{t} \rvert^{p_{t}}],\label{eq:lemmayoung1}                            \\
     & \mathbb{E} [\lvert Y_{1} \rvert]\cdots\mathbb{E} [\lvert Y_{t} \rvert] \leq \frac{1}{p_{1}}\mathbb{E} [\lvert Y_{1} \rvert^{p_{1}}]+\cdots+ \frac{1}{p_{t}}\mathbb{E} [\lvert Y_{t} \rvert^{p_{t}}].\label{eq:lemmayoung2}
  \end{align}
  In this goal, note that Young's inequality is stated as follows: For any $a_{1},\cdots,a_{t}\geq 0$, and $p_{1},\cdots,p_{t}>1$ such that $1/p_{1}+\cdots+1/p_{t}=1$, we have
  \begin{equation*}
    a_{1}\cdots a_{t}\leq \frac{1}{p_{1}}a_{1}^{p_{1}}+\cdots+\frac{1}{p_{t}}a_{t}^{p_{t}}.
  \end{equation*}
  Thus, by Young's inequality we know that
  \begin{equation*}
    \lvert Y_{1}\cdots Y_{t}\rvert\leq \frac{1}{p_{1}}\lvert Y_{1} \rvert^{p_{1}}+\cdots+\frac{1}{p_{t}}\lvert Y_{t} \rvert^{p_{t}}.
  \end{equation*}
  Taking the expectation, we have
  \begin{equation*}
    \mathbb{E} [\lvert Y_{1}\cdots Y_{t} \rvert]\leq \frac{1}{p_{1}}\mathbb{E} [\lvert Y_{1} \rvert^{p_{1}}]+\cdots+\frac{1}{p_{t}}\mathbb{E} [\lvert Y_{t} \rvert^{p_{t}}].
  \end{equation*}
  % Thus, \eqref{eq:tech} is true for $t=2$.
  Again by Young's inequality, we obtain that
  \begin{equation*}
    \mathbb{E} [\lvert Y_{1} \rvert]\cdots \mathbb{E} [\lvert Y_{t} \rvert]\leq \frac{1}{p_{1}}\mathbb{E} [\lvert Y_{1} \rvert]^{p_{1}}+\cdots+\frac{1}{p_{t}}\mathbb{E} [\lvert Y_{t} \rvert]^{p_{t}}.
  \end{equation*}
  %   Thus, \eqref{eq:lemmayoung1} and \eqref{eq:lemmayoung2} hold for $t=2$.
  By Jensen's inequality, $\mathbb{E} [\lvert Y_{i} \rvert]^{p_{i}}\leq \mathbb{E} [\lvert Y_{i} \rvert^{p_{i}}]$ for $i\in [t]$.
  This implies that
  \begin{equation*}
    \mathbb{E} [\lvert Y_{1} \rvert]\cdots \mathbb{E} [\lvert Y_{t} \rvert]\leq \frac{1}{p_{1}}\mathbb{E} [\lvert Y_{1} \rvert^{p_{1}}]+\cdots+\frac{1}{p_{t}}\mathbb{E} [\lvert Y_{t} \rvert^{p_{t}}].
  \end{equation*}

  Finally, we prove \eqref{eq:lemmayoung3}. Let $1/q_{j}:=\sum_{i=\eta_{j-1}+1}^{\eta_{j}}1/p_{i}$ for $1\leq j\leq k$.
  \begin{align*}
                                           & [\eta_{1},\cdots,\eta_{\ell}]\triangleright (\lvert Y_{1} \rvert,\cdots,\lvert Y_{k} \rvert)                                                                                                                                                                                                                                                       \\
    =                                      & \mathbb{E} \bigl[\bigl\lvert Y_{1}\cdots Y_{\eta_{1}}\bigr\rvert\bigr]\ \mathbb{E} \bigl[\bigl\lvert Y_{\eta_{1}+1}\cdots Y_{\eta_{2}} \bigr\rvert\bigr]\ \cdots\ \mathbb{E} \bigl[\bigl\lvert Y_{\eta_{1}+\cdots+\eta_{\ell-1}+1}\cdots Y_{k}\bigr\rvert\bigr]                                                                                    \\
    \overset{\eqref{eq:lemmayoung2}}{\leq} & \frac{1}{q_{1}}\mathbb{E} \bigl[\bigl\lvert Y_{1}\cdots Y_{\eta_{1}} \bigr\rvert^{q_{1}}\bigr]+\cdots+\frac{1}{q_{k}}\mathbb{E} \bigl[\bigl\lvert Y_{\eta_{1}+\cdots+\eta_{\ell-1}+1}\cdots Y_{k} \bigr\rvert^{q_{k}}\bigr]                                                                                                                        \\
    \overset{\eqref{eq:lemmayoung1}}{\leq} & \frac{1}{p_{1}}\mathbb{E} [\lvert Y_{1} \rvert^{p_{1}}]+\cdots+\frac{1}{p_{\eta_{1}}}\mathbb{E} [\lvert Y_{\eta_{1}} \rvert^{p_{\eta_{1}}}]+\cdots\\
    &\ +\frac{1}{p_{\eta_{1}+\cdots+\eta_{\ell-1}+1}}\mathbb{E} [\lvert Y_{k+1-u _{\ell}} \rvert^{p_{\eta_{1}+\cdots+\eta_{\ell-1}+1}}]+\cdots+\frac{1}{p_{k}}\mathbb{E} [\lvert Y_{k} \rvert^{p_{k}}].
  \end{align*}
\end{proof}

Now we are ready to prove \cref{thm:controlbracketnew}.

\begin{proof}[Proof of \cref{thm:controlbracketnew}]
  By \eqref{eq:comparecomp1}, we only need to prove that the following inequality holds for any $k\in\mathbb{N}_{+}$:
  \begin{equation*}
    \mathcal{R}_{\omega} [0,\pm 1,\cdots,\pm k]\lesssim M^{k-1+\omega }\sum_{i\in I}\mathbb{E} [\lvert X_{i} \rvert^{k+\omega }].
  \end{equation*}
  Once again we write $z:=\bigl\lvert\{ j:t_{j}>0 \}\bigr\rvert$. If $z\geq 1$, we write $\{ j:t_{j}>0 \}=\{ q_{1},\cdots,q_{z} \}$, where $2\leq q_{1}<\cdots<q_{z}\leq k$ is increasing. Further let $q_{0}:=1$ and $q_{z+1}:=k+1$. 

  Noticing that
  \begin{equation*}
    \underbrace{\frac{1}{k+\omega }+\cdots+\frac{1}{k+\omega }}_{k\text{ times}}+\frac{\omega }{k+\omega }=1,
  \end{equation*}
  we apply \cref{thm:lemmayoung} and obtain that
  \begin{align}\label{eq:control001}
             & [q_{1}-q_{0},\cdots,q_{z+1}-q_{z}]\triangleright \biggl(\lvert X_{i_{1}} \rvert,\cdots,\lvert X_{i_{k}} \rvert, \biggl(\frac{1}{M}\sum_{i_{k+1}\in N(i_{1:k})}\bigl\lvert X_{i_{k+1}}\bigr\rvert\biggr)^{\omega }\biggr)                \\
    \lesssim & \mathbb{E} [\lvert X_{i_{1}}\rvert^{k+\omega }]+\ \cdots\ +\mathbb{E} [\lvert X_{i_{k}}\rvert^{k+\omega } ] +\mathbb{E} \biggl[ \biggl(\frac{1}{M}\sum_{i_{k+1}\in N(i_{1:k})}\bigl\lvert X_{i_{k+1}}\bigr\rvert\biggr)^{k+\omega }\biggr].\nonumber
  \end{align}
  
  Now by Jensen's inequality and the fact that $\bigl\lvert N(i_{1:k})\bigr\rvert\leq M$, we get that
  \begin{equation*}
    \mathbb{E} \biggl[ \biggl(\frac{1}{M}\sum_{i_{k+1}\in N(i_{1:k})}\bigl\lvert X_{i_{k+1}}\bigr\rvert\biggr)^{k+\omega }\biggr]\leq \frac{1}{M}\sum_{i_{k+1}\in N(i_{1:k})}\mathbb{E} [\lvert X_{i_{k+1}} \rvert^{k+\omega}].
  \end{equation*}
  Moreover, we remark that \begin{equation}
    \begin{aligned}
    &M^{\omega}    [q_{1}-q_{0},\cdots,q_{z+1}-q_{z}]\triangleright \biggl(\lvert X_{i_{1}} \rvert,\cdots,\lvert X_{i_{k}} \rvert, \biggl(\frac{1}{M}\sum_{i_{k+1}\in N(i_{1:k})}\bigl\lvert X_{i_{k+1}}\bigr\rvert\biggr)^{\omega }\biggr)\\
    =&[q_{1}-q_{0},\cdots,q_{z+1}-q_{z}]\triangleright \biggl(\lvert X_{i_{1}} \rvert,\cdots,\lvert X_{i_{k}} \rvert, \biggl(\sum_{i_{k+1}\in N(i_{1:k})}\bigl\lvert X_{i_{k+1}}\bigr\rvert\biggr)^{\omega }\biggr)
    \end{aligned}
  \end{equation}
  Thus, this implies that
  \begin{align}\label{eq:control005}
             &\mathcal{R}_{\omega}[0,\pm 1,\cdots,\pm k]\\
             =&\sum_{i_{1}\in I}\cdots \!\!\!\!\!\sum_{i_{k}\in N(i_{1:(k-1)})}\!\!\!\![q_{1}-q_{0},\cdots,q_{z+1}-q_{z}]\triangleright \biggl(\lvert X_{i_{1}} \rvert,\cdots,\lvert X_{i_{k}} \rvert, \biggl(\sum_{i_{k+1}\in N(i_{1:k})}\bigl\lvert X_{i_{k+1}}\bigr\rvert\biggr)^{\omega }\biggr)\nonumber
    \\
    \lesssim & M^{\omega }\sum_{i_{1}\in I}\cdots \!\!\!\!\!\!\!\sum_{i_{k}\!\in N(i_{1:(k-1)})}\biggl(\mathbb{E} [\lvert X_{i_{1}} \rvert^{k+\omega }]+\cdots+\mathbb{E} [\lvert X_{i_{k}} \rvert^{k+\omega }]+\frac{1}{M}\!\sum_{i_{k+1}\in N(i_{1:k})}\!\!\mathbb{E} [\lvert X_{i_{k+1}} \rvert^{k+\omega }]\biggr).\nonumber
  \end{align}

  Since the cardinality of $N(i_{1}),\cdots,N(i_{1:k})$ are bounded by $M$, for $j=1$ we have
  \begin{equation}\label{eq:control006}
    \sum_{i_{1}\in I}\sum_{i_{2}\in N(i_{1})}\cdots\sum_{i_{k}\in N(i_{1:(k-1)})}\mathbb{E} [\lvert X_{i_{j}} \rvert^{k+\omega }]\leq M^{k-1}\sum_{i\in I}\mathbb{E} [\lvert X_{i} \rvert^{k+\omega }].
  \end{equation}

  Now we bound
  \begin{equation*}
    \sum_{i_{1}\in I}\sum_{i_{2}\in N(i_{1})}\cdots\sum_{i_{k}\in N(i_{1:(k-1)})}\mathbb{E} [\lvert X_{i_{j}} \rvert^{k+\omega }],
  \end{equation*}
  where $j=2,\cdots,k$.

  By \cref{thm:lemmacontrolset}, for any tuple $(i_{1},\cdots,i_{k})$ in the summation, there exists a permutation $\pi$ such that $\pi(1)=j$, $i_{\pi(2)}\in N\bigl(i_{\pi(1)}\bigr)$, $\cdots$, $i_{\pi(k)}\in N\bigl(i_{\pi(1)},\cdots,i_{\pi(k-1)}\bigr)$. Let $\phi_{j}$ be a map that sends $(i_{1},\cdots,i_{k})$ to $\bigl(i_{\pi(1)},\cdots,i_{\pi(k)}\bigr)$. Then no more than $(k-1)!$ tuples are mapped to the same destination since $(i_{1},\cdots,i_{k})$ is a permutation of $\bigl(i_{\pi(1)},\cdots,i_{\pi(k)}\bigr)$ and $i_{j}$ is fixed to be $i_{\pi(1)}$. Thus, we obtain that
  \begin{align}\label{eq:control007}
         & \sum_{i_{1}\in I}\sum_{i_{2}\in N(i_{1})}\cdots\sum_{i_{k}\in N(i_{1:(k-1)})}\mathbb{E} [\lvert X_{i_{j}} \rvert^{k+\omega }]\nonumber                                                                                               \\
    \leq & (k-1)!\sum_{\pi:\pi(1)=j} \sum_{i_{\pi(1)}\in I}\sum_{i_{\pi(2)}\in N\left(i_{\pi(1)}\right)}\cdots\sum_{i_{\pi(k)}\in N\left(i_{\pi(1)},\cdots,i_{\pi(k-1)}\right)}\mathbb{E} [\lvert X_{i_{\pi(1)}} \rvert^{k+\omega }]\nonumber \\
    \leq & (k-1)!\sum_{\pi:\pi(1)=j} \sum_{i_{1}\in I}\sum_{i_{2}\in N(i_{1})}\cdots\sum_{i_{k}\in N(i_{1:(k-1)})}\mathbb{E} [\lvert X_{i_{j}} \rvert^{k+\omega }]\nonumber                                                                     \\
    \leq & ((k-1)!)^{2}M^{k-1}\sum_{i\in I}\mathbb{E} [\lvert X_{i} \rvert^{k+\omega }]\lesssim M^{k-1}\sum_{i\in I}\mathbb{E} [\lvert X_{i} \rvert^{k+\omega }].
  \end{align}

  Similarly,
  \begin{equation}\label{eq:control008}
    \sum_{i_{1}\in I}\sum_{i_{2}\in N(i_{1})}\cdots\sum_{i_{k+1}\in N(i_{1:k})}\mathbb{E} [\lvert X_{i_{k+1}} \rvert^{k+\omega }]\lesssim M^{k}\sum_{i\in I}\mathbb{E} [ \lvert X_{i} \rvert^{k+\omega }].
  \end{equation}
  Substituting \eqref{eq:control006}, \eqref{eq:control007}, and \eqref{eq:control008} into \eqref{eq:control005}, we conclude
  \begin{align*}
    \mathcal{R}_{\omega}[t_{1},t_{2},\cdots,t_{k}]
    \leq &\mathcal{R}_{\omega}\bigl[0,\operatorname{sgn}(t_{2}),2\cdot \operatorname{sgn}(t_{3}),\cdots,(k-1)\operatorname{sgn}(t_{k-1})\bigr]\\
    \lesssim & M^{k-2+\omega }\sum_{i\in I}\mathbb{E} [ \lvert X_{i} \rvert^{k-1+\omega }].
  \end{align*}
\end{proof}

\begin{proof}[Proof of \cref{THM:LEMMACONTROLBRACKET}]
  By \cref{thm:controlbracketnew}, we have
  \begin{align*}
    R_{k,\omega}\overset{\eqref{eq:rkalphanew}}{=} & \sum_{t_{1:(k+2)}\in \mathcal{M}_{1,k+2}}\ \mathcal{R}_{\omega}[t_{1},t_{2},\cdots,t_{k+2}]
    \lesssim \sum_{t_{1:(k+2)}\in \mathcal{M}_{1,k+2}}M^{k+\omega}\sum_{i\in I}\mathbb{E} [ \lvert X_{i} \rvert^{k+1+\omega }].
  \end{align*}
  Noting that $\lvert \mathcal{M}_{1,k+2} \rvert< 2^{k+1}$ \cite{heubach2009combinatorics}, we conclude that
  \begin{equation*}
    R_{k,\omega}\lesssim M^{k+\omega }\sum_{i\in I}\mathbb{E} [ \lvert X_{i} \rvert^{k+1+\omega }].
  \end{equation*}
\end{proof}

The proof of \cref{THM:LOCALWP2} relies on \cref{THM:LOCALWP} and \cref{THM:LEMMACONTROLBRACKET}. 

\begin{proof}[Proof of \cref{THM:LOCALWP2}]
  Let $k:=\lceil p\rceil$. Then $p=k+\omega -1$. Without loss of generality, we assume $\sigma_{n}= 1$.
  By \cref{THM:LEMMACONTROLBRACKET},
  \begin{equation*}
    R_{j,\omega,n}\lesssim M_{n}^{j+\omega}\sum_{i\in I_{n}}\mathbb{E} \bigl[\bigl\lvert X^{\scalebox{0.6}{$(n)$}}_{i} \bigr\rvert^{j+1+\omega }\bigr].
  \end{equation*}
  If we let $q_{1}=(k-1)/(k-j)$ and $q_{2}=(k-1)/(j-1)$, then $1/q_{1}+1/q_{2}=1$ and $(2+\omega )/q_{1}+(k+1+\omega )/q_{2}=j+1+\omega$.
  Thus,
  \begin{equation*}
    \bigl\lvert X^{\scalebox{0.6}{$(n)$}}_{i} \bigr\rvert^{j+1+\omega }=\bigl\lvert X^{\scalebox{0.6}{$(n)$}}_{i} \bigr\rvert^{(2+\omega )/q_{1}}\cdot \bigl\lvert X^{\scalebox{0.6}{$(n)$}}_{i} \bigr\rvert^{(k+1+\omega )/q_{2}}.
  \end{equation*}
  By Hölder's inequality,
  \begin{align*}
    &M_{n}^{j+\omega}\sum_{i\in I_{n}}\mathbb{E} \bigl[\bigl\lvert X^{\scalebox{0.6}{$(n)$}}_{i} \bigr\rvert^{j+1+\omega }\bigr]\\
    \leq & \Bigl(M_{n}^{1+\omega}\sum_{i\in I_{n}}\mathbb{E} \bigl[\bigl\lvert X^{\scalebox{0.6}{$(n)$}}_{i} \bigr\rvert^{2+\omega }\bigr]\Bigr)^{1/q_{1}}\Bigl(M_{n}^{k+\omega}\sum_{i\in I_{n}}\mathbb{E} \bigl[\bigl\lvert X^{\scalebox{0.6}{$(n)$}}_{i} \bigr\rvert^{k+1+\omega }\bigr]\Bigr)^{1/q_{2}}          \\
    =    & \Bigl(M_{n}^{1+\omega}\sum_{i\in I_{n}}\mathbb{E} \bigl[\bigl\lvert X^{\scalebox{0.6}{$(n)$}}_{i} \bigr\rvert^{2+\omega }\bigr]\Bigr)^{(k-j)/(k-1)}\Bigl(M_{n}^{k+\omega}\sum_{i\in I_{n}}\mathbb{E} \bigl[\bigl\lvert X^{\scalebox{0.6}{$(n)$}}_{i} \bigr\rvert^{k+1+\omega }\bigr]\Bigr)^{(j-1)/(k-1)}.
  \end{align*}

  Since
  \begin{equation*}
    \frac{\omega  (k-j)}{(k-1)(j+\omega -1)}+\frac{(j-1)(k+\omega-1 )}{(k-1)(j+\omega -1)}=1,
  \end{equation*}
  by Young's inequality (See \cref{thm:lemmayoung} for details), we get
  \begin{align*}
         & \Bigl(M_{n}^{1+\omega }\sum_{i\in I_{n}}\mathbb{E} \bigl[\bigl\lvert X^{\scalebox{0.6}{$(n)$}}_{i} \bigr\rvert^{2+\omega }\bigr]\Bigr)^{\frac{k-j}{(k-1)(j+\omega -1)}}\Bigl(M_{n}^{k+\omega }\sum_{i\in I_{n}}\mathbb{E} \bigl[\bigl\lvert X^{\scalebox{0.6}{$(n)$}}_{i} \bigr\rvert^{k+1+\omega }\bigr]\Bigr)^{\frac{j-1}{(k-1)(j+\omega -1)}}                                                    \\
    \leq & \frac{\omega  (k-j)}{(k-1)(j+\omega -1)}\Bigl(M_{n}^{1+\omega }\sum_{i\in I_{n}}\mathbb{E} \bigl[\bigl\lvert X^{\scalebox{0.6}{$(n)$}}_{i} \bigr\rvert^{2+\omega }\bigr]\Bigr)^{1/\omega }\\
    &\ +\frac{(j-1)(k+\omega -1)}{(k-1)(j+\omega -1)}\Bigl(M_{n}^{k+\omega }\sum_{i\in I_{n}}\mathbb{E} \bigl[\bigl\lvert X^{\scalebox{0.6}{$(n)$}}_{i} \bigr\rvert^{k+1+\omega }\bigr]\Bigr)^{1/(k+\omega -1)}.
  \end{align*}
  Thus, we have
  \begin{equation*}
    R_{j,\omega,n}^{1/(j+\omega -1)}\lesssim\Bigl(M_{n}^{1+\omega }\sum_{i\in I_{n}}\mathbb{E} \bigl[\bigl\lvert X^{\scalebox{0.6}{$(n)$}}_{i} \bigr\rvert^{2+\omega }\bigr]\Bigr)^{1/\omega }+\Bigl(M_{n}^{k+\omega }\sum_{i\in I_{n}}\mathbb{E} \bigl[\bigl\lvert X^{\scalebox{0.6}{$(n)$}}_{i} \bigr\rvert^{k+1+\omega }\bigr]\Bigr)^{1/(k+\omega -1)}.
  \end{equation*}
  Similarly, we derive that
  \begin{align*}
    & R_{j,1,n}^{1/j}
    \lesssim  \Bigl(M_{n}^{j+1}\sum_{i\in I_{n}}\mathbb{E} \bigl[\bigl\lvert X^{\scalebox{0.6}{$(n)$}}_{i} \bigr\rvert^{j+2}\bigr]\Bigr)^{1/j}\\
    \leq &\Bigl(M_{n}^{1+\omega }\sum_{i\in I_{n}}\mathbb{E} \bigl[\bigl\lvert X^{\scalebox{0.6}{$(n)$}}_{i} \bigr\rvert^{2+\omega }\bigr]\Bigr)^{\frac{k+\omega -j-1}{kj}}\Bigl(M_{n}^{k+\omega }\sum_{i\in I_{n}}\mathbb{E} \bigl[\bigl\lvert X^{\scalebox{0.6}{$(n)$}}_{i} \bigr\rvert^{k+1+\omega }\bigr]\Bigr)^{\frac{j-\omega }{(k-1)j}} \\
    \lesssim & \Bigl(M_{n}^{1+\omega }\sum_{i\in I_{n}}\mathbb{E} \bigl[\bigl\lvert X^{\scalebox{0.6}{$(n)$}}_{i} \bigr\rvert^{2+\omega }\bigr]\Bigr)^{1/\omega }+\Bigl(M_{n}^{k+\omega }\sum_{i\in I_{n}}\mathbb{E} \bigl[\bigl\lvert X^{\scalebox{0.6}{$(n)$}}_{i} \bigr\rvert^{k+1+\omega }\bigr]\Bigr)^{1/(k+\omega -1)}.
  \end{align*}
  Since by assumption $M_{n}^{1+\omega }\sum_{i\in I_{n}}\mathbb{E}\bigl[\bigl\lvert X^{\scalebox{0.6}{$(n)$}}_{i}\bigr\rvert^{\omega +2}\bigr]\to 0$ and $M_{n}^{k+\omega }\sum_{i\in I_{n}}\mathbb{E}\bigl[\bigl\lvert X^{\scalebox{0.6}{$(n)$}}_{i}\bigr\rvert^{p+2}\bigr]\to 0$ as $n\to \infty$, we have that $R_{j,1,n}\to 0$ as $n\to \infty$.
  Therefore, by \cref{THM:LOCALWP} and noting the fact that $p=k+\omega -1$, we conclude
  \begin{align*}
    &\mathcal{W}_{p}(\mathcal{L}(W_{n}),\mathcal{N}(0,1))\\
    \leq & C_{p}\biggl(\Bigl(M_{n}^{1+\omega }\sum_{i\in I_{n}}\mathbb{E}\bigl[\bigl\lvert X^{\scalebox{0.6}{$(n)$}}_{i}\bigr\rvert^{\omega +2}\bigr] \Bigr)^{1/\omega }+\Bigl(M_{n}^{p+1}\sum_{i\in I_{n}}\mathbb{E}\bigl[\bigl\lvert X^{\scalebox{0.6}{$(n)$}}_{i}\bigr\rvert^{p+2}\bigr] \Bigr)^{1/p}\biggr),
  \end{align*}
  where $C_{p}$ only depends on $p$.
\end{proof}

\addtocontents{toc}{\SkipTocEntry}
\subsection{Proofs of Corollaries \ref{THM:MDEPFIELD} and \ref{THM:USTATWP}}

\begin{proof}[Proof of \cref{THM:MDEPFIELD}]
  Define the graph $(T_n,E_n)$ to be such that there is an edge between $i_{1},i_{2}\in T_n$ if and only if $\lVert i_{1}-i_{2} \rVert\leq m$. From the definition of the $m$-dependent random field, $\bigl(X_{i}^{\scalebox{0.6}{$(n)$}}\bigr)_{i\in T_n}$ satisfies \textup{[LD*]}. We will therefore apply \cref{THM:LOCALWP2} to obtain the desired result. We remark that $j\in N_n\bigl(i_{1:(\lceil p\rceil +1)}\bigr)$ only if there is $\ell\in[\lceil p\rceil +1]$ such that $\|i_{\ell}-j\|\le m$, which directly implies that $\bigl|N_n\bigl(i_{1:(\lceil p\rceil +1)}\bigr)\bigr|\le (2m+1)^{d}(\lceil p\rceil+1)$. 
 
   Moreover, by Hölder's inequality, we have
   \begin{align*}
     \sum_{i\in T_{n}}\mathbb{E} \bigl[\bigl\lvert X^{\scalebox{0.6}{$(n)$}}_{i} \bigr\rvert^{\omega+2}\bigr]
     \leq &\biggl(\sum_{i\in T_{n}}\mathbb{E} \bigl[\bigl\lvert X^{\scalebox{0.6}{$(n)$}}_{i} \bigr\rvert^{2}\bigr]\biggr)^{(p-\omega)/p}\biggl(\sum_{i\in T_{n}}\mathbb{E} \bigl[\bigl\lvert X^{\scalebox{0.6}{$(n)$}}_{i} \bigr\rvert^{p+2}\bigr]\biggr)^{\omega/p}\\
     \overset{(a)}{\leq}& M^{(p-\omega)/p}\sigma^{2(p-\omega)/p}\biggl(\sum_{i\in T_{n}}\mathbb{E} \bigl[\bigl\lvert X^{\scalebox{0.6}{$(n)$}}_{i} \bigr\rvert^{p+2}\bigr]\biggr)^{\omega/p}.
   \end{align*}
   Here $(a)$ is due to the non-degeneracy condition. And this directly implies that
   \begin{align*}
     &m^{(1+\omega)d/\omega}\biggl(\sigma_n^{-(\omega+2)}\sum_{i\in T_n}\mathbb{E}\bigl[\bigl\lvert X^{\scalebox{0.6}{$(n)$}}_{i}\bigr\rvert^{\omega +2}\bigr] \biggr)^{1/\omega }\\
     \leq & m^{\frac{(1+\omega)d}{\omega}}M^{\frac{p-\omega}{p\omega}} \biggl(\sigma_{n}^{-(p+2)}\sum_{i\in T_{n}}\mathbb{E}\bigl[\bigl\lvert X^{\scalebox{0.6}{$(n)$}}_{i}\bigr\rvert^{p+2}\bigr] \biggr)^{1/p}\to 0\quad\text{ as }n\to\infty.
   \end{align*}
   Therefore, by \cref{THM:LOCALWP2}, there exists $C_{p,d}>0$ such that for $n$ large enough we have
   \begin{equation*}
     \mathcal{W}_{p}(\mathcal{L}(W_n),\mathcal{N}(0,1))\leq C_{p,d}m^{\frac{(1+\omega)d}{\omega}}M^{\frac{p-\omega}{p\omega}}\sigma_n^{-\frac{p+2}{p}}\biggl(\sum_{i\in T_n}\mathbb{E}\bigl[\bigl\lvert X^{\scalebox{0.6}{$(n)$}}_{i}\bigr\rvert^{p+2}\bigr] \biggr)^{1/p}.
   \end{equation*}
 
 Moreover, if $\bigl(X^{\scalebox{0.6}{$(n)$}}_i\bigr)$ is in addition assumed to be stationary, then by assumption there is a constant $K$ such that $\liminf_{n\to \infty}\sigma_n^{2}/\lvert T_n \rvert\geq K$. Therefore, we get that
   \begin{equation*}
     \sigma_n^{-(p+2)}\sum_{i\in T_n}\mathbb{E} \bigl[\bigl\lvert X^{\scalebox{0.6}{$(n)$}}_{i} \bigr\rvert^{p+2}\bigr]\asymp \lvert T_n \rvert^{-(p+2)/2}\cdot \lvert T_n \rvert=\lvert T_n \rvert^{-p/2}\to 0,
   \end{equation*}
   and
   \begin{equation*}
    \begin{aligned}
     \mathcal{W}_{p}(\mathcal{L}(W_n),\mathcal{N}(0,1))
     \leq &C_{p,d}m^{\frac{(1+\omega)d}{\omega}}M^{\frac{p-\omega}{p\omega}}\sigma_n^{-\frac{p+2}{p}}\biggl(\sum_{i\in T_n}\mathbb{E}\bigl[\bigl\lvert X^{\scalebox{0.6}{$(n)$}}_{i}\bigr\rvert^{p+2}\bigr] \biggr)^{1/p}\\
     =&\mathcal{O}(\lvert T_n \rvert^{-1/2}).
    \end{aligned}
   \end{equation*}
 \end{proof}
 
 \begin{proof}[Proof of \cref{THM:USTATWP}]
  Consider the index set $I_n=\{ \boldsymbol{i}=(i_{1},\cdots,i_{m}):1\leq i_{1}\leq \cdots\leq i_{m}\leq n \}\subseteq\mathbb{Z}^m$. For each $\boldsymbol{i} \in I_n$, let $\xi_{\boldsymbol{i}}:=h(X_{i_{1}},\cdots,X_{i_{m}})$. Then $W_{n}=\sigma_{n}^{-1}\sum_{\boldsymbol{i}\in I}\xi_{\boldsymbol{i}}$. Let $(I_n,E_n)$ be the graph such that there is an edge between $\boldsymbol{i},\boldsymbol{j}\in I_n$ if and only if $\{ i_{1},\cdots,i_{m} \}\cap \{ j_{1},\cdots,j_{m} \}\neq \emptyset$.

  Then we remark that the conditions \textup{[LD*]} holds. Moreover, note that $\boldsymbol{j}$ is in\\ $N_n(\boldsymbol{i}_{1:(\lceil p\rceil+1)})$ only if there is $\ell\in [\lceil p\rceil+1]$ and $k_1,k_2\in [m]$ such that $j_{k_1}=(\boldsymbol{i}_{\ell})_{k_2}$, where $(\boldsymbol{i}_{\ell})_{k_2}$ denotes the $k_{2}$-th component of the vector $\boldsymbol{i}_{\ell}$. This directly implies that the cardinality of the dependency neighborhoods are bounded by $n^{m}-\bigl(n-m(\lceil p\rceil+1)\bigr)^{m}\asymp n^{m-1}$. Moreover, the non-degeneracy condition of the U-statistic implies that $\sigma_{n}^{2}\asymp n^{2m-1}$ \cite{chen2007normal}. Applying \cref{THM:LOCALWP2}, we get that
  \begin{align*}
             & \mathcal{W}_{p}(\mathcal{L}(W_{n}),\mathcal{N}(0,1))                                                                                                                                                                  \\
    \lesssim & \Bigl(n^{m}(n^{m-1})^{1+\omega }\frac{1}{\sigma_{n}^{\omega +2}}\mathbb{E} \bigl[\bigl\lvert h(X_{1},\cdots,X_{m}) \bigr\rvert^{\omega +2}\bigr]\Bigr)^{1/\omega }                                                    \\
             & \ +\Bigl(n^{m}(n^{m-1})^{p+1}\frac{1}{\sigma_{n}^{p+2}}\mathbb{E} \bigl[\bigl\lvert h(X_{1},\cdots,X_{m}) \bigr\rvert^{p+2}\bigr]\Bigr)^{1/p}                                                                         \\
    \lesssim & n^{-1/2}\Bigl(\mathbb{E} \bigl[\bigl\lvert h(X_{1},\cdots,X_{m}) \bigr\rvert^{\omega +2}\bigr]\Bigr)^{1/\omega }+n^{-1/2}\Bigl(\mathbb{E} \bigl[\bigl\lvert h(X_{1},\cdots,X_{m}) \bigr\rvert^{p+2}\bigr]\Bigr)^{1/p} \\
    \leq     & n^{-1/2}\bigl\lVert h(X_{1},\cdots,X_{m}) \bigr\rVert _{p+2}^{(\omega+2)/\omega}+n^{-1/2}\bigl\lVert h(X_{1},\cdots,X_{m}) \bigr\rVert _{p+2}^{(p+2)/p}.
  \end{align*}
  By the moment condition, $\bigl\lVert h(X_{1},\cdots,X_{m}) \bigr\rVert _{p+2}<\infty$. Thus, we conclude $$\mathcal{W}_{p}(\mathcal{L}(W_{n}),\mathcal{N}(0,1))=\mathcal{O}(n^{-1/2}).$$
\end{proof}

% Fixed a few typos June 3rd -Tyler
\addtocontents{toc}{\SkipTocEntry}
\subsection{Proof of Theorem~\ref{gillette}}

\begin{proof}
% Firstly we remark that according to \cref{THM:LOCALWP2} we know that there is a constant $K_{p,n}$ that we can compute such that $$\Phi^c(t) -\frac{K}{t\sqrt{|I_n|}^{1-\frac{1}{p+1}}}\varphi(t(1-\frac{1}{p+1}))\le \mathcal{W}_p(\mathcal{L}(W_n),\mathcal{N}(0,1))\le \frac{K_{p,n}M_n^{1+1/\omega}}{\sqrt{|I_n|}};$$ where we shorthanded $K_{p,n}:=K_p\sqrt{|I_n|}\Big(|I_n|^{1/\omega}\max_{i\in I}\|X_i^{(n)}/\sigma_n\|_{2+\omega}^{1+\frac{2}{\omega}},~|I_n|^{1/p}\max_{i\in I}\|X_i^{(n)}/\sigma_n\|_{2+p}^{1+\frac{2}{p}}\Big).$
% Firstly using \cref{THM:LOCALWP2}, we remark that under the conditions of the theorem there is a constant $K$ that does not depend n $(X_i^{(n)})$ such that 
% \begin{align}
%     \mathcal{W}_p\big(\mathcal{L}(W_n), ~\mathcal{N}(0,1)\big)\le \frac{K}{\sqrt{|I_n|}} \max\Big(\frac{M_n^{}}{C_n^{1+\frac{2}{\omega}}}\Big)
% \end{align}
For ease of notation we write $\omega_p:=\mathcal{W}_p(\mathcal{L}(W_n),\mathcal{N}(0,1)).$ Choose $\rho\in (0,1)$. Then remark that for all $\epsilon>0$ there is $G\sim \mathcal{N}(0,1)$ such that $\|G-W_n\|_{p}\le \mathcal{W}_p(\mathcal{L}(W_n),\mathcal{N}(0,1))+\epsilon$. Therefore, by the union bound we have 
\begin{align*}
    \mathbb{P}\bigl(W_n\ge t\bigr)&=    \mathbb{P}\bigl(W_n-G+G\ge t\bigr)
    \\&\le     \mathbb{P}\bigl(W_n-G\ge (1-\rho)t\bigr)+    \mathbb{P}\bigl(G\ge \rho t\bigr)
    \\&\overset{(a)}\le \Phi^c(\rho t)+ \frac{\|W_n-G\|_p^p}{((1-\rho) t)^p}
    \\&\le  \Phi^c(\rho t)+ \frac{\bigl( \omega_p+\epsilon\bigr)^p}{((1-\rho) t)^p}
\end{align*}
where to obtain $(a)$ we have used Markov's inequality. Now as this holds for any arbitrary choice of $\epsilon>0$ we conclude that 
\begin{align*}
    \mathbb{P}\bigl(W_n\ge t\bigr)&
    \le  \Phi^c(\rho t)+ \frac{\omega_p^p}{((1-\rho) t)^p}.
\end{align*} Define the function $g_t:x\mapsto (1-x)^{p+1} e^{-\frac{(xt)^2}{2}},$ then we can remark that $g_t:[0,1]\rightarrow[0,1]$ is a bijection. Choose $\rho^*_t:=g_t^{-1}\Big(\frac{\sqrt{2\pi}p\omega_p^p}{t^{p+1}}\Big)$.  
Moreover, we obtain that  \begin{align}\label{ts13}
    \mathbb{P}\bigl(W_n\ge t\bigr)&
    \le  \Phi^c(\rho^*_t t)+ \varphi(\rho^*_t t)(1-\rho^*_t)t\frac{\omega_p^p}{t^{p+1}(1-\rho_{t}^{*})^{p+1} \varphi(\rho^*_t t)}
    \\&\overset{(a)}{\le }\Phi^c(t)+ (1-\rho^*_t)t\varphi(\rho^*_t t)\big(1+\frac{1}{p}\big)\nonumber
    \\&\le \Phi^c(t)+ \frac{p^{\frac{1}{p+1}}\omega_p^{1-\frac{1}{p+1}}}{t}\varphi\Bigl(\rho^*_t t\Bigl(1-\frac{1}{p+1}\Bigr)\Bigr)\Bigl(1+\frac{1}{p}\Bigr)\nonumber
    % \\&\le  \Phi^c(t)+{p^{\frac{1}{p+1}}\omega_p^{1-\frac{1}{p}}}\varphi(t)^{-1/(p+1)}\varphi(\rho t)+\frac{\omega_p^p}{(p\omega_p^p)^{p/(p+1)}}\varphi(t)^{1-\frac{1}{p}}
\end{align}where to obtain (a) we used the fact that $\Phi^c(\rho^*_t t)\le \Phi^c(t)+(1-\rho^*_t)t\sup_{x\in [\rho^*_t t,t]}\varphi(x).$ 
\\ Suppose that $t\ge 1$  and satisfies $1-\frac{\sqrt{2\beta\log t}}{t}\le 1$. Define $$x:= \frac{\sqrt{2\beta \log t }}{t},$$ we notice that $x\in [0,1]$.  We remark that if     $$\omega_p\le (\sqrt{2\pi }p)^{\frac{1}{p+1}}\Big(1-\frac{\sqrt{2\beta\log t }}{t}\Big)t^{1-\frac{\beta}{p+1}}.$$ then we have $g_t^{-1}(x)\ge \frac{\sqrt{2\pi}p\omega_p^p}{t^{p+1}}$. 
Therefore as $g_t^{-1}(\cdot)$ is a decreasing function we have that $x\le \rho^*_t$ which implies that $$\mathbb{P}(W_n\ge t)\le \Phi^c(t)+ \frac{1}{t^{1+\beta\big(1-\frac{1}{p+1}\big)}}{p^{\frac{1}{p+1}}\omega_p^{1-\frac{1}{p+1}}}\Bigl(1+\frac{1}{p}\Bigr).$$
% \\Now as $x\rightarrow\varphi(xt)$ is decreasing  and as $\rho\ge \max\Big(1-\frac{\big(p\omega^{p}_p\big)^{\frac{1}{p+1}}}{t\varphi(t)^{\frac{1}{p+1}} } ,~\Big)$we obtain that:
% \begin{align*}
%     \mathbb{P}\bigl(W_n\ge t\bigr)&
%   \overset{}{\le }\Phi^c(t)+ \big[(1-\rho)^{p+1}\varphi(\rho t)\big]^{1/(p+1)}\varphi(\rho t)^{1-\frac{1}{p}}
%   \\&\le \Phi^c(t)+ \frac{p^{\frac{1}{p+1}}\omega_p^{1-\frac{1}{p+1}}}{t}\varphi\Big((1-\frac{\big(p\omega^{p}_p\big)^{\frac{1}{p+1}}}{t\varphi(t)^{\frac{1}{p+1}} })t\Big)^{1-\frac{1}{p+1}}
%     % \\&\le  \Phi^c(t)+{p^{\frac{1}{p+1}}\omega_p^{1-\frac{1}{p}}}\varphi(t)^{-1/(p+1)}\varphi(\rho t)+\frac{\omega_p^p}{(p\omega_p^p)^{p/(p+1)}}\varphi(t)^{1-\frac{1}{p}}
% \end{align*}

\noindent Moreover, similarly we can remark that 
\begin{align*}
    \mathbb P(G\ge (1+\rho) t)&\le \mathbb P(W_n\ge t )+\mathbb P(G-W_n\ge \rho t)
  \\  &\le  \mathbb P(W_n\ge t )+\frac{\big(\omega_p+\epsilon)^p}{\rho^pt^p}
    \end{align*} Therefore, as this holds for any arbitrary $\epsilon>0$ we obtain that
    \begin{align*}
    \Phi^c((1+\rho) t)&\le \mathbb P(W_n\ge t )+\frac{\omega_p^p}{\rho^pt^p}.
    \end{align*} 
    Moreover, we can definite $\widetilde g_{t}:x\mapsto e^{-(1+x)^2t^2}x^{p+1}$ then choose $\widetilde \rho_t^*:=\widetilde g_t^{-1}\Big(\frac{\sqrt{2\pi}p\omega_p^p}{t^{p+1}}\Big)$. We similarly obtain that 
    \begin{align}\label{ts132}
       \mathbb P(W_n\ge t)%&\ge \Phi^c(t)-\widetilde \rho_t^*t\varphi((1+\widetilde \rho_t^*)t)\big(1+\frac{1}{p}\big)
        \ge \Phi^c(t)-\frac{p^{\frac{1}{p+1}}\omega_p^{1-\frac{1}{p+1}}}{t}\varphi\Bigl(t\Bigl({1-\frac{1}{p+1}}\Bigr)\Bigr)\Bigl(1+\frac{1}{p}\Bigr).\nonumber
    \end{align}
%     By combining \cref{ts13,ts132}
% we obtain that    \begin{align}
%       \Big|  P(W_n\ge t)- \Phi^c(t)\Big|\le  \omega_p^{1-\frac{1}{p}}\frac{p^{\frac{1}{p+1}}\big(1+\frac{1}{p}\big)}{t}\begin{cases}
%          \varphi\Big(t\rho^*_t (1-\frac{1}{p+1})\Big)\\\varphi\Big(t({1-\frac{1}{p+1}})\Big)
%       %  \\&\ge \Phi^c(t)-\frac{p^{\frac{1}{p+1}}\omega_p^{1-\frac{1}{p}}}{t}\varphi\Big(t(1+\widetilde \rho_t^*)({1-\frac{1}{p+1}})\Big)\big(1+\frac{1}{p}\big)
%       \end{cases}
%     \end{align}Thefore using \textcolor{red}{cite }we know that there exists a constant $K_p$ such that $\omega_p\le \frac{K_p}{\sqrt{n}}$ therefore we obtain that there s a constant $C>0$ such that 
%      \begin{align}
%       \Big|  P(W_n\ge t)- \Phi^c(t)\Big|\le  \frac{C}{t\sqrt{n}^{1-\frac{1}{p+1}}}\begin{cases}
%          \varphi\Big(t\rho^*_t (1-\frac{1}{p+1})\Big)\\\varphi\Big(t(1+\widetilde \rho_t^*)({1-\frac{1}{p+1}})\Big)
%       %  \\&\ge \Phi^c(t)-\frac{p^{\frac{1}{p+1}}\omega_p^{1-\frac{1}{p}}}{t}\varphi\Big(t(1+\widetilde \rho_t^*)({1-\frac{1}{p+1}})\Big)\big(1+\frac{1}{p}\big)
%       \end{cases}
%     \end{align}N
\end{proof}

% \begin{supplement}
% \stitle{Title of Supplement A.}
% \sdescription{Short description of Supplement A.}
% \end{supplement}
% \begin{supplement}
% \stitle{Title of Supplement B.}
% \sdescription{Short description of Supplement B.}
% \end{supplement}

%%%%%%%%%%%%%%%%%%%%%%%%%%%%%%%%%%%%%%%%%%%%%%%%%%%%%%%%%%%%%%%%%%%
%%                                                               %%
%% Use the two commands below for producing your bibliography    %%
%% with bibtex, then comment again the commands and include the  %%
%% content of the .bbl file in this file below the commands.     %%
%%                                                               %%
%%%%%%%%%%%%%%%%%%%%%%%%%%%%%%%%%%%%%%%%%%%%%%%%%%%%%%%%%%%%%%%%%%%

\bibliographystyle{amsalpha}
\addtocontents{toc}{\SkipTocEntry}
\bibliography{biblio}

\end{document}